\documentclass[thmsd]{article}
\usepackage{amssymb}
\usepackage{latexsym} %% added by arXiv admin for \Box defn. 11May2000
\usepackage{sw20lart}
\input tcilatex
\begin{document}
\title{On the Chow ring of the classifying stack of $PGL_{3,\mathbf{C}}$}
\author{Gabriele Vezzosi \\
%EndAName
\emph{Dipartimento di Matematica,}\\
\emph{\ Universit\`a di Bologna, }\\
\emph{\ Piazza di Porta S. Donato 5,}\\
\emph{\ 40127 Bologna (Italy)}\\
\textrm{vezzosi@dm.unibo.it}}
\maketitle
\tableofcontents

\QSubDoc{Include PGL_31}{%%%%%%%%%%%%%%%%%%%%%%%%%% Start PGL_31.tex %%%%%%%%%%%%%%%%%%%%%%%%%%%%%%

\LaTeXparent{PGL_3}
\ChildStyles{amsfonts}
\ChildDefaults{chapter:1,page:1}

\section{Introduction\ }

Equivariant intersection theory is similar to Borel's equivariant
cohomology. The common basic idea is simple. Let $X$ be an algebraic scheme
over a field $k$ and let $G$ be an algebraic group acting on $X$. Since
invariant cycles are often too few to get a full-fledged intersection theory
(e.g. to have a ring structure in smooth cases) we decide to enlarge this
class to include invariant cycles not only on $X$ but on $X\times V$ where $%
V $ is any linear representation of $G$. If $k=\mathbf{C}$, equivariant
cohomology. can be defined along these lines and this definition agrees with
the usual one given using the classifying space of $G$.

In particular, we get a non trivial equivariant intersection theory $%
A_G^{*}=A_G^{*}(pt)$ on $pt=Speck$ which can be interpreted naturally as an
intersection theory on the classifying stack of the group in the same way as
equivariant cohomology of a point is naturally viewed as cohomology of the
classifying space of the group.

Equivariant intersection theory (in the sense sketched above) was first
defined by Totaro in \cite{To2} for $X=Speck$ and then extended to general $%
X $ by Edidin and Graham in \cite{EG1}. Totaro himself (\cite{To2}) and
Pandharipande (\cite{Pa1}, \cite{Pa2}) computed $A_G^{*}$ in many
interesting cases, for example $G=GL_n,$ $SL_n$ (these two cases are
trivial), $O(n),$ $SO(2n+1)$ and $SO(4)$.

Moreover, Totaro (\cite{To1}, \cite{To2}) was able to define a remarkable
refining of the classical cycle map from the Chow ring to the cohomology
ring. In particular, he proved that for any complex algebraic group $G$, the
equivariant version of the cycle class map $\mathrm{cl}_G:A_G^{*}\rightarrow
H^{*}\left( \mathrm{B}G,\mathbf{Z}\right) $, factors as 
\[
A_G^{*}\stackrel{\widetilde{\mathrm{cl}}_G}{\longrightarrow }MU^{*}\left( 
\mathrm{B}G\right) \otimes _{MU^{*}}\mathbf{Z}\longrightarrow H^{*}\left( 
\mathrm{B}G,\mathbf{Z}\right) 
\]
where $MU^{*}\left( \mathrm{B}G\right) $ is the complex cobordism of the
classifying space of $G$, 
\[
MU^{*}\left( pt\right) \equiv MU^{*}=\mathbf{Z}\left[ x_1,x_2,x_3,\ldots
\right] 
\]
where $\deg x_i=-2i$ and $\mathbf{Z}$ is viewed as a $MU^{*}$-module via the
map sending each $x_i$ to zero. He conjectured that, if $MU^{*}\left( 
\mathrm{B}G\right) \otimes _{MU^{*}}\mathbf{Z}$ is concentrated in even
degrees, then $\widetilde{\mathrm{cl}}_G$ is an isomorphism.

The case $G=PGL_n$ is of particular interest. One reason is its connection
with Brauer-Severi varieties, whose Chow groups are quite mysterious (see 
\cite{Ka1} and \cite{Ka2} for some results on codimension $2$ cycles). Also,
many parameter spaces of interest are quotients of free actions of $PGL_n$,
so the calculation of $A_{PGL_n}^{*}$ would be a necessary first step to
determine the Chow ring of some of these spaces.

Unfortunately, the ring $A_{PGL_n}^{*}$ for general $n$ seems extremely
difficult to compute. It is a general principle that among all families of
classical groups the series $PGL_n$ is often the hardest to study. Thus, for
example, while the cohomology and the complex cobordism ring of most
classical groups have been determined, very little is known about the
torsion part in the cohomology of the classifying space of $PGL_n$ for $%
n\geq 4$. Of course, given how much harder than cohomology the Chow ring
usually is, this is not encouraging. On the other hand, the cohomology with $%
\mathbf{Z}/3$ coefficients of the classifying space of $PGL_3$, as well as
its Brown-Peterson cohomology (relative to the prime $3$) have been computed
by Kono, Mimura and Shimada (\cite{KMS}) and by Kono and Yagita (\cite{KY}).

The ring $A_{PGL_2}^{*}$ was first computed by Pandharipande (\cite{Pa1})
through the isomorphism $PGL_2\simeq SO\left( 3\right) $. Pandharipande's
method does not seem to extend to $PGL_3$.

In this paper we study $A_{PGL_3}^{*}$. Our approach is completely
different. The idea is that the adjoint representation $\mathrm{sl}_n$ of $%
PGL_n$ can be stratified, using Jordan canonical form, in such a way that
the equivariant Chow ring of each stratum is amenable to study. This
determines completely $A_{PGL_2}^{*}$ (\cite{Ve}) and works fairly well for $%
n=3$ yielding generators of $A_{PGL_3}^{*}$. In principle this method could
give generators for $A_{PGL_n}^{*}$ for any $n$, but the calculations become
extremely involved as $n$ grows. Moreover, as usual, the stratification
method is not very good for finding the relations. In the case $n=3$, using
also a recent general result by Totaro (Th. \ref{gottliebtotaro}), we find
some of the relations in section 5, but unfortunately we are not able to
prove that our relations are sufficient.

We also prove some properties of the cycle map and of Totaro's refined cycle
map. In particular, we are able to prove that $A_{PGL_3}^{*}$, unlike $%
A_{PGL_2}^{*}$, is not generated by Chern classes of representations, a
result conjectured by Totaro in \cite{To2}. We have two proofs of this fact,
one (Th. \ref{result1}), relying on results of \cite{To2}, \cite{KMS}, \cite
{KY}, carries more informations on the cycle and refined cycle maps while
the other (Appendix) is self-contained not depending on cohomological
arguments.

Most of the results in this paper constitute the core of \cite{Ve}.

\smallskip\ 

Now we state the main results of this work in more detail.

It is already clear ''rationally'', that Chern classes of the adjoint
representation alone do not generate $A_{PGL_3}^{*}$. So, if $E$ is the
standard representation of $GL_3$, we also consider $Sym^3E$, the $PGL_3$%
-representation defined by 
\[
\left[ g\right] \cdot \left( v_1\cdot v_2\cdot v_3\right) \doteq \det
g^{-1}\left( gv_1\cdot gv_2\cdot gv_3\right) \text{.} 
\]
We prove the following (Theorem \ref{generators}):

\begin{theorem}
\label{duos}There exist elements $\rho $ and $\chi $ with\emph{\ }$\deg \rho
=4$\emph{\ }and $\deg \chi =6$, such that $A_{PGL_3}^{*}$ is generated by 
\[
\left\{ \lambda \doteq 2c_2\left( \mathrm{sl}_3\right) -c_2\left(
Sym^3E\right) ,c_3\left( Sym^3E\right) ,\rho ,\chi ,c_6\left( \mathrm{sl}%
_3\right) ,c_8\left( \mathrm{sl}_3\right) \right\} \text{.}
\]
\end{theorem}

The question of determining all the relations between this generators is
hard. In this direction, we can prove the following (Th. \ref{relations}):

\begin{proposition}
\label{tres}The generators above satisfy 
\[
3\rho =3\chi =3c_8\left( \mathrm{sl}_3\right) =0
\]
\[
\rho ^2=c_8\left( \mathrm{sl}_3\right) 
\]
\[
3\left( 27c_6\left( \mathrm{sl}_3\right) -c_3\left( Sym^3E\right)
^2-4\lambda ^3\right) =0\text{.}
\]
\end{proposition}

Moreover, if $R^{*}$ denotes the ring 
\[
\frac{\mathbf{Z}\left[ \lambda ,c_3\left( Sym^3E\right) ,\rho ,\chi
,c_6\left( \mathrm{sl}_3\right) ,c_8\left( \mathrm{sl}_3\right) \right] }{%
\frak{R}} 
\]
where $\frak{R}$ is the ideal generated by the relations in Prop. \ref{tres}
and\emph{\ }$\deg \rho =4$\emph{, }$\deg \chi =6$, we have (Theorem \ref
{riassunto})

\begin{theorem}
The composition 
\[
R^{*}\longrightarrow A_{PGL_3}^{*}\stackrel{\widetilde{\mathrm{cl}}}{%
\longrightarrow }MU^{*}\left( \mathrm{B}PGL_3\right) \otimes _{MU^{*}}%
\mathbf{Z}
\]
is surjective and its kernel is $3$-torsion.
\end{theorem}

Note that this also proves that the $R^{*}\left[ \frac 13\right] \simeq
A_{PGL_3}^{*}\left[ \frac 13\right] $. We also prove that while $\rho $ is
nonzero in cohomology, $\chi $ is zero in cohomology. Thus, by Remark \ref
{totaro}, we also have $\widetilde{\mathrm{cl}}\left( \chi \right) =0$. Note
that if one was able to prove that $\chi \neq 0$ then Totaro's conjecture
would be false. However, despite many efforts, we still do not know whether $%
\chi $ is zero or not.

By a result of Kono and Yagita (\cite{KY}), Totaro's conjecture predicts
that $\widetilde{\mathrm{cl}}$ is actually an isomorphism. We are able to
show that the generator $\rho $ of Theorem \ref{duos} is not in the Chern
subring\footnote{%
I.e. in the subring generated by Chern classes of representations.} of $%
A_{PGL_3}^{*}$, thus proving the following consequence of Totaro's
conjecture (Theorem \ref{result1}):

\begin{theorem}
$A_{PGL_3}^{*}$ is not generated by Chern classes.
\end{theorem}

This same result is proved in the Appendix without using cohomology
computations.

\bigskip\ 

\textbf{Conventions and notations. }The word ''scheme'' will most of the
time mean ''algebraic scheme over a field $k$''. In section 1, where we try
to give some of the results in greater generality, we will allow a different
base scheme $S$ and the finiteness conditions needed will be properly
specified.

We freely use the functorial point of view for schemes and group schemes
(e.g. \cite{DG}) to be able to express maps, actions etc. as sending
''elements to elements''.

If $s$ is a section of a vector bundle, we denote by $Z\left( s\right) $ its
zero scheme.

Algebraic groups over a field $k$ will always be linear. If $G$ is an
algebraic group over a field $k$, $T_G$ (or simply $T$ if no confusion is
possible) denotes a maximal torus of $G$ and $\widehat{T_G}$ its character
group.

If $\varphi :G\rightarrow H$ is a morphism of algebraic groups over a field $%
k$ and $V$ is a representation of $H$, we denote by $V_{\left( \varphi
\right) }$ or $V_{\left( G\right) }$ the obvious associated $G$%
-representation.

If $E$ denotes the standard $GL_3$-representation, $Sym^3E$ becomes a $PGL_3$%
-representation via 
\[
\left[ g\right] \cdot \left( v_1\cdot v_2\cdot v_3\right) \doteq \det
g^{-1}\left( gv_1\cdot gv_2\cdot gv_3\right) \text{.} 
\]

\smallskip\ 

\textbf{Acknowledgments.} I wish to thank my thesis advisor Angelo Vistoli
for his friendship, patience and constant attention to this work.

I also wish to thank Burt Totaro for generous advice and for many
illuminating discussions we had in Cambridge. Among other things, he also
explained to me the argument in Remark \ref{totaro} and allowed me to
include his still unpublished result Theorem \ref{gottliebtotaro}. I hope
this paper shows all its debts to his deep and original work.

%%%%%%%%%%%%%%%%%%%%%%%%%%% End PGL_31.tex %%%%%%%%%%%%%%%%%%%%%%%%%%%%%%%
}

\QSubDoc{Include PGL_32}{%%%%%%%%%%%%%%%%%%%%%%%%%% Start PGL_32.tex %%%%%%%%%%%%%%%%%%%%%%%%%%%%%%

\LaTeXparent{PGL_3}
\ChildStyles{amssymb} 
\ChildDefaults{chapter:2,page:1}

\section{ Basic notations and results}

In this section we mainly fix notations and collect some miscellaneous
results on equivariant Chow groups we will need in the sequel; most of them
(with one possible exception) are elementary or well known but we simply
could not find proper references in the literature. For intersection theory
the standard refernce is \cite{Fu} while for equivariant intersection theory
we refer to \cite{EG1} and \cite{To2}.

\subsection{Equivariant intersection theory and Totaro's refined cycle map\ }

If $G$ be an algebraic group over a field $k$ and $X$ a smooth\footnote{%
We restrict our attention to smooth schemes for simplicity.} $G$-scheme.
Edidin and Graham (\cite{EG1}), following an idea of Totaro (\cite{To2}),
defined a $G$-equivariant version, $A_G^{*}\left( X\right) $, of the Chow
ring $A^{*}\left( X\right) $. \label{BG}We will simply write $A_G^{*}$ for $%
A_G^{*}\left( Speck\right) $. As a rule, if we do not mention explicitly the
base field $k$, we are assuming $k=\mathbf{C}$.

We say that a pair $\left( U,V\right) $, consisting of a $k$-representation $%
V$ of $G$ and an open subset $U$ of $V$ on which $G$ acts freely, is a good
pair (or simply a pair) relative to $G$ if the codimension of $V\backslash U$
has sufficiently high codimension (see \cite{EG1}, 2.2,
Definition-Proposition 1).

All the basic properties and constructions (Chern classes, localization
sequence, proper pushforwards, Gysin maps, vector and projective bundle
theorems, projection formula, self intersection formula, cycle class map,
operational Chow groups etc.) of ordinary intersection theory (\cite{Fu})
have their equivariant counterparts. Moreover, there are additional
constructions one can do in the equivariant setting which simply do not
exist in the ordinary case, for example those related to morphisms of
algebraic groups. If 
\[
\varphi :G\longrightarrow G^{\prime } 
\]
is a morphism of algebraic groups and $X$ a $G^{\prime }$-scheme (which we
suppose smooth just in order to state each result for Chow rings), then $X$
is a $G$-scheme via $\varphi $ and if $\left( U,V\right) $ (respectively, $%
\left( U^{\prime },V^{\prime }\right) $) is a good pair relative to $G\,$%
(resp., relative to $G^{\prime }$), we let $G$ act on $V\times V^{\prime }$
as 
\[
g\cdot \left( v,v^{\prime }\right) =\left( g\cdot v,\pi \left( g\right)
\cdot v^{\prime }\right) ,\quad g\in G,\text{ }v\in V,\text{ }v^{\prime }\in
V^{\prime } 
\]
and the projection

\[
X\times U\times U^{\prime }\longrightarrow X\times U^{\prime } 
\]
induces a flat map 
\[
\left( X\times U\times U^{\prime }\right) /G\longrightarrow \left( X\times
U^{\prime }\right) /G^{\prime }. 
\]
Its pullback induces a restriction ring morphism 
\[
A_{G^{\prime }}^{*}\left( X\right) \longrightarrow A_G^{*}\left( X\right) 
\]
denoted by $\varphi _X^{*}$ (or by $\mathrm{res}_{G^{\prime },X}^G$ if $%
\varphi $ is injective). Note that the same construction made in the
topological case, define the functoriality in $G$ of the equivariant
cohomology ring $H_G^{*}\left( X;\mathbf{Z}\right) .$

Another construction which appears only in the equivariant setting is the
following transfer construction for Chow groups; we will frequently use it.
Let 
\[
1\rightarrow H\stackrel{\phi }{\longrightarrow }G\longrightarrow
F\rightarrow 1
\]
be an exact sequence of algebraic groups over a field $k$,$\,$with $F$
finite. If $X$ is an algebraic smooth $G$-scheme then $p_1:X\times
F\longrightarrow X$ is proper $G$-equivariant and there is an equivariant
push-forward 
\[
p_{1*}:A_{*}^G\left( X\times F\right) \rightarrow A_{*}^G\left( X\right) .
\]
If $\left( U,V\right) $ is a good pair for $G$, we have: 
\begin{equation}
\frac{\left( X\times F\right) \times U_{}}G\simeq \left( \frac{\left(
X\times F\right) \times U}H\right) \diagup F\simeq   \nonumber
\end{equation}
\[
\simeq \left( \frac{X\times U}H\times F\right) \diagup F\simeq \frac{X\times
U}H;
\]
hence $A_G^{*}\left( X\times F\right) \simeq A_H^{*}\left( X\right) $ and $%
p_{1*}$ induces a \emph{transfer} morphism of graded groups 
\begin{equation}
\mathrm{tsf}_{H,X}^G:A_H^{*}\left( X\right) \rightarrow A_G^{*}\left(
X\right) .  \label{tsf}
\end{equation}
which is natural in $X$ with respect to pullbacks.

Observe that the pullback $A_G^{*}\left( X\right) \rightarrow A_H^{*}\left(
X\right) $ has actually values in the $F$-invariant subring of $%
A_H^{*}\left( X\right) $ (\cite{To2}) 
\[
\mathrm{res}_{G,X}^H:A_G^{*}\left( X\right) \longrightarrow \left(
A_H^{*}\left( X\right) \right) ^F. 
\]
In exactly the same way as for group cohomology (e.g. \cite{B}, Prop. 9.5),
we have 
\[
\mathrm{tsf}_{H,X}^G\circ \mathrm{res}_{G,X}^H=\left( \#F\right) \cdot 
\]
(by projection formula) and, since $H$ is normal in $G$, 
\begin{equation}
\left( \mathrm{res}_{G,X}^H\circ \mathrm{tsf}_{H,X}^G\right) _{\mid \left(
A_H^{*}\left( X\right) \right) ^F}=\left( \#F\right) \cdot  \label{media1}
\end{equation}
If we do not restrict to $\left( A_H^{*}\left( X\right) \right) ^F$, we get 
\begin{equation}
\mathrm{res}_{G,X}^H\circ \mathrm{tsf}_{H,X}^G\left( \xi \right) =\sum_{f\in
F}f_{*}\xi  \label{media2}
\end{equation}
for any $\xi $ in $A_H^{*}\left( X\right) $.

\begin{remark}
For a general action of $G$ on $X$, the quotient $\left[ X/G\right] $ exists
only as an Artin stack\footnote{%
Not necessarily separated.} (\cite{LMB}). Edidin and Graham (\cite{EG1},
5.3, Prop. 16, 17) showed that if $\mathcal{F}$ is a quotient Artin stack $%
\mathcal{F}\simeq \left[ X/G\right] $, then the corresponding equivariant
Chow groups do not depend on the presentation chosen for the quotient,
enabling one to \emph{define }$A_G^{*}\left( X\right) $ to be the (integral)
Chow group of the stack $\mathcal{F}$. If moreover $\mathcal{F}$ is smooth,
there is a ring structure on this Chow group and this applies to the
classifying stack $\mathcal{B}G$ of any algebraic group $G$ (\cite{LMB}),
viewed as the quotient $\left[ \mathrm{pt}/G\right] $, 
\[
A_G^{*}=A^{*}\left( \mathcal{B}G\right) \text{. }
\]
\end{remark}

\begin{theorem}
\label{gottliebtotaro}(Gottlieb; Totaro) Let $G$ be an algebraic group over $%
\mathbf{C}$, $T$ a maximal torus of $G$ and $N_G\left( T\right) $ its
normalizer in $G$. The restriction maps 
\begin{equation}
A_G^{*}\longrightarrow A_{N_G\left( T\right) }^{*}  \label{totaro}
\end{equation}

\begin{equation}
H^{*}\left( \mathrm{B}G,\mathbf{Z}\right) \longrightarrow H^{*}\left( 
\mathrm{B}N_G\left( T\right) ,\mathbf{Z}\right)   \label{gottlieb}
\end{equation}
are injective.
\end{theorem}

\TeXButton{Proof}{\proof} (\ref{gottlieb}) is proved in \cite{Gtl}. We
sketch the proof of (\ref{totaro}) from \cite{To3}. If $f:Y\longrightarrow B$
is a smooth proper morphism of relative dimension $r$ between smooth,
separated schemes of finite type over $k$, let us consider the following
''modified'' pushforward 
\[
f_{\sharp }\left( \alpha \right) \doteq f_{*}\circ \left( c_r\left( \mathcal{%
T}_f\right) \cdot \alpha \right) \in A^j\left( S\right) 
\]
for any $\alpha \in A^j\left( B\right) $, where $\mathcal{T}_f$ denotes the
relative tangent bundle; by projection formula, we have 
\begin{equation}
f_{\sharp }\circ f^{*}=\chi \left( F\right)  \label{formula}
\end{equation}
where $\chi \left( F\right) $ denotes the Euler characteristic of ''the
fiber'' of $f$ (equal to the degree of the top Chern class of its tangent
bundle). Now, let $g:X\rightarrow B$ be a smooth morphism between smooth
schemes over a field $k$ which admits a smooth relative compactification 
\[
\begin{tabular}{lll}
$X$ & $\hookrightarrow $ & $\overline{X}$ \\ 
& $\searrow $ & $\downarrow $ \\ 
&  & $B$%
\end{tabular}
\]
having divisors with normal crossing $\left\{ D_i\right\} _{i=1,\ldots ,n}$
as complement (smooth over $B$). If $D_{ij}\doteq D_i\cap D_j$, $%
D_{ijk}\doteq D_i\cap D_j\cap D_k$, etc., the previous construction yields
modified pushforwards 
\[
\text{ }f_{\sharp }:A^{*}\left( \overline{X}\right) \rightarrow A^{*}\left(
B\right) ,\text{ }f_{\sharp }^{\left( 1\right) }:\oplus _iA^{*}\left(
D_i\right) \rightarrow A^{*}\left( B\right) ,\text{ }f_{\sharp }^{\left(
2\right) }:\oplus _{i<j}A^{*}\left( D_{ij}\right) \rightarrow A^{*}\left(
B\right) ,\ldots 
\]
satisfying (\ref{formula}). If $x\in A^{*}\left( X\right) $, lift it to some 
$\overline{x}\in A^{*}\left( \overline{X}\right) $ and set 
\[
g_{\sharp }\left( x\right) \doteq f_{\sharp }\left( \overline{x}\right)
-\sum_if_{\sharp }^{\left( 1\right) }\left( \overline{x}_{\mid D_i}\right)
+\sum_{i<j}f_{\sharp }^{\left( 2\right) }\left( \overline{x}_{\mid
D_{ij}}\right) -\cdots 
\]
(alternating sum) which is an element in $A^{*}\left( B\right) $. This can
be shown to be independent on the choice of the lifting and (\ref{formula})
holds for $g$ by well-known properties of the Euler characteristic.

To prove (\ref{totaro}), apply this construction to any approximation of the 
$G/N_G\left( T\right) $-torsor 
\[
\mathrm{B}N_G\left( T\right) \longrightarrow \mathrm{B}G 
\]
recalling that $\chi \left( G/N_G\left( T\right) \right) =1$. Note that this
proof works over any algebraically closed field $k$. \TeXButton{End Proof}
{\endproof}

In \cite{To2}, Totaro proved the remarkable fact that, if $G$ is a complex
algebraic group, the cycle map 
\[
\mathrm{cl}_{\mathrm{B}G}:A_G^{*}\longrightarrow H^{*}\left( \mathrm{B}G,%
\mathbf{Z}\right) 
\]
factors as 
\begin{equation}
A_{\mathrm{B}G}^{*}\stackrel{\widetilde{\mathrm{cl}}_{\mathrm{B}G}}{%
\longrightarrow }MU^{*}\left( \mathrm{B}G\right) \otimes _{MU^{*}}\mathbf{Z}%
\stackrel{\underline{\mathrm{cl}}_{\mathrm{B}G}}{\longrightarrow }%
H^{*}\left( \mathrm{B}G,\mathbf{Z}\right) .  \label{factor}
\end{equation}
where $MU^{*}\left( \mathrm{B}G\right) $ is the complex cobordism ring of $%
\mathrm{B}G$ (\cite{St}) and $\underline{\mathrm{cl}}_{\mathrm{B}G}$ is the
natural morphism (since 
\[
MU^{*}\equiv MU^{*}\left( \mathrm{pt}\right) =\mathbf{Z}\left[
x_1,x_2,\ldots ,x_n,\ldots \right] 
\]
with $\deg x_i=-2i$, here $\mathbf{Z}$ is viewed as an $MU^{*}$-module via
the map sending each generator $x_i$ to zero). We call $\widetilde{\mathrm{cl%
}}_{\mathrm{B}G}$ \emph{Totaro's refined cycle map} for $G$. The kernel and
cokernel of $\underline{\mathrm{cl}}_{\mathrm{B}G}$, $\widetilde{\mathrm{cl}}%
_{\mathrm{B}G}$ and $\mathrm{cl}_{\mathrm{B}G}$ are torsion.

In \cite{To2}, Totaro conjectures that if $G$ is a complex algebraic group
such that $MU^{*}\left( \mathrm{B}G\right) $, localized at some prime $p$,
is concentrated in even degrees, then the $p$-localization of $\widetilde{%
\mathrm{cl}}_{\mathrm{B}G}$ should be an isomorphism. As a consequence of
this conjecture, $A_{PGL_3}^{*}$ should not be generated by Chern classes
since, by \cite{KY}, $MU^{*}\left( \mathrm{B}PGL_3\right) $ is concentrated
in even degrees but not generated by Chern classes. This consequence of
Totaro's conjecture will be proved in section 4 (see also the Appendix for a
different proof).

\begin{remark}
Voevodsky (\cite{V1}, \cite{V2}) defined an \emph{algebraic cobordism}%
\textrm{\ }for an algebraic scheme over an \emph{arbitrary} field $k$, so it
would be interesting to investigate if there exists a generalization of
Totaro's refined cycle map with values in (a quotient of) algebraic
cobordism, for any algebraic group $G$ over $k$. As M. Levine suggested to
me, one may also ask more generally if Totaro's refined cycle map extends to
a map from the entire Voevodsky's motivic cohomology to algebraic cobordism.
\end{remark}

\subsection{Miscellaneous results}

\begin{proposition}
\label{rational} Let $k$ be algebraically closed. The pullback $%
A_{PGL_{n,k}}^{*}\otimes \mathbf{Q}\rightarrow A_{SL_{n,k}}^{*}\otimes 
\mathbf{Q}$ is an isomorphism.
\end{proposition}

\TeXButton{Proof}{\proof} By \cite{EG2}, Th. 1 (c), 
\[
A_G^{*}\otimes \mathbf{Q}\simeq \mathrm{Sym}_{\mathbf{Z}}\left( \widehat{T}%
\right) ^W\otimes \mathbf{Q=}\left( \mathbf{A}_T^{*}\right) ^W\otimes 
\mathbf{Q} 
\]
for any connected reductive algebraic group $G$ with maximal torus $T$ and
Weyl group $W$ and $\mathrm{Sym}_{\mathbf{Z}}\left( \widehat{T}\right)
^W\otimes \mathbf{Q}$ is the same for a group $G$ and a quotient of $G$ by a
finite central subgroup. \TeXButton{End Proof}{\endproof}

\begin{remark}
\label{equivalence}Let $S$ be a locally noetherian base scheme. Since $%
\mathrm{Aut}\left( \mathbf{P}_S^n\right) \simeq PGL_{n+1,S}$ as
group-functors, for any $S$ (\cite{DG} or \cite{GIT}, p. 20-21), the
category of Brauer-Severi schemes (\cite{Mi}, p. 134) of relative dimension $%
n$ over $X$ for the \'etale (or \textrm{fppf}) topology is equivalent to
that of $PGL_{n+1}$-torsors over $X$ for the same topology and this
equivalence actually extends to a $1$-isomorphism of $\mathcal{BS}_{n,S}$
with the classifying stack $\mathcal{B}\left( PGL_{n+1,S}\right) $, where $%
\mathcal{BS}_{n,S}$ denotes the stack over $S$ whose fibre category over $%
X/S $ is the category of Brauer-Severi schemes of relative dimension $n$
over $X$. Under this $1$-isomorphism trivial\footnote{%
I.e. of the form $P\left( E\right) \rightarrow X$ for some vector bundle $E$
over $X$.} Brauer-Severi schemes correspond to $PGL_{n+1}$-torsors induced
by $GL_{n+1}$-torsors.
\end{remark}

\begin{proposition}
Let $k$ be algebraically closed. Then \label{torsion}$\ker
(A_{PGL_{n,k}}^{*}\rightarrow A_{SL_{n,k}}^{*})$ is $n$-torsion.
\end{proposition}

\TeXButton{Proof}{\proof} By Prop. \ref{rational}, our kernel is torsion and
so it is enough to prove that $\ker (p^{*}:A_{PGL_{n,k}}^{*}\rightarrow
A_{GL_{n,k}}^{*})$ is annihilated by $n$, $A_{GL_{n,k}}^{*}$ being torsion
free.

By \cite{To2}, Th. 1.3 or \cite{EG2} Th. 1, for any reductive algebraic
group $G$, $A_G^{*}$ can be identified with the ring $\mathcal{C}_G^{*}$ of
characteristic classes for (\'etale) $G$-torsors over smooth, separated
schemes of finite type over $k$. Via this identification $p^{*}$ translates
to 
\begin{eqnarray*}
p^{*} &:&\mathcal{C}_{PGL_{n.k}}^{*}\longrightarrow \mathcal{C}%
_{GL_{n,k}}^{*} \\
F &\longmapsto &p^{*}(F):\left( 
\begin{array}{c}
E \\ 
\downarrow \\ 
X
\end{array}
\right) \mapsto F\left( 
\begin{tabular}{l}
$P_E$ \\ 
$\downarrow $ \\ 
$X$%
\end{tabular}
\right)
\end{eqnarray*}
where $P_E\rightarrow X$ is the $PGL_{n,k}$-torsor associated to $\mathbf{P}(%
\widetilde{E})\rightarrow X$, $\widetilde{E}\rightarrow X$ being the vector
bundle associated to the $GL_{n,k}$-torsor $E\rightarrow X$ and slightly
abusing notation in the argument of $F$ 
\[
p^{*}F=0\text{ }\Longleftrightarrow \text{ }F(\mathbf{P}(\widetilde{E}%
)\rightarrow X)=0,\text{ }\forall E\rightarrow X\text{ vector bundle of }rk%
\text{ }n. 
\]
Now we use the $1$-isomorphism of stacks $\mathcal{B}\left( PGL_{n,k}\right)
\simeq \mathcal{BS}_{n-1,k}$ (Remark \ref{equivalence}). If $f:P\rightarrow
X $ is a $PGL_{n,k}$-torsor and $\overline{f}:\overline{P}\rightarrow X$ the
associated Brauer-Severi scheme, the base change of $f$ via $\overline{f}$
is a $PGL_{n,k}$-torsor induced by a $GL_{n,k}$-torsor. Since $\chi \left(
P_k^{n-1}\right) =n$, formula (\ref{formula}) in the proof of Theorem \ref
{gottliebtotaro} yields 
\[
nF\left( 
\begin{tabular}{l}
$P$ \\ 
$\downarrow $ \\ 
$X$%
\end{tabular}
\right) =\overline{f}_{\sharp }\overline{f}^{*}F\left( 
\begin{tabular}{l}
$P$ \\ 
$\downarrow $ \\ 
$X$%
\end{tabular}
\right) =\overline{f}_{\sharp }F(\overline{f}^{*}\left( 
\begin{tabular}{l}
$P$ \\ 
$\downarrow $ \\ 
$X$%
\end{tabular}
\right) )=0 
\]
(by projection formula) if $p^{*}F=0$. \TeXButton{End Proof}{\endproof}

\begin{corollary}
\label{cortorsion}$A_{PGL_{n,k}}^{*}$ has only $n$-torsion.
\end{corollary}

We conclude this section collecting some elementary results on equivariant
Chow groups we will use in the sequel.

\begin{proposition}
$A_{\mathbf{\mu }_{n,k}}^{*}\simeq \mathbf{Z}\left[ t\right] \diagup \left(
nt\right) $.
\end{proposition}

\TeXButton{Proof}{\proof} From Kummer exact sequence 
\[
1\rightarrow \mathbf{\mu }_{n,k}\longrightarrow \mathbf{G}_{m,k}\stackrel{(%
\text{ })^n}{\longrightarrow }\mathbf{G}_{m,k}\rightarrow 1, 
\]
for any $N>0$ we get a $\mathbf{G}_{m,k}$-torsor 
\[
\frac{\mathbf{A}_k^{N+1}\backslash \left\{ 0\right\} }{\mathbf{\mu }_{n,k}}%
\longrightarrow \frac{\mathbf{A}_k^{N+1}\backslash \left\{ 0\right\} }{%
\mathbf{G}_{m,k}}=\mathbf{P}_k^N 
\]
whose associated line bundle is just $\mathcal{O}_{\mathbf{P}_k^N}(-n)$. By 
\cite{Gr1}, Remark p. 4-35, we get 
\[
A^{*}\left( \frac{\mathbf{A}_k^{N+1}\backslash \left\{ 0\right\} }{\mathbf{%
\mu }_{n,k}}\right) \simeq \frac{A^{*}\left( \mathbf{P}_k^N\right) }{\left(
c_1\left( \mathcal{O}_{\mathbf{P}_k^N}(-n)\right) \right) } 
\]
which implies the assert for $N\gg 0$. \TeXButton{End Proof}{\endproof}

\begin{proposition}
\label{unipotent}If $G$ is a unipotent algebraic group over a field $k$ of
characteristic zero, then $A_G^{*}\simeq \mathbf{Z}$.
\end{proposition}

\TeXButton{Proof}{\proof}Since $G$ is unipotent it has a central composition
series 
\[
G=G_0 > G_1 > G_2 > \cdots
> G_n=1 
\]
such that $G_i/G_{i+1}\simeq \mathbf{G}_{a,k}$ (\cite{DG}, IV, \S\ 2, 2.5
(vii)). We proceed by induction on the length $n$ of the composition series.

If $n=1$, $G\simeq \mathbf{G}_{a,k}$; if $U$ is a $G$-free open subset of a $%
G$-representation such that $\pi :U\rightarrow U/G$ is a (\textrm{fppf} or
\'etale) $G$-torsor then $\pi $ is a Zariski $G$-torsor ($\mathbf{G}_{a,k}$
being special, \cite{Se}) and in particular a Zariski affine bundle with
fiber $\mathbf{A}_k^1$ so that $\pi ^{*}$ is an isomorphism (\cite{Gr1}, p.
4-35).

Suppose the assert true for any unipotent group whose central composition
series has length $\leq n$. If $G$ is unipotent with a central decomposition
series 
\[
G=G_0 > G_1 > G_2 > \cdots
 > G_{n+1}=1 
\]
then $G_1$ is unipotent (\cite{DG}, IV, \S\ 2, 2.3) and we have a short
exact sequence 
\[
1\rightarrow G_1\longrightarrow G\longrightarrow G/G_1\simeq \mathbf{G}%
_{a,k}\rightarrow 1. 
\]
Therefore, if $U$ is a $G$-free open subset of a $G$-representation which
has a $G$-torsor quotient $U\rightarrow U/G$, 
\[
U/G_1\rightarrow U/G 
\]
is a $G/G_1\simeq \mathbf{G}_{a,k}^{}$-torsor. As in case $n=1$, the
pullback is an isomorphism $A_G^{*}\simeq A_{G_1}^{*}$ and we conclude since 
$G_1$ has a central decomposition series of length $n$. \TeXButton{End Proof}
{\endproof}

\begin{proposition}
\label{corunip}Let 
\begin{equation}
1\rightarrow H\longrightarrow G\stackrel{\rho }{\stackunder{\sigma }{%
\rightleftarrows }}\mathbf{G}_m\rightarrow 1  \label{extunip}
\end{equation}
be a split exact sequence of algebraic groups over a field $k$ of
characteristic zero, with $H$ unipotent. Then the pullback 
\[
\rho ^{*}:A_{\mathbf{G}_m}^{*}\longrightarrow A_G^{*} 
\]
is an isomorphism.
\end{proposition}

\TeXButton{Proof}{\proof} Let $U$ be a $G$-free open subset of a $G$%
-representation with complement of sufficiently high codimension and with a $%
G$-torsor quotient $U\rightarrow U/G$. Then 
\[
U/H\rightarrow U/G 
\]
is a $\mathbf{G}_m$-torsor which corresponds to some line bundle $L$ over $%
U/G$ (Th. \ref{equivalence}) and by \cite{Gr1}, Remark p. 4-35, 
\[
A^{*}\left( U/H\right) \simeq \frac{A^{*}\left( U/G\right) }{c_1(L)}. 
\]
Since $A_H^{*}\simeq \mathbf{Z}$, by Proposition (\ref{unipotent}), $A_G^{*}$
is then generated by $c_1(L)$. But the pullback $\mathbf{Z}[u]\simeq A_{%
\mathbf{G}_{m,k}}^{*}\rightarrow A_G^{*}$ sends $u$ to $c_1(L)$, therefore $%
\rho ^{*}$ is surjective. Injectivity follows from the hypothesis that (\ref
{extunip}) is split. \TeXButton{End Proof}{\endproof}

\begin{proposition}
\label{productbygm}If $G$ is an algebraic group over $k$, then $A_{G\times 
\mathbf{G}_{m,k}}^{*}\simeq A_G^{*}\otimes A_{\mathbf{G}_{m,k}}^{*}$.
\end{proposition}

\TeXButton{Proof}{\proof} Straightforward using $\left( \mathbf{A}%
_k^{N+1}\backslash \left\{ 0\right\} ,\mathbf{A}_k^{N+1}\right) $ as a good
pair for $G_{m,k}$, $N\gg 0$, and \cite{Fu}, Example 8.3.7. 
\TeXButton{End Proof}{\endproof}

\begin{proposition}
\label{homogeneous}Let $G$ be an algebraic group over $k$. If $H$ is a
closed algebraic subgroup of $G$, then there is a canonical isomorphism $%
A_G^{*}\left( G/H\right) \simeq A_H^{*}$.
\end{proposition}

\TeXButton{Proof}{\proof} Straightforward. \TeXButton{End Proof}{\endproof}

\begin{proposition}
\label{utrivial}Let $G$ be an algebraic group over a field $k$ and $X$ a
smooth $G$-scheme. If $U\subset \mathbf{A}_k^n$ is an open subscheme with
the trivial $G$-action, the pull-back $pr_2^{*}:A_G^{*}\left( X\right)
\simeq A_G^{*}\left( U\times X\right) $ is an isomorphism .
\end{proposition}

\TeXButton{Proof}{\proof} Since $G$ acts trivially on $U$, we can reduce to
the case of trivial $G$. By \cite{Fu}, Prop. 1.9, the pull back via $\mathbf{%
A}_k^n\times X\rightarrow X$ is surjective and so is $pr_2^{*}$ by the
localization exact sequence (\cite{Fu} Prop. 1.8).

If $k$ is infinite then $pr_2$ has always a section so that $pr_2^{*}$ is
injective. If $k$ is finite, let $p\in U$ be a closed point with $r\doteq $ $%
\left[ k(p):k\right] $. From the commutative diagram 
\[
\begin{tabular}{lll}
& $\quad U\times X$ &  \\ 
& $\nearrow \quad \downarrow ^{pr_2}$ &  \\ 
$p\times X$ & $\stackunder{\phi }{\longrightarrow }X$ & 
\end{tabular}
\]
and projection formula we get that $\ker \left( pr_2^{*}\right) $ is $r$%
-torsion. Now observe that we can always find two closed points $p$ and $%
p^{\prime }$ in $U$ with residue fields of relatively prime degrees $r$ and $%
r^{\prime }$ over $k$, so that $\ker \left( pr_2^{*}\right) $ is indeed
zero. \TeXButton{End Proof}{\endproof}

\begin{proposition}
\label{change}Let $G,H$ be algebraic groups having commuting actions on a
smooth scheme $X$ and suppose $G$ acts freely. Then there is a canonical
isomorphism $A_H^{*}\left( X/G\right) \simeq A_{G\times H}^{*}\left(
X\right) $.
\end{proposition}

\TeXButton{Proof}{\proof} If $\left( U,V\right) $ is a good pair for $H$,
with $\mathrm{codim}\left( V\backslash U\right) >i$, we have 
\[
A_H^i\left( X\diagup G\right) \simeq A^i\left( \left( U\times \frac
XG\right) \diagup H\right) \simeq 
\]
\[
\simeq A^i\left( \left( U\times X\right) \diagup G\times H\right) \simeq
A_{G\times H}^i\left( X\times U\right) ,\text{ } 
\]
by \cite{EG1}, Prop. 8. By the localization sequence, 
\[
A_{G\times H}^i\left( X\times U\right) \simeq A_{G\times H}^i\left( X\times
V\right) 
\]
for $i-$\textrm{codim}$(V\backslash U)<0$ and we conlude since for any $%
G\times H$-representation $E$, we have a pullback ring isomorphism $%
A_{G\times H}^{*}\left( X\right) \simeq A_{G\times H}^{*}\left( X\times
E\right) $.\TeXButton{End Proof}{\endproof}

%%%%%%%%%%%%%%%%%%%%%%%%%%% End PGL_32.tex %%%%%%%%%%%%%%%%%%%%%%%%%%%%%%%
}

\QSubDoc{Include PGL_33}{%%%%%%%%%%%%%%%%%%%%%%%%%% Start PGL_33.tex %%%%%%%%%%%%%%%%%%%%%%%%%%%%%%

\LaTeXparent{PGL_3}
\ChildStyles{amssymb} 
\ChildDefaults{chapter:3,page:1}

\section{ Generators for $A_{PGL_3}^{*}$}

From now on, our base field will be $\mathbf{C}$.

By Prop. \ref{rational}, we have 
\[
A_{PGL_3}^{*}\otimes \mathbf{Q}\simeq A_{SL_3}^{*}\otimes \mathbf{Q=Q}\left[
c_2\left( E\right) ,c_3\left( E\right) \right] 
\]
($E=$ standard representation of $SL_3$) and an easy computation shows that $%
c_3\left( E\right) $ is not in the image of the subring of $%
A_{PGL_3}^{*}\otimes \mathbf{Q}$ generated by the Chern classes of $\mathrm{%
sl}_3$. Therefore the Chern classes of the adjoint representation will
certainly not suffice to generate $A_{PGL_3}^{*}$.

In this section we find generators of $A_{PGL_3}^{*}$ (Prop. \ref{prelgen})
by stratifying the adjoint representation $\mathrm{sl}_3$ using Jordan
canonical forms.

Let $G$ be a complex algebraic group. For our purposes a finite $G$%
-stratification of a $G$-scheme $X$ will be a collection $\left\{
X_i\right\} _{i=1,\ldots ,n}$ of disjoint smooth $G$-invariant subschemes,
whose union is $X$ and such that for each $i$ the natural immersion 
\[
j_i:X_i\hookrightarrow X\diagdown \bigcup_{k>i}X_k\doteq X^i 
\]
is closed. In particular, $X_n$ is a closed subscheme of $X$, each $X_i$ is
topologically a locally closed subspace of $X$ and all the maps 
\[
X_1=X^1\hookrightarrow X^2\hookrightarrow X^3\hookrightarrow \cdots
\hookrightarrow X^{n-1}\hookrightarrow X^n\hookrightarrow X 
\]
are open immersions. Any stratification $\left\{ X_i\right\} _{i=1,\ldots
,n} $ gives then rise to the following exact sequences (of graded abelian
groups, $\deg \left( j_i\right) _{*}=$ $\mathrm{codim}_{X^i}\left(
X_i\right) $): 
\[
A_G^{*}\left( X_2\right) \stackrel{\left( j_2\right) _{*}}{\longrightarrow }%
A_G^{*}\left( X^2\right) \stackrel{i_2^{*}}{\longrightarrow }A_G^{*}\left(
X_1\right) \rightarrow 0 
\]
\begin{equation}
A_G^{*}\left( X_3\right) \stackrel{\left( j_3\right) _{*}}{\longrightarrow }%
A_G^{*}\left( X^3\right) \stackrel{i_3^{*}}{\longrightarrow }A_G^{*}\left(
X^2\right) \rightarrow 0  \label{uno}
\end{equation}
\[
\vdots 
\]
\[
A_G^{*}\left( X_n\right) \stackrel{\left( j_n\right) _{*}}{\longrightarrow }%
A_G^{*}\left( X=X^n\right) \stackrel{i_n^{*}}{\longrightarrow }A_G^{*}\left(
X^{n-1}=X\setminus X_n\right) \rightarrow 0. 
\]
Note that if $X$ is smooth, each graded group above is indeed a graded ring.
This will be our case.

Let 
\begin{eqnarray*}
U &\doteq &\left\{ A\in \mathrm{sl}_3\backslash \left\{ 0\right\} \mid \text{
}A\text{ has distinct eigenvalues}\right\} \stackunder{\text{open}}{\subset }%
\mathrm{sl}_3\backslash \left\{ 0\right\} , \\
\text{ }Z_{1,1} &\doteq &\left\{ A\in \mathrm{sl}_3\backslash \left\{
0\right\} \mid \text{ }A\text{ has Jordan form }\left( \text{%
\begin{tabular}{lll}
$\lambda $ & $0$ & $0$ \\ 
$1$ & $\lambda $ & $0$ \\ 
$0$ & $0$ & $-2\lambda $%
\end{tabular}
}\right) \text{, }\lambda \in \mathbf{C}^{*}\right\} , \\
Z_{1,0} &\doteq &\left\{ A\in \mathrm{sl}_3\backslash \left\{ 0\right\} \mid 
\text{ }A\text{ has Jordan form }\left( \text{%
\begin{tabular}{lll}
$\lambda $ & $0$ & $0$ \\ 
$0$ & $\lambda $ & $0$ \\ 
$0$ & $0$ & $-2\lambda $%
\end{tabular}
}\right) \text{, }\lambda \in \mathbf{C}^{*}\right\} , \\
Z_1 &\doteq &Z_{1,1}\cup Z_{1,0}\text{ ,} \\
Z_{0,1} &\doteq &\left\{ PGL_3\text{-orbit of }\left( \text{%
\begin{tabular}{lll}
$0$ & $0$ & $0$ \\ 
$1$ & $0$ & $0$ \\ 
$0$ & $1$ & $0$%
\end{tabular}
}\right) \right\} , \\
Z_{0,0} &\doteq &\left\{ PGL_3\text{-orbit of }\left( \text{%
\begin{tabular}{lll}
$0$ & $0$ & $0$ \\ 
$1$ & $0$ & $0$ \\ 
$0$ & $0$ & $0$%
\end{tabular}
}\right) \right\} , \\
Z_0 &\doteq &Z_{0.1}\cup Z_{0,0}\text{ ,} \\
&& \\
&&
\end{eqnarray*}
(note that $Z_1\cup Z_0=\mathrm{sl}_3\backslash (U\cup \left\{ 0\right\} )$%
). Then 
\begin{equation}
\left\{ U,Z_{1,1},Z_{1,0},Z_{0,1},Z_{0,0},\left\{ 0\right\} \right\}
\label{due}
\end{equation}
is a finite $PGL_3$-stratification of $\mathrm{sl}_3$. In this case the
first associated exact sequence of (\ref{uno}) is 
\begin{equation}
A_{PGL_3}^{*}\left( Z_{1,1}\right) \stackrel{\left( j_{1,1}\right) _{*}}{%
\longrightarrow }A_{PGL_3}^{*}\left( \mathrm{sl}_3\backslash (Z_{1,0}\cup
Z_0\cup \left\{ 0\right\} )\right) \stackrel{i_{1,1}^{*}}{\longrightarrow }%
A_{PGL_3}^{*}\left( U\right) \rightarrow 0  \label{tre}
\end{equation}
where $i_{1,1}:U\hookrightarrow \mathrm{sl}_3\backslash (Z_{1,0}\cup Z_0\cup
\left\{ 0\right\} )$ and $j_{1,1}:Z_{1,1}\hookrightarrow \mathrm{sl}%
_3\backslash (Z_{1,0}\cup Z_0\cup \left\{ 0\right\} )$ are the natural
immersions (open and closed, respectively).

To begin with, let us study $A_{PGL_3}^{*}\left( U\right) $.

\subsection{Generators coming from the open subset $U\subset \mathrm{sl}_3$}

Let $T$ be the maximal torus of $PGL_3$ and $\Gamma _3\doteq N_{PGL_3}\left(
T\right) =S_3\ltimes T$ its normalizer in $PGL_3$. Let $S_3\hookrightarrow
PGL_3:\sigma \mapsto \underline{\sigma }$ be the obvious inclusion (which
identifies permutations with permutation matrices). $\Gamma _3$ acts on the
subscheme $Diag_{\mathrm{sl}_3}^{*}\subset \mathrm{sl}_3\backslash \left\{
0\right\} $ of diagonal matrices with distinct eigenvalues, through $%
S_3\hookrightarrow PGL_3$ 
\[
\left( \sigma ,\left[ \underline{t}\right] \right) \cdot diag\left( \lambda
_1,\lambda _2,\lambda _3\right) \doteq \underline{\sigma }\cdot diag\left(
\lambda _1,\lambda _2,\lambda _3\right) \cdot \underline{\sigma }^{-1} 
\]
and we have\footnote{%
This Proposition holds (with the same proof given below) for any $PGL_n$.}:

\begin{proposition}
\label{p3.1}The composition of natural maps 
\[
A_{PGL_3}^{*}\left( U\right) \longrightarrow A_{\Gamma _3}^{*}\left(
U\right) \longrightarrow A_{\Gamma _3}^{*}\left( Diag_{\mathrm{sl}%
_3}^{*}\right) 
\]
is a ring isomorphism.
\end{proposition}

\TeXButton{Proof}{\proof} Let $T$ act by multiplication on the right of $%
PGL_3$ and $\frac{PGL_3}T$ be the corresponding quotient. $S_3$ acts on the
left of $\frac{PGL_3}T$ via $\sigma \cdot \left[ g\right] \doteq \left[
g\sigma ^{-1}\right] ,$ $g\in PGL_3$, and on $Diag_{\mathrm{sl}_3}^{*}$ as
above. If we let $PGL_3$ act on $Diag_{\mathrm{sl}_3}^{*}\times \frac{PGL_3}%
T $ by left multiplication on $\frac{PGL_3}T$ only, there is a $PGL_3$-
equivariant isomorphism 
\begin{eqnarray*}
U &\simeq &\left( Diag_{\mathrm{sl}_3}^{*}\times \frac{PGL_3}T\right)
\diagup S_3 \\
A &\longmapsto &\left[ \Delta ,\left[ g\right] _T\right] _{S_3}
\end{eqnarray*}
where $g^{-1}Ag=\Delta $ (diagonal).

Since $S_3$ acts freely on $Diag_{\mathrm{sl}_3}^{*}\times \frac{PGL_3}T$,
from Lemma \ref{change}, we get 
\[
A_{PGL_3}^{*}\left( U\right) \simeq A_{PGL_3\times S_3}^{*}\left( Diag_{%
\mathrm{sl}_3}^{*}\times \frac{PGL_3}T\right) . 
\]
Now, if $W$ is a free open subset of a $PGL_3\times S_3$-representation with
complement of sufficiently high codimension, we let $\Gamma _3$ act on $W$
via the inclusion 
\[
\left( i,\pi \right) :\Gamma _3\hookrightarrow PGL_3\times S_3:\left( \sigma
,\left[ \underline{t}\right] \right) \longmapsto \left( \left[ \underline{t}%
\right] \underline{\sigma },\sigma \right) 
\]
$i$ being the natural inclusion $\Gamma _3\hookrightarrow PGL_3$. Then the
morphisms 
\begin{eqnarray*}
&&\frac{W\times Diag_{\mathrm{sl}_3}^{*}\times \frac{PGL_3}T}{PGL_3\times S_3%
}\stackrel{\psi }{\stackunder{\varphi }{\leftrightarrows }}\frac{W\times
Diag_{\mathrm{sl}_3}^{*}}{\Gamma _3} \\
\varphi &:&\left[ w,\Delta ,\left[ g\right] _T\right] _{PGL_3\times
S_3}\longmapsto \left[ w\cdot \left( g,1\right) ,\Delta \right] _{\Gamma _3}
\\
\psi &:&\left[ w,\Delta \right] _{\Gamma _3}\longmapsto \left[ w,\Delta
,\left[ 1\right] _T\right] _{PGL_3\times S_3}
\end{eqnarray*}
are mutually inverse and we conclude. \TeXButton{End Proof}{\endproof}

\begin{lemma}
\label{invariants}If $T$ denotes the maximal torus of $PGL_3$ and $A_T^{*}\ $%
is viewed as a subring of $A_{T_{GL_3}}^{*}=\mathbf{Z}\left[
x_1,x_2,x_3\right] $, then the Weyl group-invariant subring $\left(
A_T^{*}\right) ^{S_3}$ is generated by 
\[
\gamma _2=s_1^2-3s_2
\]
\[
\gamma _3=2s_1^3-9s_1s_2+27s_3
\]
\[
\gamma _6=\Delta \equiv (x_1-x_2)^2(x_1-x_3)^2(x_2-x_3)^2
\]
where $s_i$ denotes the $i$-th elementary symmetric function on the $x_j$'s
and $\Delta $ is the discriminant.
\end{lemma}

\TeXButton{Proof}{\proof} We have $T=T_{PGL_3}\simeq T_{GL_3}/\mathbf{G}_m$,
where $T_{GL_3}=\left( \mathbf{G}_m\right) ^3$ and $G_m\hookrightarrow
T_{GL_3}$ diagonally. Therefore 
\[
A_T^{*}=Sym_{\mathbf{Z}}\left( \widehat{T}\right) \subset
A_{T_{GL_3}}^{*}=Sym_{\mathbf{Z}}\left( \widehat{T_{GL_3}}\right) =\mathbf{Z}%
\left[ x_1,x_2,x_3\right] 
\]
is the subring of polynomials $f\left( x_1,x_2,x_3\right) $ such that

\[
f\left( x_1+t,x_2+t,x_3+t\right) =f\left( x_1,x_2,x_3\right) . 
\]
Then 
\[
\left( A_T^{*}\right) ^{S_3}=\left\{ f\in \mathbf{Z}\left[
x_1,x_2,x_3\right] ^{S_3}\mid f\left( x_1+t,x_2+t,x_3+t\right) =f\left(
x_1,x_2,x_3\right) \right\} \equiv 
\]
\[
\equiv \left( \mathbf{Z}\left[ s_1,s_2,s_3\right] \right) ^{inv} 
\]
where $S_3$ permutes the $x_i$'s. Now, if for any polynomial $f\in \mathbf{Z}%
\left[ x_1,x_2,x_3\right] $ we let $f^t=f\left( x_1+t,x_2+t,x_3+t\right) $,
we get

\begin{eqnarray}
s_1^t &=&s_1+3t  \label{s1} \\
s_2^t &=&s_2+2s_1t+3t^2  \label{s2} \\
s_3^t &=&s_3+s_2t+s_1t^2+t^3  \label{s3}
\end{eqnarray}
and it is then easy to verify that $\gamma _2,\gamma _3$ and $\gamma _6$ are
indeed in $\left( A_T^{*}\right) ^{S_3}$.

Now, let $\varphi \in \left( A_T^{*}\right) ^{S_3}$. We first claim that
there exists $n_\varphi \geq 0$ such that $3^{n_\varphi }\varphi \in \mathbf{%
Z}\left[ \gamma _2,\gamma _3,\gamma _6\right] $. By definition of $\gamma _2$
and $\gamma _3$, we have 
\[
\left( A_T^{*}\right) ^{S_3}\left[ \frac 13\right] \equiv \left( \mathbf{Z}%
\left[ \frac 13\right] \left[ s_1,s_2,s_3\right] \right) ^{inv}=\left( 
\mathbf{Z}\left[ \frac 13\right] \left[ s_1,\gamma _2,\gamma _3\right]
\right) ^{inv}\text{.} 
\]
If 
\[
P\left( s_1,\gamma _2,\gamma _3\right) =P_0\left( \gamma _2,\gamma _3\right)
+P_1\left( \gamma _2,\gamma _3\right) s_1+\cdots +P_m\left( \gamma _2,\gamma
_3\right) s_1^m 
\]
is in $\in \left( \mathbf{Z}\left[ \frac 13\right] \left[ s_1,\gamma
_2,\gamma _3\right] \right) ^{inv}$, using (\ref{s1}) and $\gamma
_2^t=\gamma _2$, $\gamma _3^t=\gamma _3$, we easily get, by induction on $m$%
, $P_i=0,$ $\forall i\geq 1$, i.e. 
\[
\left( A_T^{*}\right) ^{S_3}\left[ \frac 13\right] =\mathbf{Z}\left[ \frac
13\right] \left[ \gamma _2,\gamma _3\right] 
\]
as claimed.

To prove that indeed $\varphi \in \mathbf{Z}\left[ \gamma _2,\gamma
_3,\gamma _6\right] $, we use induction on $n_\varphi $.

Suppose\footnote{%
Note that we allow an explicit dependence of $p$ on $\gamma _6$!} $3\varphi
=p\left( \gamma _2,\gamma _3,\gamma _6\right) $, for some polynomial $p$.
Expanding $p$ in powers of $\gamma _6$, we get 
\[
3\varphi =p_0\left( \gamma _2,\gamma _3\right) +p_1\left( \gamma _2,\gamma
_3\right) \gamma _6+\cdots 
\]
and reducing $\limfunc{mod}3$%
\[
0\equiv p_0\left( s_1^2,-s_1^3\right) +p_1\left( s_1^2,-s_1^3\right) \gamma
_6+\cdots \qquad \left( \limfunc{mod}3\right) . 
\]
But $s_1$ and $\gamma _6=\Delta $ are algebraically independent (over $%
\mathbf{Z}/3$), so $p_i\left( s_1^2,-s_1^3\right) \equiv 0$ $\left( \limfunc{%
mod}3\right) $, $\forall i$, i.e. 
\[
p_i\left( s_1^2,-s_1^3\right) =\left( \left( s_1^2\right) ^3-\left(
s_1^3\right) ^2\right) \cdot q_i\left( s_1^2,s_1^3\right) \qquad \left( 
\limfunc{mod}3\right) 
\]
then 
\[
p_i\left( \gamma _2,\gamma _3\right) =\left( \gamma _2^3-\gamma _3^2\right)
q_i\left( \gamma _2,-\gamma _3\right) +3r_i\left( \gamma _2,\gamma _3\right) 
\]
for each $i$. Thus 
\[
3\varphi =3r\left( \gamma _2,\gamma _3,\gamma _6\right) +\left( \gamma
_2^3-\gamma _3^2\right) q\left( \gamma _2,\gamma _3,\gamma _6\right) . 
\]
with an obvious notation. Straightforward computations yield 
\[
\left( \gamma _2^3-\gamma _3^2\right) =-3\left( \gamma _2^3-9\gamma
_6\right) \text{,} 
\]
and the case $n_\varphi =1$ is settled. The inductive step follows easily
from the fact that we included a possible dependence of $p$ from $\gamma _6$
in the above argument. \TeXButton{End Proof}{\endproof}

\begin{remark}
\label{twistaction}Note that there is a (non canonical) isomorphism 
\[
T\longrightarrow \left( \mathbf{G}_m\right) ^2
\]
\[
\left[ t_1,t_2,t_3\right] \longmapsto \left( t_1/t_3,t_2/t_3\right) 
\]
so that $A_T^{*}\simeq \mathbf{Z}\left[ x,y\right] $, with action of the
Weyl group given by 
\[
\left( 12\right) x=y,\text{ }\left( 12\right) y=x
\]
\begin{equation}
\left( 123\right) x=-y,\text{ }\left( 123\right) y=x-y.  \label{twistaction1}
\end{equation}
Under this isomorphism, with the same notations as in lemma \ref{invariants}%
, we have 
\[
\gamma _2=\left( x+y\right) ^2-3xy,
\]
\begin{equation}
\gamma _3=-9\left( x+y\right) xy+2\left( x+y\right) ^3  \label{twovariables}
\end{equation}
\[
\gamma _6=\left( x+y\right) ^2x^2y^2-4x^3y^3.
\]
Moreover, there is an isomorphism of $T$ with $T_{SL_3}$, the maximal torus
of $SL_3$%
\begin{equation}
T\rightarrow T_{SL_3}:\left[ t_1,t_2,t_3\right] \longmapsto \left(
t_2/t_3,t_3/t_1,t_1/t_2\right)   \label{twistaction2}
\end{equation}
and an induced isomorphism $A_T^{*}\simeq A_{T_{SL_3}}^{*}=\mathbf{Z}\left[
u_1,u_2,u_3\right] \diagup \left( u_1+u_2+u_3\right) $. The Weyl groups are
isomorphic to $S_3$ in both cases but the isomorphism above on Chow rings is
not $S_3$-equivariant, only $A_3$-equivariant. Rather, the action of $S_3$
on $A_{T_{SL_3}}^{*}$ inherited from the Weyl group action on $A_T^{*}$ via
this isomorphism, is given by 
\begin{eqnarray*}
\left( 12\right) u_1=-u_2,\text{ }\left( 12\right) u_2=-u_1,\,\left(
12\right) u_3=-u_3 \\
\left( 123\right) u_1=u_3,\text{ }\left( 123\right) u_2=u_1,\,\left(
123\right) u_3=u_1.
\end{eqnarray*}
\end{remark}

\begin{corollary}
\label{hsurj}The canonical morphism $h:A_{\Gamma _3}^{*}\rightarrow \left(
A_T^{*}\right) ^{S_3}$ is surjective.
\end{corollary}

\TeXButton{Proof}{\proof} Let $\phi :A_{PGL_3}^{*}\rightarrow A_{\Gamma
_3}^{*}$ be the restriction morphism, $E$ the standard representation of $%
GL_3$ and $Sym^3E$ be the $PGL_3$-representation: 
\[
\left[ g\right] \cdot \left( v\cdot _1v_2\cdot v_3\right) \doteq \left( \det
g^{-1}\right) \left( gv_1\cdot gv_2\cdot gv_3\right) .
\]
It is not difficult to verify that 
\begin{eqnarray*}
h\circ \phi \left( c_2\left( \mathrm{sl}_3\right) \right)  &=&-2\gamma _2 \\
\text{ }h\circ \phi \left( c_2\left( Sym^3E\right) \right)  &=&-5\gamma _2
\end{eqnarray*}
\begin{eqnarray*}
\text{ }h\circ \phi \left( c_3\left( Sym^3E\right) \right)  &=&\gamma _3 \\
\text{ }h\circ \phi \left( c_6\left( \mathrm{sl}_3\right) \right)  &=&\gamma
_6
\end{eqnarray*}
and the Corollary follows from Lemma \ref{invariants}. \TeXButton{End Proof}
{\endproof}

Now consider the subgroup $A_3\ltimes T\hookrightarrow \Gamma _3=S_3\ltimes
T $; there is a transfer morphism (see (\ref{tsf}), Section 2) 
\[
\mathrm{tsf}=\mathrm{tsf}_{A_3\ltimes T}^{\Gamma _3}\left( Diag_{\mathrm{sl}%
_3}^{*}\right) :A_{A_3\ltimes T}^{*}\left( Diag_{\mathrm{sl}_3}^{*}\right)
\rightarrow A_{\Gamma _3}^{*}\left( Diag_{\mathrm{sl}_3}^{*}\right) 
\]
and a restriction morphism: 
\[
\mathrm{res}=\mathrm{res}_{A_3\ltimes T}^{\Gamma _3}\left( Diag_{\mathrm{sl}%
_3}^{*}\right) :A_{\Gamma _3}^{*}\left( Diag_{\mathrm{sl}_3}^{*}\right)
\rightarrow \left( A_{A_3\ltimes T}^{*}\left( Diag_{\mathrm{sl}%
_3}^{*}\right) \right) ^{C_2}. 
\]

\begin{lemma}
\label{transfertrick}\emph{(transfer-trick)} $\mathrm{res}$ induces an
isomorphism 
\[
_3A_{\Gamma _3}^{*}\left( Diag_{\mathrm{sl}_3}^{*}\right) \rightarrow \text{ 
}_3\left( A_{A_3\ltimes T}^{*}\left( Diag_{\mathrm{sl}_3}^{*}\right) \right)
^{C_2}
\]
with inverse $(-\mathrm{tsf})$.
\end{lemma}

\TeXButton{Proof}{\proof} By projection formula, $\mathrm{tsf}\circ \mathrm{%
res}=2$; so if $\xi $ is $3$-torsion, we have $\mathrm{tsf}\circ \mathrm{res}%
\left( \xi \right) =-\xi $. On the other hand, if $C_2=\left\{ 1,\varepsilon
\right\} $, we have $\mathrm{res}\circ \mathrm{tsf}\left( \eta \right) =\eta
+\eta ^\varepsilon $; so, if $\eta $ is $C_2$-invariant and $3$-torsion, we
have $\mathrm{res}\circ \mathrm{tsf}\left( \eta \right) =-\eta $ and
conclude. \TeXButton{End Proof}{\endproof}

The isomorphism (\ref{twistaction2}) of Rmk. \ref{twistaction} induces an
isomorphism 
\[
A_3\ltimes T\simeq A_3\ltimes T_{SL_3} 
\]
and hence an isomorphism 
\begin{equation}
A_{A_3\ltimes T}^{*}\simeq A_{A_3\ltimes T_{SL_3}}^{*}.  \label{bingo}
\end{equation}
We will consider the $C_2$-action on $A_{A_3\ltimes T_{SL_3}}^{*}$ induced
by the canonical action on $A_{A_3\ltimes T}^{*}$ via this isomorphism. As
already in Rmk. \ref{twistaction}, we warn the reader that this is not the
canonical action induced by the inclusion $A_3\ltimes
T_{SL_3}\hookrightarrow N_{SL_3}\left( T_{SL_3}\right) $.

If $A_{A_3}^{*}=\mathbf{Z}\left[ \alpha \right] /(3\alpha )$ \footnote{%
We use that $A_3\simeq \mu _3$, which is true over any algebraically closed
field of characteristic $\neq 3$. Note that in characteristic $3$, it is no
longer true that $A_{A_3}^{*}\simeq \mathbf{Z}\left[ \alpha \right] /\left(
3\alpha \right) $.}, we still denote by $\alpha $ the image of $\alpha $ in $%
A_{A_3\ltimes T_{SL_3}}^{*}$ via the pullback induced by the projection $%
A_3\ltimes T_{SL_3}\rightarrow A_3$. We also recall the isomorphism $%
A_{T_{SL_3}}^{*}\simeq \mathbf{Z}\left[ u_1,u_2,u_3\right] \diagup \left(
u_1+u_2+u_3\right) $.

Then, if $W\simeq \mathbf{C}^3$ denotes the $A_3\ltimes T_{SL_3}$%
-representation 
\begin{equation}
\text{ }\left( \sigma ,\left( \underline{s}\right) \right) \cdot \left( 
\underline{x}\right) \doteq \left( s_1x_{\sigma ^{-1}\left( 1\right)
},s_2x_{\sigma ^{-1}\left( 2\right) },s_3x_{\sigma ^{-1}\left( 3\right)
}\right) \text{,}  \label{vudoppio}
\end{equation}
we have the following basic result

\begin{proposition}
\label{basico}The ring $A_{A_3\ltimes T_{SL_3}}^{*}$ is generated by 
\[
\left\{ \alpha ,\text{ }c_2\left( W\right) ,\text{ }c_3\left( W\right) ,%
\text{ }\theta \doteq \mathrm{tsf}_{T_{SL_3}}^{A_3\ltimes T_{SL_3}}\left(
u_2^2u_3\right) \right\} .
\]
\end{proposition}

\TeXButton{Proof}{\proof}Throughout the proof we identify $A_3\ltimes T$
with $A_3\ltimes T_{SL_3}$ (Remark \ref{twistaction}). $A_3\ltimes T_{SL_3}$
acts on $\mathbf{P}\left( W\right) $ with a dense orbit $U\doteq D_{+}\left(
x_1x_2x_3\right) $ with stabilizer isomorphic to $A_3\times \mathbf{\mu }_3$%
. If $j_2:Y_2\hookrightarrow \mathbf{P}\left( W\right) $ denotes the
(closed) orbit of $\left[ 1,0,0\right] \in \mathbf{P}\left( W\right) $ and 
\[
Y_1\doteq \mathbf{P}\left( W\right) \backslash U\cup Y_2\stackrel{j_1}{%
\stackunder{^{closed}}{\hookrightarrow }}\mathbf{P}\left( W\right)
\backslash Y_2, 
\]
the orbit of $\left[ 1,1,0\right] \in \mathbf{P}\left( W\right) $, then $%
\left\{ U,Y_1,Y_2\right\} $ is a finite $A_3\ltimes T_{SL_3}$-stratification
of $\mathbf{P}\left( W\right) $ and the exact sequences (\ref{uno}) are

\begin{equation}
A_{\mathbf{G}_m}^{*}\simeq A_{A_3\ltimes T_{SL_3}}^{*}\left( Y_1\right) 
\stackrel{\left( j_1\right) _{*}}{\longrightarrow }A_{A_3\ltimes
T_{SL_3}}^{*}\left( \mathbf{P}\left( W\right) \backslash Y_2\right) 
\stackrel{i^{*}}{\longrightarrow }A_{A_3\ltimes T_{SL_3}}^{*}\left( U\right)
\simeq A_{A_3\times \mathbf{\mu }_3}^{*}\rightarrow 0  \label{prima}
\end{equation}
\begin{equation}
A_{T_{SL_3}}^{*}\simeq A_{A_3\ltimes T_{SL_3}}^{*}\left( Y_2\right) 
\stackrel{\left( j_2\right) _{*}}{\longrightarrow }A_{A_3\ltimes
T_{SL_3}}^{*}\left( \mathbf{P}\left( W\right) \right) \stackrel{i_2^{*}}{%
\longrightarrow }A_{A_3\ltimes T_{SL_3}}^{*}\left( \mathbf{P}\left( W\right)
\backslash Y_2\right) \rightarrow 0.  \label{seconda}
\end{equation}
where we used Prop. \ref{homogeneous} together with the fact that $Y_1$
(resp. $Y_2$, $U$) has stabilizer isomorphic to $\mathbf{G}_m$ (resp. $%
T_{SL_3}$, $A_3\times \mathbf{\mu }_3$). By (\cite{Fu}, Th. 3.3 (b)), we
have ($c_1\left( W\right) =0$) : 
\[
A_{A_3\ltimes T_{SL_3}}^{*}\left( \mathbf{P}\left( W\right) \right) \simeq
A_{A_3\ltimes T_{SL_3}}^{*}\left[ \ell \right] \diagup \left( \ell
^3+c_2\left( W\right) \ell +c_3\left( W\right) \right) 
\]
where $\ell =c_1\left( \mathcal{O}_{\mathbf{P}\left( W\right) }\left(
1\right) \right) $. Moreover, the canonical K\"unneth morphism 
\[
A_{A_3}^{*}\otimes A_{\mathbf{\mu }_3}^{*}\simeq \mathbf{Z}\left[ \alpha
\right] /\left( 3\alpha \right) \otimes \mathbf{Z}\left[ \beta \right]
/\left( 3\beta \right) \rightarrow A_{A_3\times \mathbf{\mu }_3}^{*} 
\]
is an isomorphism (e.g. \cite{To2}, \S\ 6). It is not difficult to show that 
\[
i^{*}\left( \ell \right) =-\beta \text{, }i^{*}\left( \alpha \right) =\alpha 
\text{, }j_1^{*}\left( \ell \right) =-u 
\]
(where $A_{\mathbf{G}_m}^{*}=\mathbf{Z}\left[ u\right] $ and with the usual
abuse of notation, we write $\ell $ for $i_2^{*}\left( \ell \right) $ and $%
\alpha $ for its pullback to $A_{A_3\ltimes T_{SL_3}}^{*}\left( \mathbf{P}%
\left( W\right) \backslash Y_2\right) $). So we can conclude the analysis of
(\ref{prima}), by computing $\left( j_1\right) _{*}\left( 1\right) =\left[
Y_1\right] $. $Y_1$ is the zero scheme of the $A_3\ltimes T_{SL_3}$%
-invariant regular section $x_1x_2x_3\in \Gamma \left( \mathcal{O}\left(
3\right) ,\mathbf{P}\left( W\right) \backslash Y_2\right) $, hence (\cite{Fu}%
, p. 61), $\left[ Y_1\right] =3\ell $ so that $A_{A_3\ltimes
T_{SL_3}}^{*}\left( \mathbf{P}\left( W\right) \backslash Y_2\right) $ is
generated by $\left\{ \alpha ,\ell \right\} $.

Now let us turn our attention to (\ref{seconda}). It is easy to verify that,
with the usual abuse of notation, $i_2^{*}\left( \alpha \right) =\alpha $
and $i_2^{*}\left( \ell \right) =\ell $, so we are left to find generators
of $A_{A_3\ltimes T_{SL_3}}^{*}\left( Y_2\right) \simeq A_{T_{SL_3}}^{*}$ as
an $A_{A_3\ltimes T_{SL_3}}^{*}\left( \mathbf{P}\left( W\right) \right) $%
-module.

First of all, we have $j_2{}^{*}\left( \ell \right) =u_1$ \footnote{%
Of course this relation depends on the choice of the isomorphism 
\[
\mathbf{Z}\left[ u_1,u_2,u_3\right] \diagup \left( u_1+u_2+u_3\right)
=A_{T_{SL_3}}^{*}\simeq A_{A_3\ltimes T_{SL_3}}^{*}\left( Y_2\right) 
\]
which in its turn depends essentially on the choice of a point 
\[
p\in Y_2=\left\{ \left[ 1,0,0\right] ,\left[ 0,1,0\right] ,\left[
0,0,1\right] \right\} . 
\]
The choice we are making here is $p=\left[ 1,0,0\right] $.}. Therefore, by
projection formula and the relation $u_1+u_2+u_3=0$, we see that $%
A_{T_{SL_3}}^{*}\simeq A_{A_3\ltimes T_{SL_3}}^{*}\left( Y_2\right) $ is
generated by 
\begin{equation}
\left\{ 1,u_2^n\mid n>0\right\}  \label{first step}
\end{equation}
as an $A_{A_3\ltimes T_{SL_3}}^{*}\left( \mathbf{P}\left( W\right) \right) $%
-module. But 
\[
j_2{}^{*}\left( c_2\left( W\right) \right) =u_1u_2+u_2u_3+u_3u_1=-\left(
u_1^2+u_2^2+u_1u_2\right) 
\]
so that, by induction on $n$, $\left( j_2\right) _{*}\left( u_2^n\right) $, $%
n>1$, belongs to the submodule generated by $\left( j_2\right) _{*}\left(
1\right) $ and $\left( j_2\right) _{*}\left( u_2\right) $ (e.g. 
\begin{eqnarray*}
\left( j_2\right) _{*}\left( u_2^2\right) &=&\left( j_2\right) _{*}\left(
j_2^{*}\left( -c_2\left( W\right) \right) \right) +\left( j_2\right)
_{*}\left( j_2^{*}\left( -\ell ^2\right) \right) -\left( j_2\right)
_{*}\left( j_2^{*}\left( \ell \right) \cdot u_2\right) = \\
&=&-c_2\left( W\right) \cdot \left( j_2\right) _{*}\left( 1\right) -\ell
^2\cdot \left( j_2\right) _{*}\left( 1\right) -\ell \cdot \left( j_2\right)
_{*}\left( u_2\right)
\end{eqnarray*}
and similarly for higher powers of $u_2$). Thus, the ideal 
\[
im\left( j_2\right) _{*}\subset A_{A_3\ltimes T_{SL_3}}^{*}\left( \mathbf{P}%
\left( W\right) \right) 
\]
is actually generated by $\left( j_2\right) _{*}\left( 1\right) $ and $%
\left( j_2\right) _{*}\left( u_2\right) $.

Let us first compute $\left( j_2\right) _{*}\left( 1\right) $ using a
transfer argument (Section 2). Consider the $A_3\ltimes T_{SL_3}$-
equivariant commutative diagram 
\[\begin{tabular}{ll}
$Y_2\stackrel{h_2}{\longrightarrow }$ & $\mathbf{P}\left( W\right) \times A_3
$ \\ 
$\quad _{j_2}\searrow $ & $\quad \downarrow ^{pr_1}$ \\ 
& $\mathbf{P}\left( W\right) $%
\end{tabular}
\]

where 
\begin{eqnarray*}
h_2\left( \left[ 1,0,0\right] \right) =\left( \left[ 1,0,0\right] ,1\right) ,%
\text{ }h_2\left( \left[ 0,1,0\right] \right) =\left( \left[ 0,1,0\right]
,\sigma \right) ,\text{ }h_2\left( \left[ 0,0,1\right] \right) =\left(
\left[ 0,0,1\right] ,\sigma ^2\right) 
\end{eqnarray*}
with $\sigma =\left( 123\right) $. Using the canonical isomorphism 
\[
A_{A_3\ltimes T_{SL_3}}^{*}\left( \mathbf{P}(W)\times A_3\right) \simeq
A_T^{*}\left( \mathbf{P}(W)\right) 
\]
we see that 
\begin{equation}
\left( j_2\right) _{*}\left( 1\right) =\left( pr_1\right) _{*}\circ \left(
h_2\right) _{*}\left( 1\right) =\mathrm{tsf}_{T_{SL_3}}^{A_3\ltimes
T_{SL_3}}\left( \mathbf{P}(W)\right) \left( \left[ \left\{ \left[
1,0,0\right] \right\} \right] \right) .  \label{alfa}
\end{equation}
But $\left[ 1,0,0\right] =Z\left( x_2\right) \cap Z\left( x_3\right) $,
where the section $x_i$, $i=2,3$ are $T_{SL_3}$-semi-invariant (\cite{SGA3.I}%
, Expos\'e VI$_{\mathrm{B}},$ p. 406) so that if we consider the $T_{SL_3}$%
-equivariant line bundles $L_i\rightarrow Spec\mathbf{C}$ associated to the
representations 
\[
\left( \underline{t}\right) x=t_ix\text{,\quad \quad }i=2,3\text{,}
\]
we have induced $T_{SL_3}$-invariant regular sections $\widetilde{x_i}\in
\Gamma \left( P(W),\mathcal{O}\left( 1\right) \otimes p^{*}\left( L_i^{\vee
}\right) \right) $ \footnote{%
Note that $p^{*}\left( L_i^{\vee }\right) $ is trivial but not $T_{SL_3}$%
-equivariantly trivial.} with, obviously, $Z\left( \widetilde{x_i}\right)
=Z\left( x_i\right) $. Then 
\begin{equation}
\left[ \left\{ \left[ 1,0,0\right] \right\} \right] =\left( \ell -u_2\right)
\left( \ell -u_3\right) =\ell ^2+\ell u_1+u_2u_3  \label{beta}
\end{equation}
in $A_{T_{SL_3}}^{*}\left( \mathbf{P}(W)\right) $. Since $\ell =\mathrm{res}%
_{A_3\ltimes T_{SL_3}}^{T_{SL_3}}\left( \mathbf{P}(W)\right) \left( \ell
\right) $ and the diagram 
\[
\begin{tabular}{ccc}
$\qquad \qquad A_{T_{SL_3}}^{*}$ & $\longrightarrow $ & $A_{T_{SL_3}}^{*}%
\left( \mathbf{P}(W)\right) $ \\ 
$^{\mathrm{tsf}_{T_{SL_3}}^{A_3\ltimes T_{SL_3}}}\downarrow $ &  & $\qquad
\quad \downarrow ^{\mathrm{tsf}_{T_{SL_3}}^{A_3\ltimes T_{SL_3}}\left( 
\mathbf{P}(W)\right) }$ \\ 
$\qquad \qquad A_{A_3\ltimes T_{SL_3}}^{*}$ & $\longrightarrow $ & $%
A_{A_3\ltimes T_{SL_3}}^{*}\left( \mathbf{P}(W)\right) $%
\end{tabular}
\]
is commutative, we have 
\begin{eqnarray}
\mathrm{tsf}_{T_{SL_3}}^{A_3\ltimes T_{SL_3}}\left( \mathbf{P}(W)\right)
\left( \ell ^2\right)  &=&3\ell ^2  \label{gamma} \\
\mathrm{tsf}_{T_{SL_3}}^{A_3\ltimes T_{SL_3}}\left( \mathbf{P}(W)\right)
\left( \ell u_1\right)  &=&\ell \cdot \mathrm{tsf}_{T_{SL_3}}^{A_3\ltimes
T_{SL_3}}\left( u_1\right) .  \label{gamma2}
\end{eqnarray}
Now we claim $\mathrm{tsf}_{T_{SL_3}}^{A_3\ltimes T_{SL_3}}\left( u_1\right)
=0$, $i=1,2,3$. In fact, let $\pi :A_{A_3\ltimes T_{SL_3}}^{*}\rightarrow A_3
$ be the projection and $\rho :A_3\hookrightarrow A_{A_3\ltimes T_{SL_3}}^{*}
$ its right inverse. Since $i_2^{*}$ in (\ref{seconda}) is an isomorphism in
degree $1$ and $A_{A_3\ltimes T_{SL_3}}^{*}\left( \mathbf{P}\left( W\right)
\backslash Y_2\right) $ is generated by $\alpha $ and $\ell $, $\mathrm{tsf}%
_{T_{SL_3}}^{A_3\ltimes T_{SL_3}}\left( u_i\right) =n_i\pi ^{*}\alpha $ for
some integer $n_i$ (in fact 
\[
\mathrm{res}_{A_3\ltimes T_{SL_3}}^{T_{SL_3}}\circ \mathrm{tsf}%
_{T_{SL_3}}^{A_3\ltimes T_{SL_3}}\left( u_i\right) =u_1+u_2+u_3=0
\]
thus $\mathrm{tsf}_{T_{SL_3}}^{A_3\ltimes T_{SL_3}}\left( u_i\right) $ is $3$%
-torsion). Since 
\[
\rho ^{*}\circ \mathrm{tsf}_{T_{SL_3}}^{A_3\ltimes T_{SL_3}}\equiv \mathrm{%
res}_{A_3}^{A_3\ltimes T_{SL_3}}\circ \mathrm{tsf}_{T_{SL_3}}^{A_3\ltimes
T_{SL_3}}=0,
\]
we get 
\[
n_i\rho ^{*}\pi ^{*}\left( \alpha \right) =n_i\alpha =0
\]
in $A_{A_3}^{*}$ and the claim follows.

Since $\left( j_2\right) _{*}\left( u_2\right) $ has degree $3$, from (\ref
{seconda}) and the computations we have just done (in particular (\ref{alfa}%
), (\ref{beta}), (\ref{gamma}) and (\ref{gamma2})), we know that the ring $%
A_{A_3\ltimes T_{SL_3}}^{*}\left( \mathbf{P}\left( W\right) \right) $ is
generated up to degree $2$ (included) by

\[
\left\{ \alpha ,\ell ,\mathrm{tsf}_{T_{SL_3}}^{A_3\ltimes T_{SL_3}}\left(
u_2u_3\right) \right\} . 
\]
We will show that:

\texttt{CLAIM}. $A_{A_3\ltimes T_{SL_3}}^{*}\left( \mathbf{P}\left( W\right)
\right) $ is generated up to degree $2$ (included) by 
\[
\left\{ \alpha ,\ell ,c_2\left( W\right) \right\} . 
\]

\underline{Proof of Claim}. We write 
\[
\eta _{\mid T_{SL_3}}\equiv \mathrm{res}_{A_3\ltimes
T_{SL_3}}^{T_{SL_3}}\left( \eta \right) , 
\]
for any $\eta \in A_{A_3\ltimes T_{SL_3}}^{*}\left( \mathbf{P}\left(
W\right) \right) $.

Observe that 
\[
\mathrm{res}_{A_3\ltimes T_{SL_3}}^{T_{SL_3}}\circ \mathrm{tsf}%
_{T_{SL_3}}^{A_3\ltimes T_{SL_3}}\left( u_2u_3\right)
=u_2u_3+u_3u_1+u_1u_2=c_2\left( W\right) _{\mid T_{SL_3}}\text{,} 
\]
therefore $\mathrm{tsf}_{T_{SL_3}}^{A_3\ltimes T_{SL_3}}\left( u_2u_3\right)
-c_2\left( W\right) =\xi $, for some $3$-torsion element\footnote{%
In fact $\mathrm{tsf}_{T_{SL_3}}^{A_3\ltimes T_{SL_3}}\circ \mathrm{res}%
_{A_3\ltimes T_{SL_3}}^{T_{SL_3}}=3$.} 
\[
\xi \in A_{A_3\ltimes T_{SL_3}}^2\left( \mathbf{P}\left( W\right) \right) . 
\]
Since the group $A_{A_3\ltimes T_{SL_3}}^2\left( \mathbf{P}\left( W\right)
\right) $ is generated by 
\[
\left\{ \alpha ^2,\ell ^2,\alpha \ell ,\mathrm{tsf}_{T_{SL_3}}^{A_3\ltimes
T_{SL_3}}\left( u_2u_3\right) =c_2\left( W\right) +\xi \right\} 
\]
we have 
\begin{equation}
c_2\left( W\right) =A\left( c_2\left( W\right) +\xi \right) +B\alpha
^2+C\ell ^2+D\alpha \ell \text{.}  \label{delta}
\end{equation}
Restricting to $T_{SL_3}$, we get 
\[
c_2\left( W\right) _{\mid T_{SL_3}}=Ac_2\left( W\right) _{\mid
T_{SL_3}}+C\ell ^2; 
\]
but from 
\[
A_{T_{SL_3}}^{*}\left( \mathbf{P}(W)\right) \simeq A_{T_{SL_3}}^{*}\left[
\ell \right] \diagup \left( \ell ^3+\ell ^2c_1\left( W\right) _{\mid
T_{SL_3}}+\ell c_2\left( W\right) _{\mid T_{SL_3}}+c_3\left( W\right) _{\mid
T_{SL_3}}\right) , 
\]
we see that $c_2\left( W\right) _{\mid T_{SL_3}}$ and $\ell ^2$ are
algebraically independent, so we must have $A=1$, $C=0$. Thus (\ref{delta})
yields $\xi =B\alpha ^2+D\alpha \ell $ and this concludes the proof of Claim.%
$\ \mathbf{\Box }$

So, the other possible generators of $A_{A_3\ltimes T_{SL_3}}^{*}\left( 
\mathbf{P}\left( W\right) \right) $ in degree $>2$ can only come from $%
\left( j_2\right) _{*}\left( u_2\right) $. Using the same arguments as in
the computation of $\left( j_2\right) _{*}\left( 1\right) $ above, we get 
\[
\left( j_2\right) _{*}\left( u_2\right) =\mathrm{tsf}_{T_{SL_3}}^{A_3\ltimes
T_{SL_3}}\left( \mathbf{P}(W)\right) \left( u_2\left( \ell -u_2\right)
\left( \ell -u_3\right) \right) \text{.} 
\]
But, since we know that $\mathrm{tsf}_{T_{SL_3}}^{A_3\ltimes T_{SL_3}}\left(
u_i\right) =0$ $\forall i$, the only new generator is $\mathrm{tsf}%
_{T_{SL_3}}^{A_3\ltimes T_{SL_3}}\left( u_2^2u_3\right) $.

To summarize, we have proved so far that $A_{A_3\ltimes T_{SL_3}}^{*}\left( 
\mathbf{P}\left( W\right) \right) $ is generated by 
\[
\left\{ \alpha ,\ell ,c_2(W),\mathrm{tsf}_{T_{SL_3}}^{A_3\ltimes
T_{SL_3}}\left( u_2^2u_3\right) \right\} \text{.} 
\]
Since 
\[
A_{A_3\ltimes T_{SL_3}}^{*}\diagup \left( c_3\left( W\right) \right) \simeq
A_{A_3\ltimes T_{SL_3}}^{*}\left( W\backslash \left\{ 0\right\} \right)
\simeq A_{A_3\ltimes T_{SL_3}}^{*}\left( \mathbf{P}\left( W\right) \right)
\diagup \left( \ell \right) , 
\]
we conclude that $A_{A_3\ltimes T_{SL_3}}^{*}$ is generated by 
\[
\left\{ \alpha ,c_2(W),c_3(W),\mathrm{tsf}_{T_{SL_3}}^{A_3\ltimes
T_{SL_3}}\left( u_2^2u_3\right) \right\} . 
\]
\TeXButton{End Proof}{\endproof}

Recall (\ref{bingo}) and the isomorphism 
\[
_3A_{\Gamma _3}^{*}\left( Diag_{\mathrm{sl}_3}^{*}\right) \simeq \text{ }%
_3\left( A_{A_3\ltimes T}^{*}\left( Diag_{\mathrm{sl}_3}^{*}\right) \right)
^{C_2} 
\]
from lemma \ref{transfertrick}. If $C_2=\left\{ 1,\varepsilon \right\} $, we
denote by $W^\varepsilon $ the $A_3\ltimes T_{SL_3}$-representation obtained
from $W$ twisting the action by $\varepsilon $. Let us also define the
element 
\begin{equation}
\underline{\chi }\doteq \left( 2\mathrm{tsf}_{T_{SL_3}}^{A_3\ltimes
T_{SL_3}}\left( u_2^2u_3\right) +3c_3\left( W\right) \right) ^2+4c_2\left(
W\right) ^3+27c_3\left( W\right) ^2\in A_{A_3\ltimes T_{SL_3}}^6  \label{chi}
\end{equation}
and denote $\mathrm{tsf}_{T_{SL_3}}^{A_3\ltimes T_{SL_3}}\left(
u_2^2u_3\right) $ simply by $\theta $.

\begin{lemma}
\label{basico2}\textrm{(i) }In $A_{A_3\ltimes T}^{*}$ we have 
\[
3\underline{\chi }=3\alpha =\alpha \theta =\alpha ^3+\alpha c_2\left(
W\right) =0.
\]

\textrm{(ii)} The kernel of the restriction map $h^{\prime }:A_{A_3\ltimes
T}^{*}\rightarrow A_{A_3\ltimes T}^{*}\left( Diag_{\mathrm{sl}_3}^{*}\right) 
$ is the ideal $\left( \alpha ^2\right) $.

\textrm{(iii)} In $A_{A_3\ltimes T}^{*}$, we have 
\[
c_2\left( W^\varepsilon \right) =c_2\left( W\right) ,\qquad c_3\left(
W^\varepsilon \right) =-c_3\left( W\right) ,
\]

\[
\theta ^\varepsilon =\theta +3c_3\left( W\right) ,\qquad \underline{\chi }%
^\varepsilon =\underline{\chi }.
\]
\textrm{(iv)} Let 
\[
q\left( c_2(W),c_3(W),\mathrm{tsf}_{T_{SL_3}}^{A_3\ltimes T_{SL_3}}\left(
u_2^2u_3\right) \right) \in \text{ }_3A_{A_3\ltimes T_{SL_3}}^{*}
\]
be a polynomial in the arguments indicated. Then there exists a polynomial 
\[
\widetilde{q}=\widetilde{q}\left( c_2(W),c_3(W),\mathrm{tsf}%
_{T_{SL_3}}^{A_3\ltimes T_{SL_3}}\left( u_2^2u_3\right) \right) 
\]
such that $q=\underline{\chi }\widetilde{q}$.
\end{lemma}

\TeXButton{Proof}{\proof} (i) Since 
\[
\left( 2\mathrm{tsf}_{T_{SL_3}}^{A_3\ltimes T_{SL_3}}\left( u_2^2u_3\right)
+3c_3\left( W\right) \right) _{\mid T_{SL_3}}^2=\Delta \left(
u_1,u_2,u_3\right) 
\]
\begin{eqnarray*}
c_2\left( W\right) _{\mid T_{SL_3}} &=&s_2\left( u_1,u_2,u_3\right) \\
c_3\left( W\right) _{\mid T_{SL_3}} &=&s_3\left( u_1,u_2,u_3\right)
\end{eqnarray*}
in $A_{T_{SL_3}}^{*}\simeq A_T^{*}$ (where $\Delta $ is the discriminant and 
$s_i$ the $i$-th elementary symmetric function), it is well known that $%
\underline{\chi }_{\mid T}=0$. Therefore $3\underline{\chi }=0$. $\alpha $
is $3$-torsion by definition and 
\[
\alpha \cdot \mathrm{tsf}_{T_{SL_3}}^{A_3\ltimes T_{SL_3}}\left(
u_2^2u_3\right) =0 
\]
by projection formula. Finally observe that 
\[
\left( c_2\left( W\right) -\mathrm{tsf}_{T_{SL_3}}^{A_3\ltimes
T_{SL_3}}\left( u_1u_3\right) \right) _{\mid T_{SL_3}}=0 
\]
and therefore (Proposition \ref{basico}) there exist $A,B\in \mathbf{Z}$
such that 
\[
c_2\left( W\right) -\mathrm{tsf}_{T_{SL_3}}^{A_3\ltimes T_{SL_3}}\left(
u_1u_2\right) =A\alpha ^2+Bc_2\left( W\right) 
\]
is a $3$-torsion element in $A_{A_3\ltimes T_{SL_3}}^{*}$. Restricting to $%
T_{SL_3}$ we get $B=0$ while restricting to $A_3$ we get $A\equiv -1$ $%
\limfunc{mod}3$. Multiplying by $\alpha $, we get 
\[
\alpha ^3+\alpha c_2\left( W\right) =0 
\]
by projection formula.

(ii) A straightforward computation yields 
\[
c_2\left( Diag_{\mathrm{sl}_3}\right) =-\alpha ^2\in A_{A_3\ltimes T}^{*}. 
\]
Consider then the two localization sequences: 
\begin{equation}
A_{A_3\ltimes T}^{*}\stackrel{(-\alpha ^2)\cdot }{\longrightarrow }%
A_{A_3\ltimes T}^{*}\left( Diag_{\mathrm{sl}_3}\right) \simeq A_{A_3\ltimes
T}^{*}\longrightarrow A_{A_3\ltimes T}^{*}\left( Diag_{\mathrm{sl}%
_3}\backslash \left\{ 0\right\} \right) \rightarrow 0  \label{luno}
\end{equation}

\begin{equation}
A_T^{*}\simeq A_{A_3\ltimes T}^{*}\left( Z\right) \stackrel{j_{*}}{%
\longrightarrow }A_{A_3\ltimes T}^{*}\left( Diag_{\mathrm{sl}_3}\backslash
\left\{ 0\right\} \right) \longrightarrow A_{A_3\ltimes T}^{*}\left( Diag_{%
\mathrm{sl}_3}^{*}\right) \rightarrow 0  \label{ldue}
\end{equation}
(where we used the obvious $A_3\ltimes T$-equivariant isomorphism $Z\simeq
A_3\times \mathbf{C}^{*}$); (\ref{luno}) shows that $\alpha ^2\in \ker
h^{\prime }$ and the reverse inclusion will be established if we show that
the push-forward $j_{*}$ is zero.

Consider the projectivization $\mathbf{P}\left( Diag_{\mathrm{sl}_3}\right)
\simeq \mathbf{P}^1$ of $Diag_{\mathrm{sl}_3}$. We have a cartesian diagram 
\[
\begin{tabular}{lll}
$Z$ & $\stackrel{j}{\hookrightarrow }$ & $Diag_{\mathrm{sl}_3}\backslash
\left\{ 0\right\} $ \\ 
$^p\downarrow $ &  & $\qquad \downarrow ^\pi $ \\ 
$Z^{\prime }$ & $\stackunder{j^{\prime }}{\hookrightarrow }$ & $\mathbf{P}%
\left( Diag_{\mathrm{sl}_3}\right) $%
\end{tabular}
\]
where 
\[
Z^{\prime }=\left\{ \left[ 1,1\right] ,\left[ -2,1\right] ,\left[
1,-2\right] \right\} \simeq A_3 
\]
$A_3\ltimes T$-equivariantly. Since 
\[
j_{*}\circ p^{*}=\pi ^{*}\circ j_{*}^{\prime } 
\]
and $p^{*}$ is obviously an isomorphism, it is enough to show that 
\begin{equation}
im\left( j_{*}^{\prime }\right) \subseteq \ker \left( \pi ^{*}\right)
=\left( \ell \right) \subset \frac{A_{A_3\ltimes T}^{*}\left[ \ell \right] }{%
\left( \ell ^2-\alpha ^2\right) }  \label{ltre}
\end{equation}
by the projective bundle theorem.

To compute $j_{*}^{\prime }$ we translate it into a transfer map. Consider
the $A_3\ltimes T$-equivariant commutative diagram 
\[
\begin{tabular}{lll}
$Z^{\prime }\stackrel{\rho }{\longrightarrow }$ & $A_3\times \mathbf{P}%
\left( Diag_{\mathrm{sl}_3}\right) $ &  \\ 
$^{j^{\prime }}\searrow $ & $\qquad \downarrow ^{pr_2}$ &  \\ 
& $\mathbf{P}\left( Diag_{\mathrm{sl}_3}\right) $ & 
\end{tabular}
\]
where ($\sigma =\left( 123\right) \in A_3$) 
\begin{eqnarray*}
\rho \left( \left[ 1,1\right] \right) &=&\left( 1,\left[ 1,1\right] \right)
\\
\rho \left( \left[ -2,1\right] \right) &=&\left( \sigma ,\left[ -2,1\right]
\right) \\
\rho \left( \left[ 1,-2\right] \right) &=&\left( \sigma ^2,\left[
1,-2\right] \right) .
\end{eqnarray*}
Since 
\[
A_{A_3\ltimes T}^{*}\left( A_3\times \mathbf{P}\left( Diag_{\mathrm{sl}%
_3}\right) \right) \simeq A_T^{*}\left( \mathbf{P}\left( Diag_{\mathrm{sl}%
_3}\right) \right) , 
\]
we have 
\[
j_{*}^{\prime }\left( \xi \right) =pr_{2*}\circ \rho _{*}\left( \xi \right) =%
\mathrm{tsf}_T^{A_3\ltimes T}\left( \mathbf{P}\left( Diag_{\mathrm{sl}%
_3}\right) \right) \left( \xi \cdot \left[ \left\{ \left[ 1,1\right]
\right\} \right] \right) 
\]
for any $\xi \in A_T^{*}\simeq A_{A_3\ltimes T}^{*}\left( Z^{\prime }\right) 
$, where $\left\{ \left[ 1,1\right] \right\} $ is a $T$-invariant cycle on $%
\mathbf{P}\left( Diag_{\mathrm{sl}_3}\right) $.

Now, $\left\{ \left[ 1,1\right] \right\} $ is the zero scheme of the $T$%
-invariant regular section $\left( x_1-x_2\right) \in \Gamma \left( \mathbf{P%
}\left( Diag_{\mathrm{sl}_3}\right) ,\mathcal{O}\left( 1\right) \right) $,
therefore 
\[
\left[ \left\{ \left[ 1,1\right] \right\} \right] =c_1\left( \mathcal{O}%
\left( 1\right) \right) \equiv \ell ^{\prime }\in A_T^{*}\left( \mathbf{P}%
\left( Diag_{\mathrm{sl}_3}\right) \right) 
\]
and, obviously, 
\[
\mathrm{res}_{A_3\ltimes T}^T\left( \mathbf{P}\left( Diag_{\mathrm{sl}%
_3}\right) \right) \left( \ell \right) =\ell ^{\prime }. 
\]
By projection formula, we then get 
\[
j_{*}^{\prime }\left( \xi \right) =\mathrm{tsf}_T^{A_3\ltimes T}\left( 
\mathbf{P}\left( Diag_{\mathrm{sl}_3}\right) \right) \left( \xi \cdot \ell
^{\prime }\right) =\ell \cdot \mathrm{tsf}_T^{A_3\ltimes T}\left( \xi
\right) 
\]
for any $\xi \in A_T^{*}\simeq A_{A_3\ltimes T}^{*}\left( Z^{\prime }\right) 
$, which proves (\ref{ltre}).

(iii) By Prop. \ref{basico}, there are integers $A,$ $B$ such that 
\[
c_2\left( W^\varepsilon \right) =A\alpha ^2+Bc_2\left( W\right) . 
\]
Restricting this to $T$, we get $B=1$ and applying the involution $\left(
\cdot \right) ^\varepsilon $ we obtain $A\equiv 0$ $\limfunc{mod}3$.

Again by Prop. \ref{basico}, there are integers $A,B,C,D$ such that 
\[
c_3\left( W^\varepsilon \right) =A\alpha ^3+B\alpha c_2\left( W\right)
+Cc_3\left( W\right) +D\theta 
\]
in $A_{A_3\ltimes T}^{*}$. Restricting to $T$, we get 
\[
\left( C+1\right) u_1u_2u_3+D\left( u_2^2u_3+u_3^2u_1+u_1^2u_2\right) =0\in
A_{T_{SL_3}}^{*}=\frac{\mathbf{Z}\left[ u_1,u_2,u_3\right] }{\left(
u_1+u_2+u_3\right) }; 
\]
but $u_2^2u_3+u_3^2u_1+u_1^2u_2$ and $u_1u_2u_3$ are linearly independent,
hence $C=-1,$ $D=0$. Now apply the involution $\left( \cdot \right)
^\varepsilon $ to get 
\[
A\alpha ^3+B\alpha c_2\left( W\right) =0. 
\]
Since (Remark \ref{twistaction}) 
\[
\theta ^\varepsilon =-\mathrm{tsf}_{T_{SL_3}}^{A_3\ltimes T_{SL_3}}\left(
u_1^2u_3\right) , 
\]
an easy computation yields 
\[
\left( \theta -\theta ^\varepsilon +3c_3\left( W\right) \right) _{\mid
T_{SL_3}}=0. 
\]
Therefore (Proposition \ref{basico} and (i) of this Lemma) there exist $%
A,B,C\in \mathbf{Z}$ such that 
\[
\theta -\theta ^\varepsilon +3c_3\left( W\right) =A\alpha ^3+Bc_3\left(
W\right) +C\theta 
\]
is $3$-torsion. Then, restricting to $T_{SL_3}$ and observing that $%
c_3\left( W\right) _{\mid T}$ and $\theta _{\mid T}$ are linearly
independent, we get $B=C=0$; restricting now to $A_3$, we obtain $A\equiv 0$ 
$\limfunc{mod}3$ (since 
\[
\mathrm{res}_{A_3\ltimes T_{SL_3}}^{A_3}\circ \mathrm{tsf}%
_{T_{SL_3}}^{A_3\ltimes T_{SL_3}}=0). 
\]
The $C_2$-invariance of $\underline{\chi }$ is a consequence of the
transformation rules of $c_2\left( W\right) $, $c_3\left( W\right) $ and $%
\theta $.

(iv) Since $q$ is $3$-torsion, we may suppose $2$ inverted. We have $q_{\mid
T_{SL_3}}=0$ because $A_{T_{SL_3}}^{*}$ is torsion-free. It is not difficult
to verify that 
\[
\left( 2\Theta _{\mid T_{SL_3}}+3c_3\left( W\right) _{\mid T_{SL_3}}\right)
^2+4c_2\left( W\right) _{\mid T_{SL_3}}^3+27c_3\left( W\right) _{\mid
T_{SL_3}}^2=0. 
\]
Then it is enough to prove that the ideal $\mathcal{I}$ of relations between 
\[
\left\{ c_2(W)_{\mid T_{SL_3}},c_3(W)_{\mid T_{SL_3}},\Theta _{\mid
T_{SL_3}}\right\} 
\]
in $A_{T_{SL_3}}^{*}\left[ \frac 12\right] $ is generated by just this one.

Now, $\Theta _{\mid T_{SL_3}}=-\frac 32c_3(W)_{\mid T_{SL_3}}+\frac 12\delta 
$, where $\delta =\left( u_1-u_2\right) \left( u_2-u_3\right) \left(
u_1-u_3\right) $, so we have to show that 
\[
\mathcal{I}=\left( \delta ^2+4c_2(W)_{\mid T_{SL_3}}^3+27c_3(W)_{\mid
T_{SL_3}}^2\right) .
\]
Let $p\in Z\left[ \frac 12\right] \left[ X,Y,Z\right] $ with 
\[
p\left( c_2(W)_{\mid T_{SL_3}},c_3(W)_{\mid T_{SL_3}},\delta \right) =0
\]
in $A_{T_{SL_3}}^{*}\left[ \frac 12\right] $. We have 
\begin{equation}
p\left( X,Y,Z\right) =p_0\left( X,Y\right) +Zp_1\left( X,Y\right) \text{
\quad mod}\left( Z^2+4X^3+27Y^2\right) .  \label{star}
\end{equation}
If we let $C_2=\left\{ 1,\varepsilon \right\} $ act on $A_{T_{Sl_3}}^{*}%
\left[ \frac 12\right] $ permuting $u_1$ and $u_2$, we get 
\begin{eqnarray*}
\left( c_2(W)_{\mid T_{SL_3}}\right) ^\varepsilon  &=&c_2(W)_{\mid T_{SL_3}}
\\
\left( c_3(W)_{\mid T_{SL_3}}\right) ^\varepsilon  &=&c_3(W)_{\mid T_{SL_3}}
\\
\delta ^\varepsilon  &=&-\delta 
\end{eqnarray*}
and then 
\[
p^\varepsilon =p\left( c_2(W)_{\mid T_{SL_3}},c_3(W)_{\mid T_{SL_3}},-\delta
\right) =0
\]
in $A_{T_{SL_3}}^{*}\left[ \frac 12\right] $ (note that $u_1+u_2+u_3$ is $C_2
$-invariant). From (\ref{star}) we get 
\[
\left\{ 
\begin{array}{c}
p_0\left( c_2(W)_{\mid T_{SL_3}},c_3(W)_{\mid T_{SL_3}}\right) +p_1\left(
c_2(W)_{\mid T_{SL_3}},c_3(W)_{\mid T_{SL_3}}\right) \delta =0 \\ 
p_0\left( c_2(W)_{\mid T_{SL_3}},c_3(W)_{\mid T_{SL_3}}\right) -p_1\left(
c_2(W)_{\mid T_{SL_3}},c_3(W)_{\mid T_{SL_3}}\right) \delta =0
\end{array}
\right. 
\]
so ($\delta \neq 0$) 
\[
p_0\left( c_2(W)_{\mid T_{SL_3}},c_3(W)_{\mid T_{SL_3}}\right) =p_1\left(
c_2(W)_{\mid T_{SL_3}},c_3(W)_{\mid T_{SL_3}}\right) =0.
\]
But $c_2(W)_{\mid T_{SL_3}}$ and $c_3(W)_{\mid T_{SL_3}}$ are algebraically
independent, thus $p_0\left( X,Y\right) =p_1\left( X,Y\right) =0$ as
polynomials, as desired. \TeXButton{End Proof}{\endproof}

Therefore both $h^{\prime }\left( \alpha c_3(W)\right) $ and $h^{\prime
}\left( \underline{\chi }\right) $ can be identified (via lemma \ref
{transfertrick}) with their transfers, which are elements of $_3A_{\Gamma
_3}^{*}\left( Diag_{\mathrm{sl}_3}^{*}\right) $.

\begin{proposition}
\label{final}The natural morphism 
\[
f:A_{\Gamma _3}^{*}\left( Diag_{\mathrm{sl}_3}^{*}\right) \longrightarrow
\left( A_T^{*}\left( Diag_{\mathrm{sl}_3}^{*}\right) \right) ^{S_3}=\left(
A_T^{*}\right) ^{S_3}
\]
is surjective with kernel $\left( h^{\prime }\left( \alpha c_3\left(
W\right) \right) ,h^{\prime }\left( \underline{\chi }\right) \right) $,
where 
\[
h^{\prime }:A_{A_3\ltimes T}^{*}\longrightarrow A_{A_3\ltimes T}^{*}\left(
Diag_{\mathrm{sl}_3}^{*}\right) 
\]
is the pullback.
\end{proposition}

\TeXButton{Proof}{\proof} Commutativity of 
\begin{tabular}{ll}
$A_{\Gamma _3}^{*}\stackrel{g}{\quad \longrightarrow }$ & $(A_T^{*}\left(
Diag_{\mathrm{sl}_3}^{*}\right) $ \\ 
$^h\downarrow $ & $\nearrow _f$ \\ 
$A_{\Gamma _3}^{*}\left( Diag_{\mathrm{sl}_3}^{*}\right) $ & 
\end{tabular}

together with lemma \ref{hsurj}, prove that $f$ is surjective. Moreover $%
h^{\prime }\left( \alpha c_3\left( W\right) \right) $ and $h^{\prime }\left( 
\underline{\chi }\right) $ are $3$-torsion so $\ker f\supseteq \left(
h^{\prime }\left( \alpha c_3\left( W\right) \right) ,h^{\prime }\left( 
\underline{\chi }\right) \right) $ since $A_T^{*}$ is torsion-free. So we
are left to prove the reverse inclusion.

\texttt{CLAIM. }$\ker f=$ $_3A_{\Gamma _3}^{*}\left( Diag_{\mathrm{sl}%
_3}^{*}\right) \simeq $ $_3\left( A_{A_3\ltimes T}^{*}\left( Diag_{\mathrm{sl%
}_3}^{*}\right) \right) ^{C_2}$.

\underline{Proof of Claim}. $A_T^{*}$ is torsion-free, so $\ker f\supseteq $ 
$_3A_{\Gamma _3}^{*}\left( Diag_{\mathrm{sl}_3}^{*}\right) $. The pullback $%
\pi :A_{PGL_3}^{*}\rightarrow \left( A_T^{*}\right) ^{S_3}$ factors as 
\[
A_{PGL_3}^{*}\stackrel{p}{\twoheadrightarrow }A_{PGL_3}^{*}\left( U\right)
\simeq A_{\Gamma _3}^{*}\left( Diag_{\mathrm{sl}_3}^{*}\right) \stackrel{f}{%
\longrightarrow }\left( A_T^{*}\left( Diag_{\mathrm{sl}_3}^{*}\right)
\right) ^{S_3}=\left( A_T^{*}\right) ^{S_3} 
\]
and from Prop. \ref{torsion}, we get $\ker \pi =$ $_3A_{PGL_3}^{*}$; so $%
\ker \left( f\circ p\right) =$ $_3A_{PGL_3}^{*}$ and we conclude since $p$
is surjective. $\mathbf{\Box }$

Now, let $\xi \in _3\left( A_{A_3\ltimes T}^{*}\left( Diag_{\mathrm{sl}%
_3}^{*}\right) \right) ^{C_2}$. Omitting to write $h^{\prime }\left( \cdot
\right) $ everywhere and denoting $\mathrm{tsf}_T^{A_3\ltimes T}\left(
u_2^2u_3\right) $ by $\theta $, we must have 
\[
\xi =\alpha \cdot p\left( c_2(W),c_3(W)\right) +\underline{\chi }\cdot
q\left( c_2(W),c_3(W),\theta \right) 
\]
for some polynomials $p$ and $q$, since $\xi $ is $3$-torsion (we used Prop. 
\ref{basico} and lemma \ref{basico2} (i), (ii), (iv) ). But $\xi $ is also $%
C_2$-invariant, so if $C_2=\left\{ 1,\varepsilon \right\} ,$ we have: 
\[
\xi =-2\xi =-\left( \xi +\xi ^\varepsilon \right) =\alpha \cdot \left(
p^\varepsilon -p\right) +\underline{\chi }\cdot \left( -\left(
q+q^\varepsilon \right) \right) 
\]
(lemma \ref{basico2} (iii)). By lemma \ref{basico2} (iii), we have 
\[
\alpha \cdot p^\varepsilon =\alpha \cdot p\left( c_2(W),-c_3(W)\right) 
\]
and we can write $\alpha \left( p-p^\varepsilon \right) $ as a polynomial of
the form 
\[
\alpha c_3\left( W\right) \cdot p^{\prime }\left( c_2(W),c_3(W)^2\right) 
\]
for some polynomial $p^{\prime }$. Bt the \texttt{CLAIM} above, we conclude
that $\ker f\subseteq \left( \alpha c_3\left( W\right) ,\underline{\chi }%
\right) $. \TeXButton{End Proof}{\endproof}\ 

Let us summarize the situation so far. We are studying the first step (\ref
{tre}) of the stratification of $\mathrm{sl}_3$. So we started studying $%
A_{PGL_3}^{*}\left( U\right) $. We have an isomorphism $A_{PGL_3}^{*}\left(
U\right) \simeq A_{\Gamma _3}^{*}\left( Diag_{\mathrm{sl}_3}^{*}\right) $
(Prop. \ref{p3.1}) and an exact sequence (Prop. \ref{final}): 
\[
0\rightarrow \left( \alpha c_3\left( W\right) ,\underline{\chi }\right)
\longrightarrow A_{\Gamma _3}^{*}\left( Diag_{\mathrm{sl}_3}^{*}\right) 
\stackrel{f}{\longrightarrow }\left( A_T^{*}\left( Diag_{\mathrm{sl}%
_3}^{*}\right) \right) ^{S_3}=\left( A_T^{*}\right) ^{S_3}\rightarrow 0.
\]
To be precise, $\alpha c_3\left( W\right) $ and $\underline{\chi }\ $belong
to $A_{A_3\ltimes T}^{*}$ but we denote by the same symbols the elements 
\begin{eqnarray*}
&&\ \ \ \left( \mathrm{tsf}_{A_3\ltimes T}^{\Gamma _3}\left( Diag_{\mathrm{sl%
}_3}^{*}\right) \circ h^{\prime }\right) \left( \alpha c_3\left( W\right)
\right)  \\
&&\ \ \ \left( \mathrm{tsf}_{A_3\ltimes T}^{\Gamma _3}\left( Diag_{\mathrm{sl%
}_3}^{*}\right) \circ h^{\prime }\right) \left( \underline{\chi }\right) 
\end{eqnarray*}
in $A_{\Gamma _3}^{*}\left( Diag_{\mathrm{sl}_3}^{*}\right) $, where 
\[
\mathrm{tsf}_{A_3\ltimes T}^{\Gamma _3}\left( Diag_{\mathrm{sl}%
_3}^{*}\right) :A_{A_3\ltimes T}^{*}\left( Diag_{\mathrm{sl}_3}^{*}\right)
\rightarrow A_{\Gamma _3}^{*}\left( Diag_{\mathrm{sl}_3}^{*}\right) 
\]
is the transfer morphism and 
\[
h^{\prime }:A_{A_3\ltimes T}^{*}\longrightarrow A_{A_3\ltimes T}^{*}\left(
Diag_{\mathrm{sl}_3}^{*}\right) 
\]
is the obvious pullback. Moreover, by the proof of lemma \ref{hsurj} and
with the same notations, the elements 
\[
\left\{ 2c_2\left( \mathrm{sl}_3\right) -c_2\left( Sym^3E\right) ,c_3\left(
Sym^3E\right) ,c_6\left( \mathrm{sl}_3\right) \right\} \subset A_{PGL_3}^{*}
\]
project to the three generators (lemma \ref{invariants}) of $\left(
A_T^{*}\right) ^{S_3}$ through the composition 
\[
A_{PGL_3}^{*}\rightarrow A_{PGL_3}^{*}\left( U\right) \simeq A_{\Gamma
_3}^{*}\left( Diag_{\mathrm{sl}_3}^{*}\right) \stackrel{f}{\longrightarrow }%
\left( A_T^{*}\left( Diag_{\mathrm{sl}_3}^{*}\right) \right) ^{S_3}=\left(
A_T^{*}\right) ^{S_3}\text{.}
\]
If we lift the elements 
\[
\ \ \ \left( \mathrm{tsf}_{A_3\ltimes T}^{\Gamma _3}\left( Diag_{\mathrm{sl}%
_3}^{*}\right) \circ h^{\prime }\right) \left( \alpha c_3\left( W\right)
\right) ,\text{ }
\]
\[
\ \left( \mathrm{tsf}_{A_3\ltimes T}^{\Gamma _3}\left( Diag_{\mathrm{sl}%
_3}^{*}\right) \circ h^{\prime }\right) \left( \underline{\chi }\right) \in
A_{\Gamma _3}^{*}\left( Diag_{\mathrm{sl}_3}^{*}\right) ,
\]
respectively to elements $\rho ,\chi \in A_{PGL_3}^{*}$, via the surjective
pullback 
\[
A_{PGL_3}^{*}\longrightarrow A_{PGL_3}^{*}\left( U\right) \simeq A_{\Gamma
_3}^{*}\left( Diag_{\mathrm{sl}_3}^{*}\right) \text{,}
\]
we find the following $5$ generators of $A_{PGL_3}^{*}$ coming from the open
subscheme $U\subset \mathrm{sl}_3$ (through the first step (\ref{tre}) of
the stratification of $\mathrm{sl}_3$) 
\begin{equation}
\left\{ 2c_2\left( \mathrm{sl}_3\right) -c_2\left( Sym^3E\right) ,c_3\left(
Sym^3E\right) ,\rho ,\chi ,c_6\left( \mathrm{sl}_3\right) \right\} \text{,}
\label{generatorsfromU}
\end{equation}
with $\deg \rho =4$ and $\deg \chi =6$.

In the following subsection we will determine the other generators of $%
A_{PGL_3}^{*}$ coming from the complement $\mathrm{sl}_3\backslash U$,
starting from $Z_{1,1}$.

%%%%%%%%%%%%%%%%%%%%%%%%%%% End PGL_33.tex %%%%%%%%%%%%%%%%%%%%%%%%%%%%%%%
}

\QSubDoc{Include PGL_34}{%%%%%%%%%%%%%%%%%%%%%%%%%% Start PGL_34.tex %%%%%%%%%%%%%%%%%%%%%%%%%%%%%%

\LaTeXparent{PGL_3}
\ChildStyles{amssymb} 
\ChildDefaults{chapter:4,page:1}

\subsection{Generators coming from the complement of $U\subset \mathrm{sl}_3$%
}

Consider again the first step of the stratification (\ref{due}): 
\begin{equation}
A_{PGL_3}^{*}\left( Z_{1,1}\right) \stackrel{\left( j_{1,1}\right) _{*}}{%
\longrightarrow }A_{PGL_3}^{*}\left( \mathrm{sl}_3\backslash (Z_{1,0}\cup
Z_0\cup \left\{ 0\right\} )\right) \stackrel{i_{1,1}^{*}}{\longrightarrow }%
A_{PGL_3}^{*}\left( U\right) \rightarrow 0  \label{tre}
\end{equation}
where $\left( j_{1,1}\right) _{*}$ has degree $1$, equal to the codimension
of $Z_{1,1}$ in $\mathrm{sl}_3$.

\begin{lemma}
If $A\in Z_{1,1}$, let $g\in PGL_3$ be such that 
\[
g^{-1}Ag=\left( 
\begin{array}{ccc}
\lambda  & 0 & 0 \\ 
1 & \lambda  & 0 \\ 
0 & 0 & -2\lambda 
\end{array}
\right) ;
\]
then, the rule 
\[
A\longmapsto \left( \lambda ,\left[ g\right] \right) 
\]
defines a $PGL_3$-equivariant isomorphism $Z_{1,1}\rightarrow \mathbf{A}%
^1\backslash \left\{ 0\right\} \times \frac{PGL_3}{\mathrm{U}_2\times 
\mathbf{G}_m}$, where $\mathrm{U}_2$ is the full unipotent subgroup of $GL_2$
and $PGL_3$ acts trivially on $\mathbf{A}^1\backslash \left\{ 0\right\} $.
\end{lemma}

\TeXButton{Proof}{\proof} Everything is a straightforward verification left
to the interested reader. We only note that the stabilizer of 
\[
\left( 
\begin{array}{ccc}
\lambda & 0 & 0 \\ 
1 & \lambda & 0 \\ 
0 & 0 & -2\lambda
\end{array}
\right) 
\]
(under the adjoint action of $PGL_3$) is 
\[
\left\{ \left[ g\right] \mid g=\left( 
\begin{array}{ccc}
\alpha & 0 & 0 \\ 
\beta & \alpha & 0 \\ 
0 & 0 & \gamma
\end{array}
\right) \text{, }\alpha ,\gamma \in \mathbf{G}_m\right\} 
\]
which is obviously isomorphic to $\mathrm{U}_2\times \mathbf{G}_m$. 
\TeXButton{End Proof}{\endproof}

By Corollary \ref{corunip}, Prop. \ref{productbygm}, \ref{homogeneous} and
lemma \ref{utrivial}, we have 
\begin{equation}
A_{PGL_3}^{*}\left( Z_{1,1}\right) \simeq A_{\mathbf{G}_m}^{*}=\mathbf{Z}%
\left[ u\right] .  \label{z11}
\end{equation}
It is not difficult to verify that 
\[
j_{1,1}^{*}\left( 2c_2\left( \mathrm{sl}_3\right) -c_2\left( Sym^3E\right)
\right) =u^2\text{ ,}
\]
where we abused notation writing $2c_2\left( \mathrm{sl}_3\right) -c_2\left(
Sym^3E\right) $ for its pullback to $A_{PGL_3}^{*}\left( \mathrm{sl}%
_3\backslash (Z_{1,0}\cup Z_0\cup \left\{ 0\right\} )\right) $. Moreover 
\[
\left( j_{1,1}\right) _{*}\left( 1\right) =\left[ Z_{1,1}\right]
=D^{*}\left( \left[ \left\{ 0\right\} \right] \right) =0
\]
where $D:\mathrm{sl}_3\backslash (Z_{1,0}\cup Z_0\cup \left\{ 0\right\}
)\rightarrow \mathbf{A}^1$ is the discriminant; so, by projection formula,
the ideal $\mathrm{im}\left( j_{1,1}\right) _{*}$ is generated by $\left(
j_{1,1}\right) _{*}\left( u\right) $.

Let $\Theta _{1,1}^{\left( 2\right) }$ be a lift of $\left( j_{1,1}\right)
_{*}\left( u\right) \in A_{PGL_3}^{*}\left( \mathrm{sl}_3\backslash
(Z_{1,0}\cup Z_0\cup \left\{ 0\right\} )\right) $ to $A_{PGL_3}^{*}$. The
analysis we made of (\ref{tre}) has the following upshot (recall (\ref
{generatorsfromU})): $A_{PGL_3}^{*}\left( \mathrm{sl}_3\backslash
(Z_{1,0}\cup Z_0\cup \left\{ 0\right\} )\right) $\emph{\ is generated by
(the images via }$A_{PGL_3}^{*}\twoheadrightarrow A_{PGL_3}^{*}\left( 
\mathrm{sl}_3\backslash (Z_{1,0}\cup Z_0\cup \left\{ 0\right\} )\right) $%
\emph{\ of) } 
\begin{equation}
\left\{ 2c_2\left( \mathrm{sl}_3\right) -c_2\left( Sym^3E\right) ,\Theta
_{1,1}^{\left( 2\right) },c_3\left( Sym^3E\right) ,\rho ,\chi ,c_6\left( 
\mathrm{sl}_3\right) \right\} \text{.}  \label{firststep}
\end{equation}
\ 

Now let us proceed one step further in the analysis of stratification (\ref
{due}); the second exact sequence of (\ref{uno}) is: 
\begin{equation}
A_{PGL_3}^{*}\left( Z_{1,0}\right) \stackrel{\left( j_{1,0}\right) _{*}}{%
\longrightarrow }A_{PGL_3}^{*}\left( \mathrm{sl}_3\backslash (Z_0\cup
\left\{ 0\right\} )\right) \stackrel{i_{1,0}^{*}}{\longrightarrow }%
A_{PGL_3}^{*}\left( \mathrm{sl}_3\backslash (Z_{1,0}\cup Z_0\cup \left\{
0\right\} )\right) \rightarrow 0  \label{quattro}
\end{equation}
where $\left( j_{1,0}\right) _{*}$ has degree $3$, equal to the codimension
of $Z_{1,0}$ in $\mathrm{sl}_3$. We omit the straightforward proof of the
following:

\begin{lemma}
If $A\in Z_{1,0}$, let $g\in PGL_3$ be such that 
\[
g^{-1}Ag=\left( 
\begin{array}{ccc}
\lambda  & 0 & 0 \\ 
0 & \lambda  & 0 \\ 
0 & 0 & -2\lambda 
\end{array}
\right) .
\]
Then, the rule 
\[
A\longmapsto \left( \lambda ,\left[ g\right] \right) 
\]
defines a $PGL_3$-equivariant isomorphism $Z_{1,0}\rightarrow \mathbf{A}%
^1\backslash \left\{ 0\right\} \times \frac{PGL_3}{GL_2}$, where $GL_2$
injects as 
\[
\left( 
\begin{array}{cc}
GL_2 & 0 \\ 
0 & 1
\end{array}
\right) 
\]
and $PGL_3$ acts trivially on $\mathbf{A}^1\backslash \left\{ 0\right\} $.
\end{lemma}

Then, by Prop. \ref{homogeneous} and lemma \ref{utrivial}, we have\footnote{$%
\lambda _i\doteq c_i\left( \text{\textrm{standard representation}}\right) $.}
\begin{equation}
A_{PGL_3}^{*}\left( Z_{1,0}\right) \simeq A_{GL_2}^{*}=\mathbf{Z}\left[
\lambda _1,\lambda _2\right] \text{.}  \label{az10}
\end{equation}

\begin{lemma}
\label{diga}$\left( j_{1,0}\right) _{*}$ is $3$-torsion.
\end{lemma}

\TeXButton{Proof}{\proof}If $\xi \in A_{PGL_3}^{*}\left( Z_{1,0}\right) $,
let $\widehat{\xi }\in A_{PGL_3}^{*}$ be a lift of $\left( j_{1,0}\right)
_{*}\left( \xi \right) $ via the surjective pullback 
\[
\pi _{1,0}:A_{PGL_3}^{*}\longrightarrow A_{PGL_3}^{*}\left( sl_3\backslash
Z_0\cup \left\{ 0\right\} \right) . 
\]
It is enough to prove that $\widehat{\xi }$ is $3$-torsion i.e. that 
\[
\widehat{\xi }\in \ker \left( A_{PGL_3}^{*}\longrightarrow \left(
A_T^{*}\right) ^{S_3}\right) , 
\]
since by \cite{EG1}, Prop. 6, the rational pullback 
\[
A_{PGL_3}^{*}\otimes \mathbf{Q}\longrightarrow \left( A_T^{*}\right)
^{S_3}\otimes \mathbf{Q} 
\]
is an isomorphism and $A_{PGL_3}^{*}$ has only $3$-torsion by Cor. \ref
{cortorsion}.

Now, observe that 
\[
\left( j_{1,0}\right) _{*}\left( \xi \right) \in \ker \left(
A_{PGL_3}^{*}\left( \mathrm{sl}_3\backslash Z_0\cup \left\{ 0\right\}
\right) \rightarrow A_{PGL_3}^{*}\left( \mathrm{sl}_3\backslash Z_{1,0}\cup
Z_0\cup \left\{ 0\right\} \right) \right) 
\]
by the obvious localization sequence and therefore 
\[
\widehat{\xi }\in \ker \left( A_{PGL_3}^{*}\longrightarrow
A_{PGL_3}^{*}\left( U\right) \right) , 
\]
by (\ref{tre}). To conclude, we note that $A_{PGL_3}^{*}\longrightarrow
\left( A_T^{*}\right) ^{S_3}$ factors as 
\[
A_{PGL_3}^{*}\longrightarrow A_{PGL_3}^{*}\left( U\right) \simeq A_{\Gamma
_3}^{*}\left( Diag_{\mathrm{sl}_3}^{*}\right) \longrightarrow \left(
A_T^{*}\right) ^{S_3}\text{.} 
\]

\TeXButton{End Proof}{\endproof}

\begin{proposition}
The ideal $\mathrm{im}\left( j_{1,0}\right) _{*}$ is generated by 
\[
\left\{ \left( j_{1,0}\right) _{*}\left( 1\right) ,\left( j_{1,0}\right)
_{*}\left( \lambda _1\right) ,\left( j_{1,0}\right) _{*}\left( \lambda
_2\right) ,\left( j_{1,0}\right) _{*}\left( \lambda _2^2\right) ,\left(
j_{1,0}\right) _{*}\left( \lambda _1\lambda _2\right) ,\left( j_{1,0}\right)
_{*}\left( \lambda _1\lambda _2^2\right) \right\} . 
\]
\end{proposition}

\TeXButton{Proof}{\proof}Identifying $A_{PGL_3}^{*}\left( Z_{1,0}\right) \ $%
with $A_{GL_2}^{*}=\mathbf{Z}\left[ \lambda _1,\lambda _2\right] $ via (\ref
{az10}) and writing $\left( \cdot \right) _{\mid GL_2}$ for $j_{1,0}^{*}$,
one can easily verify that 
\begin{equation}
\left( \lambda \equiv 2c_2\left( \mathrm{sl}_3\right) -c_2\left(
Sym^3E\right) \right) _{\mid GL_2}=\lambda _1^2-3\lambda _2\doteq \tau _2
\label{alef}
\end{equation}

\begin{equation}
c_6\left( \mathrm{sl}_3\right) _{\mid GL_2}=-\lambda _1^2\lambda
_2^2+4\lambda _2^3\doteq \tau _6.  \label{alef1}
\end{equation}
Therefore 
\begin{equation}
\lambda _2^3=\tau _6+\tau _2\lambda _2^2,  \label{alef3}
\end{equation}
and if $\frak{I}$ denotes the ideal generated by 
\[
\begin{array}{c}
\left\{ \left( j_{1,0}\right) _{*}\left( 1\right) ,\left( j_{1,0}\right)
_{*}\left( \lambda _1\right) ,\left( j_{1,0}\right) _{*}\left( \lambda
_2\right) ,\right. \\ 
\left. \left( j_{1,0}\right) _{*}\left( \lambda _2^2\right) ,\left(
j_{1,0}\right) _{*}\left( \lambda _1\lambda _2\right) ,\left( j_{1,0}\right)
_{*}\left( \lambda _1\lambda _2^2\right) \right\} ,
\end{array}
\]
we have 
\begin{equation}
\left( j_{1,0}\right) _{*}\left( \lambda _2^m\right) \in \frak{I},\text{ }%
\forall m\geq 0,  \label{alef2}
\end{equation}
by an easy induction on $m$, using projection formula.

Now, consider the general monomial $\lambda _1^n\lambda _2^m$. If $n=2r$, we
have 
\begin{equation}
\lambda _1^n\lambda _2^m=\left( \tau _2+3\lambda _2\right) ^r\lambda
_2^m\equiv \tau _2^r\lambda _2^m\text{ }\left( \limfunc{mod}3\right)
\label{alef4}
\end{equation}
and then 
\[
\left( j_{1,0}\right) _{*}\left( \lambda _1^{2r}\lambda _2^m\right) =\lambda
^r\cdot \left( j_{1,0}\right) _{*}\left( \lambda _2^m\right) , 
\]
by lemma \ref{diga} and projection formula; thus 
\[
\left( j_{1,0}\right) _{*}\left( \lambda _1^{2r}\lambda _2^m\right) \in 
\frak{I},\text{ }\forall m,r\geq 0, 
\]
by (\ref{alef2}). If $n=2r+1$, (\ref{alef4}), (\ref{alef3}) and projection
formula easily reduce the assert 
\[
\left( j_{1,0}\right) _{*}\left( \lambda _1^{2r+1}\lambda _2^m\right) \in 
\frak{I},\text{ }\forall m,r\geq 0 
\]
to the assert 
\[
\left( j_{1,0}\right) _{*}\left( \lambda _1\lambda _2^m\right) \in \frak{I},%
\text{ }\forall m\geq 0, 
\]
which is easily proved by induction on $m$.

Since the monomials $\lambda _1^n\lambda _2^m$ generates $A_{GL_2}^{*}$ as a 
$\mathbf{Z}$-module, we conclude that 
\[
\frak{I}=\mathrm{im}\left( j_{1,0}\right) _{*}\text{ .} 
\]
\TeXButton{End Proof}{\endproof}

Therefore, if we denote by $\Theta _{1,0}^{\left( 3\right) }$ (respectively, 
$\Theta _{1,0}^{\left( 4\right) }$, $\Theta _{1,0}^{\left( 5\right) }$, $%
\Theta _{1,0}^{\left( 6\right) }$, $\Theta _{1,0}^{\left( 7\right) }$, $%
\Theta _{1,0}^{\left( 8\right) }$) a lift of $\left( j_{1,0}\right)
_{*}\left( 1\right) $ (respectively, of $\left( j_{1,0}\right) _{*}\left(
\lambda _1\right) $, $\left( j_{1,0}\right) _{*}\left( \lambda _2\right) $, $%
\left( j_{1,0}\right) _{*}\left( \lambda _1\lambda _2\right) $, $\left(
j_{1,0}\right) _{*}\left( \lambda _2^2\right) $, $\left( j_{1,0}\right)
_{*}\left( \lambda _1\lambda _2^2\right) $) to $A_{PGL_3}^{*}$, from (\ref
{firststep}) and (\ref{quattro}) we get that $A_{PGL_3}^{*}\left( \mathrm{sl}%
_3\backslash (Z_0\cup \left\{ 0\right\} )\right) $\emph{\ is generated by
(the images via }$A_{PGL_3}^{*}\twoheadrightarrow A_{PGL_3}^{*}\left( 
\mathrm{sl}_3\backslash (Z_0\cup \left\{ 0\right\} )\right) $\emph{\ of) } 
\begin{equation}
\begin{array}{c}
\left\{ 2c_2\left( \mathrm{sl}_3\right) -c_2\left( Sym^3E\right) ,\Theta
_{1,1}^{\left( 2\right) },c_3\left( Sym^3E\right) ,\Theta _{1,0}^{\left(
3\right) },\right. \\ 
\left. \rho ,\Theta _{1,0}^{\left( 4\right) },\Theta _{1,0}^{\left( 5\right)
},\chi ,\Theta _{1,0}^{\left( 6\right) },c_6\left( \mathrm{sl}_3\right)
,\Theta _{1,0}^{\left( 7\right) },\Theta _{1,0}^{\left( 8\right) }\right\} .
\end{array}
\label{secondstep}
\end{equation}

Let us proceed one step further in the analysis of stratification (\ref{due}%
); the third exact sequence of (\ref{uno}), in our case is: 
\begin{equation}
A_{PGL_3}^{*}\left( Z_{0,1}\right) \stackrel{\left( j_{0,1}\right) _{*}}{%
\longrightarrow }A_{PGL_3}^{*}\left( \mathrm{sl}_3\backslash (Z_{0,0}\cup
\left\{ 0\right\} )\right) \stackrel{i_{0,1}^{*}}{\longrightarrow }%
A_{PGL_3}^{*}\left( \mathrm{sl}_3\backslash (Z_0\cup \left\{ 0\right\}
)\right) \rightarrow 0  \label{cinque}
\end{equation}
where $\left( j_{0,1}\right) _{*}$ has degree $2$, equal to the codimension
of $Z_{0,1}$ in $\mathrm{sl}_3$. $Z_{0,1}$ is a $PGL_3$-orbit with
stabilizer 
\[
\left\{ \left[ g\right] \mid g=\left( 
\begin{array}{ccc}
1 & 0 & 0 \\ 
\alpha & 1 & 0 \\ 
\beta & \alpha & 1
\end{array}
\right) ,\alpha ,\beta \in \mathbf{C}\right\} 
\]
which is unipotent and then, by Prop. \ref{unipotent}, we have $%
A_{PGL_3}^{*}\left( Z_{0,1}\right) =\mathbf{Z}$. If we denote by $\Theta
_{0,1}^{\left( 2\right) }$ a lift of $\left( j_{0,1}\right) _{*}\left(
1\right) =\left[ Z_{0,1}\right] \in A_{PGL_3}^2\left( \mathrm{sl}%
_3\backslash (Z_{0,0}\cup \left\{ 0\right\} )\right) $ via the surjective
pullback 
\[
A_{PGL_3}^{*}\longrightarrow A_{PGL_3}^{*}\left( \mathrm{sl}_3\backslash
(Z_{0,0}\cup \left\{ 0\right\} )\right) , 
\]
from (\ref{cinque}) and (\ref{secondstep}) we get that $A_{PGL_3}^{*}\left( 
\mathrm{sl}_3\backslash (Z_{0,0}\cup \left\{ 0\right\} )\right) $\emph{\ is
generated by (the images via }$A_{PGL_3}^{*}\twoheadrightarrow
A_{PGL_3}^{*}\left( \mathrm{sl}_3\backslash (Z_{0,0}\cup \left\{ 0\right\}
)\right) $\emph{\ of) } 
\begin{equation}
\begin{array}{c}
\left\{ 2c_2\left( \mathrm{sl}_3\right) -c_2\left( Sym^3E\right) ,\Theta
_{1,1}^{\left( 2\right) },\Theta _{0,1}^{\left( 2\right) },c_3\left(
Sym^3E\right) ,\Theta _{1,0}^{\left( 3\right) },\right. \\ 
\left. \rho ,\Theta _{1,0}^{\left( 4\right) },\Theta _{1,0}^{\left( 5\right)
},\chi ,\Theta _{1,0}^{\left( 6\right) },c_6\left( \mathrm{sl}_3\right)
,\Theta _{1,0}^{\left( 7\right) },\Theta _{1,0}^{\left( 8\right) }\right\} .
\end{array}
\label{secondbisstep}
\end{equation}

We have come to the second-last step of stratification (\ref{due}): 
\begin{equation}
A_{PGL_3}^{*}\left( Z_{0,0}\right) \stackrel{\left( j_{0,0}\right) _{*}}{%
\longrightarrow }A_{PGL_3}^{*}\left( \mathrm{sl}_3\backslash \left\{
0\right\} \right) \stackrel{i_{0,0}^{*}}{\longrightarrow }%
A_{PGL_3}^{*}\left( \mathrm{sl}_3\backslash (Z_{0,0}\cup \left\{ 0\right\}
)\right) \rightarrow 0  \label{sei}
\end{equation}
where $\left( j_{0,0}\right) _{*}$ has degree $4$, equal to the codimension
of $Z_{0,0}$ in $\mathrm{sl}_3$. $Z_{0,0}$ is a $PGL_3$-orbit with
stabilizer 
\[
\left\{ \left[ g\right] \mid g=\left( 
\begin{array}{ccc}
1 & 0 & 0 \\ 
\alpha & \delta & 0 \\ 
\beta & \gamma & 1
\end{array}
\right) ,\alpha ,\beta ,\gamma \in \mathbf{C,}\text{ }\delta \in \mathbf{G}%
_m\right\} 
\]
which is a split extension of $\mathbf{G}_m$ by the full unipotent group $%
\mathrm{U}_3\subset GL_3$. By Cor. \ref{corunip} we get an isomorphism $%
A_{PGL_3}^{*}\left( Z_{0,0}\right) \simeq A_{\mathbf{G}_m}^{*}=\mathbf{Z}%
\left[ u\right] $. Since 
\[
j_{0,0}^{*}\left( 2c_2\left( sl_3\right) -c_2\left( Sym^3E\right) \right)
=u^2, 
\]
$A_{PGL_3}^{*}\left( Z_{0,0}\right) $ is generated by $\left\{ \left(
j_{0,0}\right) _{*}\left( 1\right) ,\left( j_{0,0}\right) _{*}\left(
u\right) \right\} $ as an $A_{PGL_3}^{*}\left( \mathrm{sl}_3\backslash
\left\{ 0\right\} \right) $-module (by projection formula) and if we denote
by $\Theta _{0,0}^{\left( 4\right) }$ (respectively, $\Theta _{0,0}^{\left(
5\right) }$) a lift of $\left( j_{0,0}\right) _{*}\left( 1\right) $
(respectively, of $\left( j_{0,0}\right) _{*}\left( u\right) $) to $%
A_{PGL_3}^{*}$, we get that $A_{PGL_3}^{*}\left( \mathrm{sl}_3\backslash
\left\{ 0\right\} \right) $\emph{\ is generated by (the images via }$%
A_{PGL_3}^{*}\twoheadrightarrow A_{PGL_3}^{*}\left( \mathrm{sl}_3\backslash
\left\{ 0\right\} \right) $\emph{\ of) } 
\begin{equation}
\begin{array}{c}
\left\{ 2c_2\left( \mathrm{sl}_3\right) -c_2\left( Sym^3E\right) ,\Theta
_{1,1}^{\left( 2\right) },\Theta _{0,1}^{\left( 2\right) },c_3\left(
Sym^3E\right) ,\Theta _{1,0}^{\left( 3\right) },\right. \\ 
\left. \rho ,\Theta _{1,0}^{\left( 4\right) },\Theta _{0,0}^{(4)},\Theta
_{1,0}^{\left( 5\right) },\Theta _{0,0}^{(5)},\chi ,\Theta _{1,0}^{\left(
6\right) },c_6\left( \mathrm{sl}_3\right) ,\Theta _{1,0}^{\left( 7\right)
},\Theta _{1,0}^{\left( 8\right) }\right\} .
\end{array}
\label{thirdstep}
\end{equation}

The last step of (\ref{uno}) for stratification (\ref{due}) is immediate
because 
\[
A_{PGL_3}^{*}\left( \mathrm{sl}_3\backslash \left\{ 0\right\} \right) \simeq
A_{PGL_3}^{*}\diagup \left( c_8\left( \mathrm{sl}_3\right) \right) 
\]
by self-intersection formula (\cite{Fu}, p. 103).

Therefore we conclude our analysis of the stratification (\ref{due}) with
the following result:

\begin{proposition}
\label{prelgen}$A_{PGL_3}^{*}$\ is generated by 
\begin{equation}
\begin{array}{c}
\left\{ 2c_2\left( \mathrm{sl}_3\right) -c_2\left( Sym^3E\right) ,\Theta
_{1,1}^{\left( 2\right) },\Theta _{0,1}^{\left( 2\right) },c_3\left(
Sym^3E\right) ,\Theta _{1,0}^{\left( 3\right) },\right. \\ 
\left. \rho ,\Theta _{1,0}^{\left( 4\right) },\Theta _{0,0}^{(4)},\Theta
_{1,0}^{\left( 5\right) },\Theta _{0,0}^{(5)},\chi ,\Theta _{1,0}^{\left(
6\right) },c_6\left( \mathrm{sl}_3\right) ,\Theta _{1,0}^{\left( 7\right)
},\Theta _{1,0}^{\left( 8\right) },c_8\left( \mathrm{sl}_3\right) \right\} ,
\end{array}
\label{laststep}
\end{equation}
where $\deg \Theta _{0,1}^{\left( 2\right) }=2,$ $\deg \Theta _{1,1}^{\left(
2\right) }=2$, $\deg \rho =4$, $\deg \chi =6$, $\deg \Theta _{1,0}^{\left(
m\right) }=m$\ and $\deg \Theta _{0,0}^{(r)}=r.$
\end{proposition}

We will make this result more precise in the following section by getting
rid of all the $\Theta $ generators.

\section{$A_{PGL_3}^{*}$ is not generated by Chern classes. Elimination of
some generators}

In this section we first prove that $A_{PGL_3}^{*}$ is not generated by
Chern classes and then that all its $\Theta $ generators are zero.

\begin{lemma}
\label{cohom}Writing $H_{PGL_3}^i$ for $H^i\left( \mathrm{B}PGL_3,\mathbf{Z}%
\right) $, we have: 
\begin{equation}
\begin{tabular}{|ll|}
\hline
\multicolumn{1}{|l|}{$H_{PGL_3}^0\simeq \mathbf{Z}$} & $H_{PGL_3}^1=0$ \\ 
\hline
\multicolumn{1}{|l|}{$H_{PGL_3}^2=0$} & $H_{PGL_3}^3\simeq \mathbf{Z}/3$ \\ 
\hline
\multicolumn{1}{|l|}{$H_{PGL_3}^4\simeq \mathbf{Z}$} & $H_{PGL_3}^5=0$ \\ 
\hline
\multicolumn{1}{|l|}{$H_{PGL_3}^6\simeq \mathbf{Z}$} & $H_{PGL_3}^7=0$ \\ 
\hline
\multicolumn{1}{|l|}{$H_{PGL_3}^8\simeq \mathbf{Z\oplus Z}/3$} & $%
H_{PGL_3}^9=0$ \\ \hline
\multicolumn{1}{|l|}{$H_{PGL_3}^{10}\simeq \mathbf{Z}$} & $%
H_{PGL_3}^{11}\simeq \mathbf{Z}/3$ \\ \hline
\multicolumn{1}{|l|}{$H_{PGL_3}^{12}\simeq \mathbf{Z\oplus Z}$} & $%
H_{PGL_3}^{13}=0$ \\ \hline
\multicolumn{1}{|l|}{$H_{PGL_3}^{14}\simeq \mathbf{Z}$} & $%
H_{PGL_3}^{15}\simeq \mathbf{Z}/3$ \\ \hline
\multicolumn{1}{|l|}{$H_{PGL_3}^{16}\simeq \mathbf{Z\oplus Z}\oplus \mathbf{Z%
}/3$} & $\cdots $ \\ \hline
\end{tabular}
\label{Hstar}
\end{equation}
\end{lemma}

\TeXButton{Proof}{\proof} It is a routine computation using the Universal
Coefficients' Formula for cohomology (e.g. \cite{Sp}), once one knows the
following facts:

\begin{enumerate}
\item[1.]  $H^{*}\left( \mathrm{B}PGL_3,\mathbf{Q}\right) \simeq H^{*}\left( 
\mathrm{B}SL_3,\mathbf{Q}\right) =\mathbf{Q}\left[ c_2\left( E\right)
,c_3\left( E\right) \right] $, $E$ being the standard representation of $SL_3
$;

\item[2.]  $H^{*}\left( \mathrm{B}PGL_3,\mathbf{Z}\right) $ has only $3$%
-torsion;

\item[3.]  there is a ring isomorphism 
\[
H^{*}\left( \mathrm{B}PGL_3,\mathbf{Z}/3\right) \simeq \mathbf{Z}/3\left[
y_2,y_8,y_{12}\right] \otimes \Lambda \left( y_3,y_7\right) \diagup \left(
y_2y_3,y_2y_7,y_2y_8+y_3y_7\right) 
\]
where $\deg y_i=i$.
\end{enumerate}

1. follows immediately from the Leray spectral sequence 
\[
H^p\left( \mathrm{B}PGL_3,H^q\left( \mathrm{B}\mathbf{\mu }_3,\mathbf{Z}%
\right) \right) \Longrightarrow H^{p+q}\left( \mathrm{B}SL_3,\mathbf{Z}%
\right) ; 
\]
2. is proved in \cite{KY}, p. 790 and 3. was computed in \cite{KMS}. 
\TeXButton{End Proof}{\endproof}

\begin{theorem}
\label{result1}$A_{PGL_3}^{*}$ is not generated by Chern classes; more
precisely, $\rho $ is not a polynomial in Chern classes.
\end{theorem}

\TeXButton{Proof}{\proof} We proceed in 4 steps:

\begin{enumerate}
\item[(I)]  First we show that $\mathrm{cl}\left( \rho \right) $ is nonzero
in $H^8(\mathrm{B}PGL_3,\mathbf{Z})_{tors}$, where $\mathrm{cl}%
:A_{PGL_3}^{*}\rightarrow H^{*}(\mathrm{B}PGL_3,\mathbf{Z})$ is the cycle
class map and $\rho $ is one of the generators of $A_{PGL_3}^{*}$ (see Prop. 
\ref{prelgen});

\item[(II)]  Then we use a spectral sequence argument to show that 
\[
\mathrm{im}(H^8(\mathrm{B}PGL_3,\mathbf{Z})\rightarrow H^8(\mathrm{B}SL_3,%
\mathbf{Z}))
\]
has index at least $9$ in $H^8(\mathrm{B}SL_3,\mathbf{Z})\simeq \mathbf{Z}$;

\item[(III)]  Next, we use the fact that $c_2\left( \mathrm{sl}_3\right)
^2\mapsto 36\alpha _2^2$ via 
\[
H^8(\mathrm{B}PGL_3,\mathbf{Z})\simeq \mathbf{Z}\oplus \mathbf{Z}%
/3\longrightarrow H^8(\mathrm{B}SL_3,\mathbf{Z})\simeq \mathbf{Z\cdot }%
\alpha _2^2
\]
to conclude that 
\[
H^8(\mathrm{B}PGL_3,\mathbf{Z}_{(3)})\simeq \mathbf{Z}_{(3)}\cdot c_2\left( 
\mathrm{sl}_3\right) ^2\oplus \mathbf{Z}/3\cdot \mathrm{cl}(\rho )
\]
(where we have written $\alpha $ in place of $j_{\bullet }\left( \alpha
\right) $, with $j_{\bullet }:H^{*}(\mathrm{B}PGL_3,\mathbf{Z})\rightarrow
H^{*}(\mathrm{B}PGL_3,\mathbf{Z}_{(3)})$ induced by the localization $j:%
\mathbf{Z}\rightarrow \mathbf{Z}_{(3)}$ and $\alpha _i\doteq c_i\left( \text{%
standard repr. of }SL_3\right) $);

\item[(IV)]  Finally we show that $\mathrm{cl}(\rho )\in H^8(\mathrm{B}PGL_3,%
\mathbf{Z})$ is not in the Chern subring of $H^{*}(\mathrm{B}PGL_3,\mathbf{Z}%
)$ (implying that $\rho $ itself is not in the Chern subring of $%
A_{PGL_3}^{*}$).
\end{enumerate}

(I) We freely use Remark \ref{twistaction}. Recall that $\rho $ is a lift to 
$A_{PGL_3}^{*}$ of 
\[
\left( \alpha c_3\left( W\right) \right) _{\mid Diag_{\mathrm{sl}_3}^{*}}\in 
\text{ }_3\left( A_{A_3\ltimes T}^{*}\left( Diag_{\mathrm{sl}_3}^{*}\right)
\right) ^{C_2}\simeq \text{ }_3A_{\Gamma _3}^{*}\left( Diag_{\mathrm{sl}%
_3}^{*}\right) \simeq \text{ }_3A_{PGL_3}^{*}\left( U\right) , 
\]
with $\alpha c_3\left( W\right) \in $ $_3A_{A_3\ltimes T}^{*}$. To prove (I)
it is then enough to show that 
\[
\mathrm{cl}\left( \alpha c_3\left( W\right) \right) _{\mid Diag_{\mathrm{sl}%
_3}^{*}}\neq 0\text{ in }H_{A_3\ltimes T}^8\left( Diag_{\mathrm{sl}_3}^{*},%
\mathbf{Z}\right) \text{.} 
\]
If $A_3\times \mathbf{\mu }_3\hookrightarrow A_3\ltimes T$, we will get this
by showing 
\begin{equation}
\mathrm{cl}\left( \alpha c_3\left( W\right) \right) _{\mid Diag_{\mathrm{sl}%
_3}^{*}}\neq 0\text{ in }H_{A_3\times \mathbf{\mu }_3}^8\left( Diag_{\mathrm{%
sl}_3}^{*},\mathbf{Z}\right)  \label{pippo}
\end{equation}
(writing again $\mathrm{cl}\left( \alpha c_3\left( W\right) \right) $ for
its the restriction to $H_{A_3\times \mathbf{\mu }_3}^8$). Let us consider
the localization exact sequences for cohomology, corresponding to $Diag_{%
\mathrm{sl}_3}\supset Diag_{\mathrm{sl}_3}\backslash \left\{ 0\right\}
\supset Diag_{\mathrm{sl}_3}^{*}$ 
\begin{equation}
H_{A_3\times \mathbf{\mu }_3}^4\stackrel{\cdot (-\alpha ^2)}{\longrightarrow 
}H_{A_3\times \mathbf{\mu }_3}^8\left( Diag_{\mathrm{sl}_3},\mathbf{Z}%
\right) \stackrel{p}{\longrightarrow }H_{A_3\times \mathbf{\mu }_3}^8\left(
Diag_{\mathrm{sl}_3}\smallsetminus \left\{ 0\right\} ,\mathbf{Z}\right)
\label{cinque.uno}
\end{equation}
\begin{equation}
H_{\mathbf{\mu }_3}^6\simeq H_{A_3\times \mathbf{\mu }_3}^6\stackrel{}{%
\left( Z,\mathbf{Z}\right) \stackrel{i_{*}}{\longrightarrow }}H_{A_3\times 
\mathbf{\mu }_3}^8\left( Diag_{\mathrm{sl}_3}\smallsetminus \left\{
0\right\} ,\mathbf{Z}\right) \stackrel{q}{\longrightarrow }H_{A_3\times 
\mathbf{\mu }_3}^8\left( Diag_{\mathrm{sl}_3}^{*},\mathbf{Z}\right)
\label{cinque.due}
\end{equation}
where $i:Z\doteq \left( Diag_{\mathrm{sl}_3}\backslash \left\{ 0\right\}
\right) \backslash Diag_{\mathrm{sl}_3}^{*}\hookrightarrow Diag_{\mathrm{sl}%
_3}\backslash \left\{ 0\right\} $ and we used that $Z\simeq A_3\times 
\mathbf{C}^{*}$, $A_3\ltimes T$-equivariantly. If $\mathbf{C}_{\chi ,\mathbf{%
\mu }_3}$ (respectively, $\mathbf{C}_{perm,A_3}^3$) denotes the $\mathbf{\mu 
}_3$-representation given by multiplication by the character $\chi =\exp
\left( i2\pi /3\right) $ (respectively, the $A_3$-permutation
representation), we have $W\simeq \mathbf{C}_{\chi ,\mathbf{\mu }%
_3}\boxtimes \mathbf{C}_{perm,A_3}^3$ as $A_3\times \mathbf{\mu }_3$%
-representations. Then, if we let $H_{\mathbf{\mu }_3}^{*}=\mathbf{Z}\left[
\beta \right] /\left( 3\beta \right) $, $H_{A_3}^{*}=\mathbf{Z}\left[ \alpha
\right] /\left( 3\alpha \right) $, the Chern roots of $W$ are $\left\{ \beta
+\alpha ,\beta -\alpha ,\beta \right\} $ and 
\[
\mathrm{cl}\left( \alpha c_3\left( W\right) \right) =\left( \beta ^2-\alpha
^2\right) \alpha \beta \in H_{A_3\times \mathbf{\mu }_3}^8\text{.} 
\]
Now we claim $i_{*}=0$ in (\ref{cinque.due}). In fact, consider the pullback 
$E$ of $\mathbf{C}_{\chi ,\mathbf{\mu }_3}$ to $Diag_{\mathrm{sl}%
_3}\backslash \left\{ 0\right\} $ and view $E$ as an $A_3\times \mathbf{\mu }%
_3$-equivariant vector bundle on $Diag_{\mathrm{sl}_3}\backslash \left\{
0\right\} $, with $A_3$ acting trivially on $E$. Obviously, $i^{*}\left(
c_1\left( E\right) \right) =\beta $. But we also have $i_{*}\left( 1\right)
=0$ since 
\[
Z=D^{-1}\left( \left\{ 0\right\} \right) , 
\]
where 
\begin{eqnarray}
D &:&Diag_{\mathrm{sl}_3}\backslash \left\{ 0\right\} \rightarrow \mathbf{A}%
^1  \nonumber \\
\left( \lambda _1,\lambda _2,\lambda _3\right) &\longmapsto &\left( \lambda
_1-\lambda _2\right) \left( \lambda _1-\lambda _3\right) \left( \lambda
_2-\lambda _3\right)
\end{eqnarray}
is the square root of the discriminant (which is $A_3\times \mathbf{\mu }_3$%
-equivariant!). By projection formula, $i_{*}=0$ and $q$ is injective.

So, we are left to show that $p\left( \left( \beta ^2-\alpha ^2\right)
\alpha \beta \right) =p\left( \alpha \beta ^3\right) \neq 0$ in (\ref
{cinque.uno}). Now observe that 
\[
H_{A_3\times \mathbf{\mu }_3}^{2n}\simeq \left( H_{A_3}^{*}\otimes H_{%
\mathbf{\mu }_3}^{*}\right) ^{2n} 
\]
by K\"unneth formula, since 
\[
\bigoplus_{p+q=2n+1}Tor_1^{\mathbf{Z}}\left( H_{A_3}^p,H_{\mathbf{\mu }%
_3}^q\right) =0 
\]
(either $p$ or $q$ being odd in every summand). So 
\begin{eqnarray*}
H_{A_3\times \mathbf{\mu }_3}^8 &=&\mathbf{Z}/3\left\langle \alpha ^4,\alpha
^3\beta ,\alpha ^2\beta ^2,\alpha \beta ^3,\beta ^4\right\rangle \\
H_{A_3\times \mathbf{\mu }_3}^4 &=&\mathbf{Z}/3\left\langle \alpha ^2,\alpha
\beta ,\beta ^2\right\rangle
\end{eqnarray*}
and $\alpha \beta ^3\notin im\left( \cdot (-\alpha ^2)\right) $ i.e. $%
p\left( \alpha \beta ^3\right) \neq 0$.

(II) Consider the Leray spectral sequence 
\[
E_2^{pq}=H^p\left( \mathrm{B}PGL_3,H^q\left( \mathrm{B}\mathbf{\mu }_3,%
\mathbf{Z}\right) \right) \Rightarrow H^{p+q}\left( \mathrm{B}SL_3,\mathbf{Z}%
\right) . 
\]
By Lemma \ref{cohom}, its (first quadrant) $E_2$-term\footnote{%
We write only the parts we'll need.} is: 
\[
\begin{tabular}{|cccccccccccc|}
\hline
\multicolumn{1}{|c|}{$\cdots $} & \multicolumn{1}{c|}{} & 
\multicolumn{1}{c|}{} & \multicolumn{1}{c|}{} & \multicolumn{1}{c|}{} & 
\multicolumn{1}{c|}{} & \multicolumn{1}{c|}{} & \multicolumn{1}{c|}{} & 
\multicolumn{1}{c|}{} & \multicolumn{1}{c|}{} & \multicolumn{1}{c|}{} & $%
\vdots $ \\ \hline
\multicolumn{1}{|c|}{$\alpha ^4$} & \multicolumn{1}{c|}{$0$} & 
\multicolumn{1}{c|}{$\alpha ^3\xi _2$} & \multicolumn{1}{c|}{$\alpha ^4y_3$}
& \multicolumn{1}{c|}{$\alpha ^4y_4$} & \multicolumn{1}{c|}{$0$} & 
\multicolumn{1}{c|}{$\alpha ^4y_6$} & \multicolumn{1}{c|}{$\alpha ^3\xi _7$}
& \multicolumn{1}{c|}{$\alpha ^4x_8$, $\alpha ^4y_8$} & \multicolumn{1}{c|}{}
& \multicolumn{1}{c|}{} &  \\ \hline
\multicolumn{1}{|c|}{$0$} & \multicolumn{1}{c|}{$0$} & \multicolumn{1}{c|}{$%
0 $} & \multicolumn{1}{c|}{$0$} & \multicolumn{1}{c|}{$0$} & 
\multicolumn{1}{c|}{$0$} & \multicolumn{1}{c|}{$0$} & \multicolumn{1}{c|}{$0$%
} & \multicolumn{1}{c|}{$0$} & \multicolumn{1}{c|}{$0$} & 
\multicolumn{1}{c|}{$0$} &  \\ \hline
\multicolumn{1}{|c|}{$\alpha ^3$} & \multicolumn{1}{c|}{$0$} & 
\multicolumn{1}{c|}{$\alpha ^2\xi _2$} & \multicolumn{1}{c|}{$\alpha ^3y_3$}
& \multicolumn{1}{c|}{$\alpha ^3y_4$} & \multicolumn{1}{c|}{$0$} & 
\multicolumn{1}{c|}{$\alpha ^3y_6$} & \multicolumn{1}{c|}{$\alpha ^2\xi _7$}
& \multicolumn{1}{c|}{$\alpha ^3x_8$, $\alpha ^3y_8$} & \multicolumn{1}{c|}{}
& \multicolumn{1}{c|}{} &  \\ \hline
\multicolumn{1}{|c|}{$0$} & \multicolumn{1}{c|}{$0$} & \multicolumn{1}{c|}{$%
0 $} & \multicolumn{1}{c|}{$0$} & \multicolumn{1}{c|}{$0$} & 
\multicolumn{1}{c|}{$0$} & \multicolumn{1}{c|}{$0$} & \multicolumn{1}{c|}{$0$%
} & \multicolumn{1}{c|}{$0$} & \multicolumn{1}{c|}{$0$} & 
\multicolumn{1}{c|}{$0$} &  \\ \hline
\multicolumn{1}{|c|}{$\alpha ^2$} & \multicolumn{1}{c|}{$0$} & 
\multicolumn{1}{c|}{$\alpha \xi _2$} & \multicolumn{1}{c|}{$\alpha ^2y_3$} & 
\multicolumn{1}{c|}{$\alpha ^2y_4$} & \multicolumn{1}{c|}{$0$} & 
\multicolumn{1}{c|}{$\alpha ^2y_6$} & \multicolumn{1}{c|}{$\alpha \xi _7$} & 
\multicolumn{1}{c|}{$\alpha ^2x_8$, $\alpha ^2y_8$} & \multicolumn{1}{c|}{}
& \multicolumn{1}{c|}{} &  \\ \hline
\multicolumn{1}{|c|}{$0$} & \multicolumn{1}{c|}{$0$} & \multicolumn{1}{c|}{$%
0 $} & \multicolumn{1}{c|}{$0$} & \multicolumn{1}{c|}{$0$} & 
\multicolumn{1}{c|}{$0$} & \multicolumn{1}{c|}{$0$} & \multicolumn{1}{c|}{$0$%
} & \multicolumn{1}{c|}{$0$} & \multicolumn{1}{c|}{$0$} & 
\multicolumn{1}{c|}{$0$} &  \\ \hline
\multicolumn{1}{|c|}{$\alpha $} & \multicolumn{1}{c|}{$0$} & 
\multicolumn{1}{c|}{$\xi _2$} & \multicolumn{1}{c|}{$\alpha y_3$} & 
\multicolumn{1}{c|}{$\alpha y_4$} & \multicolumn{1}{c|}{$0$} & 
\multicolumn{1}{c|}{$\alpha y_6$} & \multicolumn{1}{c|}{$\xi _7$} & 
\multicolumn{1}{c|}{$\alpha x_8$, $\alpha y_8$} & \multicolumn{1}{c|}{} & 
\multicolumn{1}{c|}{} &  \\ \hline
\multicolumn{1}{|c|}{$0$} & \multicolumn{1}{c|}{$0$} & \multicolumn{1}{c|}{$%
0 $} & \multicolumn{1}{c|}{$0$} & \multicolumn{1}{c|}{$0$} & 
\multicolumn{1}{c|}{$0$} & \multicolumn{1}{c|}{$0$} & \multicolumn{1}{c|}{$0$%
} & \multicolumn{1}{c|}{$0$} & \multicolumn{1}{c|}{$0$} & 
\multicolumn{1}{c|}{$0$} &  \\ \hline
\multicolumn{1}{|c|}{$\mathbf{Z}$} & \multicolumn{1}{c|}{$0$} & 
\multicolumn{1}{c|}{$0$} & \multicolumn{1}{c|}{$y_3\mathbf{Z}/3$} & 
\multicolumn{1}{c|}{$y_4\mathbf{Z}/3$} & \multicolumn{1}{c|}{$0$} & 
\multicolumn{1}{c|}{$y_6\mathbf{Z}$} & \multicolumn{1}{c|}{$0$} & 
\multicolumn{1}{c|}{$x_8\mathbf{Z}\oplus y_8\mathbf{Z}/3$} & 
\multicolumn{1}{c|}{$0$} & \multicolumn{1}{c|}{$\mathbf{Z}$} & $\cdots $ \\ 
\hline
\end{tabular}
\]
where from the second row up, the coefficients are in $\mathbf{Z}/3$.

One of the edge maps is 
\[
H^8\left( \mathrm{B}PGL_3,\mathbf{Z}\right) =E_2^{8,0}\twoheadrightarrow
E_\infty ^{8,0}=F^8H^8\left( \mathrm{B}SL_3,\mathbf{Z}\right)
\hookrightarrow H^8\left( \mathrm{B}SL_3,\mathbf{Z}\right) 
\]
so we have to show that $F^8H^8\left( \mathrm{B}SL_3,\mathbf{Z}\right) \ $%
has index at least $9$ in 
\[
H^8(\mathrm{B}SL_3,\mathbf{Z})\simeq \mathbf{Z}\cdot \alpha _2^2. 
\]

First of all, note that $d_{(3)}(\alpha )=\pm y_3$ since 
\[
E_\infty ^{3,0}=F^3H^3\left( \mathrm{B}SL_3,\mathbf{Z}\right)
\hookrightarrow H^3\left( \mathrm{B}SL_3,\mathbf{Z}\right) =0 
\]
and both $\alpha $ and $y_3$ are $3$-torsion; we choose $y_3$ to have the
plus sign. Therefore 
\[
d_{(3)}\left( \alpha ^2y_3\right) =2\alpha y_3^2+\alpha ^2d_{(3)}(y_3)=0 
\]
since $y_3^2$ is $3$-torsion in $H^6\left( \mathrm{B}PGL_3,\mathbf{Z}\right)
\simeq \mathbf{Z}$, hence is zero.

Then 
\[
E_2^{62}=E_3^{62}=E_4^{62}=E_\infty ^{62}\simeq \mathbf{Z}/3 
\]
and we have the first $3$ factor of the desired index. Finally we have 
\[
d_{(3)}\left( \alpha ^2y_4\right) =2\alpha y_3y_4+\alpha ^2d_{(3)}(y_4)=0 
\]
since $y_3y_4\in H^7\left( \mathrm{B}PGL_3,\mathbf{Z}\right) =0$; then 
\[
E_2^{44}=E_3^{44}=E_4^{44}=E_\infty ^{44}\simeq \mathbf{Z}/3 
\]
yielding the other $3$ factor in the index of $F^8H^8\left( \mathrm{B}SL_3,%
\mathbf{Z}\right) \ $in $H^8(\mathrm{B}SL_3,\mathbf{Z})\simeq \mathbf{Z}%
\cdot \alpha _2^2$.

(III) As already observed, we have $c_2\left( \mathrm{sl}_3\right) ^2\mapsto
36\alpha _2^2$ via the pull back (use (I)) 
\[
\phi :\mathbf{Z}\oplus \mathrm{cl}\left( \rho \right) \cdot \left( \mathbf{Z/%
}3\right) \simeq H^8\left( \mathrm{B}PGL_3,\mathbf{Z}\right) \longrightarrow
H^8\left( \mathrm{B}SL_3,\mathbf{Z}\right) \simeq \mathbf{Z}\cdot \alpha
_2^2. 
\]
whose kernel is $3$-torsion; combining this with (II), we get that the image
of $\phi $ has exactly index $9$. Therefore 
\[
H^8\left( \mathrm{B}PGL_3,\mathbf{Z}_{(3)}\right) \simeq \mathbf{Z}%
_{(3)}\cdot j_{\bullet }\left( c_2\left( \mathrm{sl}_3\right) ^2\right)
\oplus \left( \mathbf{Z/}3\right) \cdot j_{\bullet }\left( \mathrm{cl}\left(
\rho \right) \right) , 
\]
where $j_{\bullet }:H^{*}(\mathrm{B}PGL_3,\mathbf{Z})\rightarrow H^{*}(%
\mathrm{B}PGL_3,\mathbf{Z}_{(3)})$ is the morphism induced by the
localization $j:\mathbf{Z}\rightarrow \mathbf{Z}_{(3)}$.

(IV) By \cite{KY} Cor. 4.7, we know that 
\[
H^8\left( \mathrm{B}PGL_3,\mathbf{Z}/3\right) \simeq \left( \mathbf{Z}%
/3\right) \cdot y_2^4\oplus \left( \mathbf{Z/}3\right) \cdot y_8, 
\]
and that the second generator $y_8$ is not in the Chern subring of $%
H^{*}\left( \mathrm{B}PGL_3,\mathbf{Z}/3\right) $. By the Bockstein exact
sequence, the natural map 
\[
\left( j_{(3)}\right) _{\bullet }:H^8\left( \mathrm{B}PGL_3,\mathbf{Z}%
_{(3)}\right) \rightarrow H^8\left( \mathrm{B}PGL_3,\mathbf{Z}/3\right) 
\]
is surjective since $H^9\left( \mathrm{B}PGL_3,\mathbf{Z}_{(3)}\right) =0$.
Therefore there exists an element $\xi =\alpha j_{\bullet }\left( c_2\left( 
\mathrm{sl}_3\right) ^2\right) +\beta j_{\bullet }\left( \mathrm{cl}\left(
\rho \right) \right) \in H^8\left( \mathrm{B}PGL_3,\mathbf{Z}_{(3)}\right) $
such that $\left( j_{(3)}\right) _{\bullet }\left( \xi \right) =y_8$. In
particular, $\mathrm{cl}\left( \rho \right) $ cannot be in the Chern subring
of $H^{*}\left( \mathrm{B}PGL_3,\mathbf{Z}\right) $. \TeXButton{End Proof}
{\endproof}

\begin{remark}
For a different proof of Theorem \ref{result1}, which does not depend on
Kono-Yagita's results on $H^{*}\left( \mathrm{B}PGL_3,\mathbf{Z}/3\right) $
(and in fact does not depend on cohomology at all), see the Appendix.
\end{remark}

\begin{lemma}
\label{teta3tors} $\Theta _{1,1}^{\left( 2\right) },\Theta _{0,1}^{\left(
2\right) },\Theta _{1,0}^{\left( 3\right) },\Theta _{1,0}^{\left( 4\right)
},\Theta _{0,0}^{(4)},\Theta _{1,0}^{\left( 5\right) },\Theta
_{0,0}^{(5)},\Theta _{1,0}^{\left( 6\right) },\Theta _{1,0}^{\left( 7\right)
}$ and $\Theta _{1,0}^{\left( 8\right) }$ are $3$-torsion.
\end{lemma}

\TeXButton{Proof}{\proof} All the $\Theta $'s are supported on the
complement of $U$ and so they all go to zero via $A_{PGL_3}^{*}\rightarrow
A_T^{*}$, since this map factors through $A_{PGL_3}^{*}\rightarrow
A_{PGL_3}^{*}\left( U\right) $. But, by \cite{EG1}, Prop. 6, the rational
pullback 
\[
A_{PGL_3}^{*}\otimes \mathbf{Q}\longrightarrow \left( A_T^{*}\right)
^{S_3}\otimes \mathbf{Q} 
\]
is an isomorphism, so the $\Theta $'s are torsion and hence $3$-torsion by
Cor. \ref{cortorsion}. \TeXButton{End Proof}{\endproof}

\begin{remark}
Note that $\mathrm{cl}\left( \chi \right) =0$ since $\chi $ is torsion while 
$H^{12}\left( \mathrm{B}PGL_3,\mathbf{Z}\right) $ is torsion free by lemma 
\ref{cohom}.
\end{remark}

\begin{lemma}
\label{grado4}$\Theta _{1,0}^{\left( 4\right) }$ and $\Theta _{0,0}^{\left(
4\right) }$ are in the kernel of the cycle map $\mathrm{cl}%
:A_{PGL_3}^{*}\rightarrow H^{*}(\mathrm{B}PGL_3,\mathbf{Z})$.
\end{lemma}

\TeXButton{Proof}{\proof} By part (I) of the proof of Th. \ref{result1}, $%
\mathrm{cl}\left( \rho \right) $ generates the $3$-torsion of $H^8(\mathrm{B}%
PGL_3,\mathbf{Z})$ and moreover $\mathrm{cl}\left( \rho \right) _{\mid
U}\neq 0$ in $H_{PGL_3}^8\left( U,\mathbf{Z}\right) $, where $U\subset 
\mathrm{sl}_3$ is the open subscheme of matrices with distinct eigenvalues.
Since $\Theta _{1,0}^{\left( 4\right) }$ and $\Theta _{0,0}^{\left( 4\right)
}$ are both $3$-torsion in $A_{PGL_3}^4$, we must have 
\begin{eqnarray*}
\mathrm{cl}\left( \Theta _{1,0}^{\left( 4\right) }\right) &=&A\cdot \mathrm{%
cl}\left( \rho \right) \\
\mathrm{cl}\left( \Theta _{0,0}^{\left( 4\right) }\right) &=&B\cdot \mathrm{%
cl}\left( \rho \right) \text{ .}
\end{eqnarray*}
But $\Theta _{1,0}^{\left( 4\right) }$ and $\Theta _{0,0}^{\left( 4\right) }$
have supports in the complement of $U$, so $A=B=0$. \TeXButton{End Proof}
{\endproof}

\begin{remark}
\label{teta8}Note that also the generator $\Theta _{1,0}^{\left( 8\right) }$
can be chosen in such a way that 
\[
\mathrm{cl}\left( \Theta _{1,0}^{\left( 8\right) }\right) =0.
\]
In fact $c_8\left( \mathrm{sl}_3\right) \neq 0$ in $H^{16}(\mathrm{B}PGL_3,%
\mathbf{Z})$ by \cite{KY}, Lemma 3.18 and 
\[
H^{16}(\mathrm{B}PGL_3,\mathbf{Z})\simeq \mathbf{Z}\oplus \mathbf{Z}\oplus 
\mathbf{Z}/3
\]
(Lemma. \ref{cohom}), therefore 
\[
\mathrm{cl}\left( \Theta _{1,0}^{\left( 8\right) }\right) =Ac_8\left( 
\mathrm{sl}_3\right) .
\]
Now observe that 
\[
c_8\left( \mathrm{sl}_3\right) _{\mid \mathrm{sl}_3\backslash Z_0\cup
\left\{ 0\right\} }=0
\]
while $\Theta _{1,0}^{\left( 8\right) }$ is a lift of $\left( j_{1,0}\right)
_{*}\left( \lambda _1\lambda _2^2\right) $ where 
\[
j_{1,0}:Z_{1,0}\hookrightarrow \mathrm{sl}_3\backslash Z_0\cup \left\{
0\right\} \text{ };
\]
thus we can choose a lift $\Theta _{1,0}^{\left( 8\right) }$ such that $A=0$.
\end{remark}

\begin{proposition}
\label{allzero}The elements 
\[
\left\{ \Theta _{1,1}^{\left( 2\right) },\Theta _{0,1}^{\left( 2\right)
},\Theta _{1,0}^{\left( 3\right) },\Theta _{1,0}^{\left( 4\right) },\Theta
_{0,0}^{(4)},\Theta _{1,0}^{\left( 5\right) },\Theta _{0,0}^{(5)},\Theta
_{1,0}^{\left( 6\right) },\Theta _{1,0}^{\left( 7\right) },\Theta
_{1,0}^{\left( 8\right) }\right\} 
\]
are all zero in $A_{PGL_3}^{*}$.
\end{proposition}

\TeXButton{Proof}{\proof}We first prove that 
\begin{equation}
\Theta _{1,1}^{\left( 2\right) }=\Theta _{0,1}^{\left( 2\right) }=\Theta
_{1,0}^{\left( 3\right) }=\Theta _{1,0}^{\left( 4\right) }=\Theta
_{0,0}^{(4)}=\Theta _{1,0}^{\left( 5\right) }=\Theta _{0,0}^{(5)}=\Theta
_{1,0}^{(7)}=0\text{.}  \label{firsthalf}
\end{equation}
Consider the commutative diagram 
\[
\begin{tabular}{lll}
$A_{PGL_3}^{*}$ & $\stackrel{cl}{\longrightarrow }$ & $H^{*}\left( \mathrm{B}%
PGL_3,\mathbf{Z}\right) $ \\ 
$\downarrow $ &  & $\downarrow $ \\ 
$A_{\Gamma _3}^{*}$ & $\stackunder{cl}{\longrightarrow }$ & $H^{*}\left( 
\mathrm{B}\Gamma _3,\mathbf{Z}\right) $%
\end{tabular}
\]
where the vertical arrows are injective by Theorem \ref{gottliebtotaro}. We
know that 
\[
\Theta _{1,1}^{\left( 2\right) },\Theta _{0,1}^{\left( 2\right) },\Theta
_{1,0}^{\left( 3\right) },\Theta _{1,0}^{\left( 4\right) },\Theta
_{0,0}^{(4)},\Theta _{1,0}^{\left( 5\right) },\Theta _{0,0}^{(5)},\Theta
_{1,0}^{(7)}
\]
are $3$-torsion and zero in cohomology (lemmas \ref{teta3tors}, \ref{cohom}
and \ref{grado4}), so (\ref{firsthalf}) will be proved if we show that 
\[
_3cl:\text{ }_3A_{\Gamma _3}^{*}\longrightarrow \text{ }_3H^{*}\left( 
\mathrm{B}\Gamma _3,\mathbf{Z}\right) 
\]
is injective up to degree $5$ and in degree $7$. But, by the usual
transfer-trick, the restriction induces isomorphisms 
\[
_3A_{\Gamma _3}^{*}\simeq \left( _3A_{A_3\ltimes T}^{*}\right) ^{C_2},\text{ 
}_3H^{*}\left( \mathrm{B}\Gamma _3,\mathbf{Z}\right) \simeq \left(
_3H^{*}\left( \mathrm{B}\left( A_3\ltimes T\right) ,\mathbf{Z}\right)
\right) ^{C_2}
\]
and it will be (more than) enough to show that 
\[
cl:A_{A_3\ltimes T}^{*}\longrightarrow H^{*}\left( \mathrm{B}\left(
A_3\ltimes T\right) ,\mathbf{Z}\right) 
\]
is injective up to degree $5$ and in degree $7$.

Recall (Prop. \ref{basico}) that $A_{A_3\ltimes T}^{*}$ is generated by 
\begin{equation}
\left\{ \alpha ,c_2\left( W\right) ,c_3\left( W\right) ,\theta \doteq 
\mathrm{tsf}_T^{A_3\ltimes T}\left( u_2^2u_3\right) \right\}  \label{gen}
\end{equation}
where $W$ is the representation defined in (\ref{vudoppio}) and we identify $%
A_T^{*}$ with 
\[
A_{T_{SL_3}}^{*}\simeq \mathbf{Z}\left[ u_1,u_2,u_3\right] \diagup \left(
u_1+u_2+u_3\right) ; 
\]
moreover (see Lemma \ref{basico2}), we have 
\begin{equation}
3\alpha =0,\text{ }\alpha \theta =0,\text{ }\alpha ^3+\alpha c_2\left(
W\right) =0,  \label{rel}
\end{equation}
\[
3\left[ \left( 2\theta +3c_3\left( W\right) \right) ^2+4c_2\left( W\right)
^3+27c_3\left( W\right) ^2\right] =0. 
\]
For the duration of this proof, we will denote $c_2\left( W\right) $ and $%
c_3\left( W\right) $ simply by $c_2$ and $c_3$; moreover, if $\xi \in
A_{A_3\ltimes T}^{*}$, we will write $\overline{\xi }$ for $cl\left( \xi
\right) $.

As shown in the proof of Th. \ref{result1}, we have 
\[
H^{2n}\left( \mathrm{B}\left( A_3\times \mathbf{\mu }_3\right) ,\mathbf{Z}%
\right) \simeq \left( H^{*}\left( \mathrm{B}A_3,\mathbf{Z}\right) \otimes
\left( \mathrm{B}\mathbf{\mu }_3,\mathbf{Z}\right) \right) ^{2n} 
\]
and 
\begin{equation}
\overline{c_2\left( W\right) }_{\mid A_3\times \mathbf{\mu }_3}=-\overline{%
\alpha }^2,\text{ }\overline{c_3\left( W\right) }_{\mid A_3\times \mathbf{%
\mu }_3}=\overline{\beta }\left( \overline{\beta }^2-\overline{\alpha }%
^2\right)  \label{a3mi3}
\end{equation}
where 
\[
H^{*}\left( \mathrm{B}A_3,\mathbf{Z}\right) =\frac{\mathbf{Z}\left[ 
\overline{\alpha }\right] }{\left( \overline{\alpha }\right) },\text{ }%
H^{*}\left( \mathrm{B}\mathbf{\mu }_3,\mathbf{Z}\right) =\frac{\mathbf{Z}%
\left[ \overline{\beta }\right] }{\left( 3\overline{\beta }\right) }. 
\]
In the following computations we will freely use that the cycle class map
respects Chern classes, restrictions and transfers and that 
\[
cl:A_T^{*}\longrightarrow H^{*}\left( \mathrm{B}T,\mathbf{Z}\right) 
\]
is an isomorphism.

If $\xi \in \ker cl\cap A_{A_3\ltimes T}^1,$ we have 
\[
\xi =A\alpha \text{ \quad and \quad }A\overline{\alpha }=0 
\]
for some $A\in \mathbf{Z}$; restricting this to $A_3$ (in cohomology) we
then get $A\equiv 0$ $\limfunc{mod}3$, hence $\xi =0$.

If $\xi \in \ker cl\cap A_{A_3\ltimes T}^2,$ we have 
\[
\xi =A\alpha ^2+Bc_2,\text{ }A\overline{\alpha }^2+B\overline{c_2}=0 
\]
for some $A,B\in \mathbf{Z}$; restricting to $T$, we get $B=0$ then,
restricting to $A_3$, we get $A\equiv 0$ $\limfunc{mod}3$. Therefore, $\xi
=0 $.

If $\xi \in \ker cl\cap A_{A_3\ltimes T}^3,$ we have 
\[
\xi =A\alpha ^3+Bc_3+C\theta 
\]
\[
A\overline{\alpha }^3+B\overline{c_3}+C\overline{\theta }=0 
\]
for some $A,B,C\in \mathbf{Z}$; restricting to $T$, we get $B=C=0$ since $%
\overline{c_3}_{\mid T}$ and $\overline{\theta }_{\mid T}$ are linearly
independent $H^{*}\left( \mathrm{B}T,\mathbf{Z}\right) $. Restricting then
to $A_3$, we get $A\equiv 0$ $\limfunc{mod}3$, hence $\xi =0$.

If $\xi \in \ker cl\cap A_{A_3\ltimes T}^4,$ we have 
\[
\xi =A\alpha ^4+B\alpha c_3+Cc_2^2 
\]
\[
A\overline{\alpha }^4+B\overline{\alpha }\overline{c_3}+C\overline{c_2}^2=0 
\]
for some $A,B,C\in \mathbf{Z}$; restricting to $T$, we get $C=0$.
Restricting then to $A_3\times \mathbf{\mu }_3$, from (\ref{a3mi3}) we get $%
B\equiv A\equiv 0$ $\limfunc{mod}3$, hence $\xi =0$.

If $\xi \in \ker cl\cap A_{A_3\ltimes T}^5,$ we have 
\[
\xi =A\alpha ^5+B\alpha ^2c_3+Cc_2c_3 
\]
\[
A\overline{\alpha }^5+B\overline{\alpha }^2\overline{c_3}+C\overline{c_2}%
\overline{c_3}=0 
\]
for some $A,B,C\in \mathbf{Z}$; restricting to $T$, we get $C=0$.
Restricting then to $A_3\times \mathbf{\mu }_3$, from (\ref{a3mi3}) we get $%
B\equiv A\equiv 0$ $\limfunc{mod}3$, hence $\xi =0$.

Finally, if $\xi \in \ker cl\cap A_{A_3\ltimes T}^7,$ we have 
\[
\xi =A\alpha ^7+B\alpha ^4c_3+Cc_2^2c_3+Dc_2^2\theta +E\alpha c_3^2 
\]
\[
A\overline{\alpha }^7+B\overline{\alpha }^4\overline{c_3}+C\overline{c_2}^2%
\overline{c_3}+D\overline{c_2}^2\overline{\theta }+E\overline{\alpha }%
\overline{c_3}^2=0 
\]
for some $A,B,C,D,E\in \mathbf{Z}$; restricting to $T$, we get $C=D=0$ since 
$\overline{c_2}_{\mid T}\neq 0$ and $\left( \overline{c_3}_{\mid T},%
\overline{\theta }_{\mid T}\right) $ are linearly independent in the domain $%
H^{*}\left( \mathrm{B}T,\mathbf{Z}\right) $. Restricting then to $A_3\times 
\mathbf{\mu }_3$, from (\ref{a3mi3}) we get $A\equiv B\equiv E\equiv 0$ $%
\limfunc{mod}3$, hence $\xi =0$. This concludes the proof of (\ref{firsthalf}%
).

Now we prove the remaining relations 
\begin{equation}
\Theta _{1,0}^{\left( 6\right) }=\Theta _{1,0}^{\left( 8\right) }=0.
\label{secondhalf}
\end{equation}
First observe that $\Theta _{1,0}^{\left( 6\right) }$ and $\Theta
_{1,0}^{\left( 8\right) }$ are $3$-torsion and zero in cohomology (with $%
\Theta _{1,0}^{\left( 8\right) }$ chosen as in Remark \ref{teta8}). Since
they are lifts of elements having supports in the complement of $U\subset 
\mathrm{sl}_3$, their restrictions to $A_{\Gamma _3}^{*}$ are in the kernel
of 
\[
A_{\Gamma _3}^{*}\longrightarrow A_{\Gamma _3}^{*}\left( Diag_{\mathrm{sl}%
_3}^{*}\right) \simeq A_{PGL_3}^{*}\left( U\right) \text{.} 
\]
and in particular: 
\[
\left\{ \Theta _{1,0\mid A_3\ltimes T}^{\left( 6\right) },\Theta _{1,0\mid
A_3\ltimes T}^{\left( 8\right) }\right\} \subset \ker \left( g:A_{A_3\ltimes
T}^{*}\longrightarrow A_{A_3\ltimes T}^{*}\left( Diag_{\mathrm{sl}%
_3}^{*}\right) \right) . 
\]
By Lemma \ref{basico2} (ii) and (\ref{gen}), (\ref{rel}), we must have 
\begin{eqnarray*}
\Theta _{1,0\mid A_3\ltimes T}^{\left( 6\right) } &=&\alpha ^2\left( A\alpha
^4+B\alpha c_3+Cc_2^2\right) \\
\Theta _{1,0\mid A_3\ltimes T}^{\left( 8\right) } &=&\alpha ^2\left( D\alpha
^6+Ec_2^3+Fc_3^2+G\alpha ^3c_3\right)
\end{eqnarray*}
for some $A,...,E\in \mathbf{Z}$. Using again (\ref{rel}), we get 
\begin{eqnarray*}
\Theta _{1,0\mid A_3\ltimes T}^{\left( 6\right) } &=&A^{\prime }\alpha
^6+B\alpha ^3c_3 \\
\Theta _{1,0\mid A_3\ltimes T}^{\left( 8\right) } &=&D^{\prime }\alpha
^6+Fc_3^2+G\alpha ^3c_3
\end{eqnarray*}
for some $A^{\prime },B,D^{\prime },F,G\in \mathbf{Z}$. Again denoting $%
cl\left( \xi \right) $ by $\overline{\xi }$ for $\xi \in A_{A_3\ltimes
T}^{*} $, we have 
\begin{eqnarray*}
0 &=&\overline{\Theta _{1,0\mid A_3\ltimes T}^{\left( 6\right) }}=A^{\prime }%
\overline{\alpha }^6+B\overline{\alpha }^3\overline{c_3} \\
0 &=&\overline{\Theta _{1,0\mid A_3\ltimes T}^{\left( 8\right) }}=D^{\prime }%
\overline{\alpha }^6+F\overline{c_3}^2+G\overline{\alpha }^3\overline{c_3}
\end{eqnarray*}
in $H^{*}\left( \mathrm{B}\left( A_3\ltimes T\right) ,\mathbf{Z}\right) $.
Restricting these relations to $A_3\times \mathbf{\mu }_3$, by (\ref{a3mi3})
we obtain: 
\begin{eqnarray*}
A^{\prime } &\equiv &B\equiv 0\text{ }\limfunc{mod}3 \\
D^{\prime } &\equiv &F\equiv G\equiv 0\text{ }\limfunc{mod}3
\end{eqnarray*}
i.e. 
\begin{equation}
\Theta _{1,0\mid A_3\ltimes T}^{\left( 6\right) }=\text{ }\Theta _{1,0\mid
A_3\ltimes T}^{\left( 8\right) }=0  \label{qfinito}
\end{equation}
in $A_{A_3\ltimes T}^{*}$. But the restriction map induces an isomorphism 
\[
_3A_{\Gamma _3}^{*}\simeq \left( _3A_{A_3\ltimes T}^{*}\right) ^{C_2} 
\]
and then, we also get 
\[
\Theta _{1,0\mid \Gamma _3}^{\left( 6\right) }=\text{ }\Theta _{1,0\mid
\Gamma _3}^{\left( 8\right) }=0 
\]
in $A_{\Gamma _3}^{*}$. By Theorem \ref{gottliebtotaro} we finally get (\ref
{secondhalf}). \TeXButton{End Proof}{\endproof}

Thus we can summarize the main result obtained so far in the following:

\begin{theorem}
\label{generators}With the notation of (\ref{laststep}), $A_{PGL_3}^{*}$ is
generated by 
\begin{equation}
\left\{ 2c_2\left( \mathrm{sl}_3\right) -c_2\left( Sym^3E\right) ,c_3\left(
Sym^3E\right) ,\rho ,\chi ,c_6\left( \mathrm{sl}_3\right) ,c_8\left( \mathrm{%
sl}_3\right) \right\}   \label{generatori}
\end{equation}
where\emph{\ }$\deg \rho =4$\emph{, }$\deg \chi =6$.
\end{theorem}

\begin{remark}
We point out that 
\[
2c_2\left( \mathrm{sl}_3\right) -c_2\left( Sym^3E\right) ,c_3\left(
Sym^3E\right) ,c_6\left( \mathrm{sl}_3\right) ,c_8\left( \mathrm{sl}%
_3\right) 
\]
are nonzero (by checking their images in $A_{SL_3}^{*}$ or in cohomology)
and we will show in the next section that $\rho \neq 0$. Unfortunately, we
do not know whether $\chi $ is zero or not.

Note also that the generators $\rho $ and $\chi $, defined originally as
lifts from the open subset $U$ (therefore not unique \emph{a priori}) are
indeed uniquely defined since they have degrees $<8$ and $c_8\left( \mathrm{%
sl}_3\right) $ is the only generator coming from the complement of $U$.
\end{remark}

%%%%%%%%%%%%%%%%%%%%%%%%%%% End PGL_34.tex %%%%%%%%%%%%%%%%%%%%%%%%%%%%%%%
}

\QSubDoc{Include PGL_35}{%%%%%%%%%%%%%%%%%%%%%%%%%% Start PGL_35.tex %%%%%%%%%%%%%%%%%%%%%%%%%%%%%%

\LaTeXparent{PGL_3}
\ChildStyles{amssymb} 
\ChildDefaults{chapter:5,page:1}

\section{Other relations and results on the cycle maps\ }

With the notations established in the preceding sections we have:

\begin{proposition}
\label{relations} The following relations hold among the generators of $%
A_{PGL_3}^{*}$: 
\[
3\rho =3\chi =3c_8\left( \mathrm{sl}_3\right) =0
\]
\[
3\left( 27c_6\left( \mathrm{sl}_3\right) -c_3\left( Sym^3E\right)
^2-4\lambda ^3\right) =0
\]
\[
\rho ^2=c_8\left( \mathrm{sl}_3\right) .
\]
\end{proposition}

\TeXButton{Proof}{\proof}The pullback $\varphi :A_{PGL_3}^{*}\rightarrow
A_T^{*}$ factors through the composition 
\[
\pi :A_{PGL_3}^{*}\twoheadrightarrow A_{PGL_3}^{*}\left( U\right) \simeq
A_{\Gamma _3}^{*}\left( Diag_{\mathrm{sl}_3}^{*}\right) \twoheadrightarrow
A_T^{*}\left( Diag_{\mathrm{sl}_3}^{*}\right) ^{S_3}=\left( A_T^{*}\right)
^{S_3}\text{, } 
\]
and, by definition of $\chi $ and $\rho $, $\pi \left( \chi \right) =\pi
\left( \rho \right) =0$. Since (\cite{EG1}, Prop. 6) the rational pullback $%
\varphi _{\mathbf{Q}}$ is an isomorphism, $\chi $ and $\rho $ are torsion
and then $3$-torsion by Cor. \ref{cortorsion}.

Since $\mathrm{sl}_3=E\otimes E^{\vee }-\mathbf{1}$, as $SL_3$%
-representations ($E$ being the standard representation), $c_8\left( \mathrm{%
sl}_3\right) $ is in the kernel of $A_{PGL_3}^{*}\rightarrow A_{SL_3}^{*}$,
so it is $3$-torsion (Prop. \ref{torsion}).

A long but straightforward computation\footnote{%
The basic fact here is that $c_6\left( \mathrm{sl}_3\right) $ restricts to
minus the discriminant, $4\alpha _2^3+27\alpha _3^2$, in $A_{SL_3}^{*}=%
\mathbf{Z}\left[ \alpha _2,\alpha _3\right] $, where $\alpha _i=c_i\left(
E\right) $.} shows that 
\[
27c_6\left( \mathrm{sl}_3\right) -c_3\left( Sym^3E\right) ^2-4\lambda ^3\in
\ker \left( A_{PGL_3}^{*}\rightarrow A_{SL_3}^{*}\right) 
\]
so that this element is $3$-torsion (again by Prop. \ref{torsion}).

By definition of $\rho $ and Lemma \ref{basico2}, we have 
\begin{equation}
\text{ }\rho _{\mid A_3\ltimes T}=\alpha c_3\left( W\right) +A\alpha ^4
\label{ro}
\end{equation}
for some $A\in \mathbf{Z}/3$. Since 
\[
_3A_{\Gamma _3}^{*}\simeq \left( _{\text{ }3}A_{A_3\ltimes T}^{*}\right)
^{C_2}, 
\]
by \label{lastrelation}Lemma \ref{basico2} (ii), $\rho ^2$ belongs to the
kernel of 
\[
A_{PGL_3}^{*}\longrightarrow A_{PGL_3}^{*}\left( U\right) \simeq A_{\Gamma
_3}^{*}\left( Diag_{\mathrm{sl}_3}^{*}\right) . 
\]
Therefore, since by Proposition \ref{allzero} all the generators of $%
A_{PGL_3}^{*}$ coming from the complement of $U$ are zero except for $%
c_8\left( \mathrm{sl}_3\right) $, we have 
\begin{equation}
\rho ^2=Bc_8\left( \mathrm{sl}_3\right)  \label{roquadro}
\end{equation}
for some $B\in \mathbf{Z}/3$.

Let us determine $A$ and $B$. Since $c_8\left( \mathrm{sl}_3\right) _{\mid
A_3}=0$, from (\ref{ro}) and (\ref{roquadro}), we get $A=0$ i.e. 
\begin{equation}
\rho _{\mid A_3\ltimes T}=\alpha c_3\left( W\right) .  \label{ro'}
\end{equation}
Straightforward computations show that 
\[
c_8\left( \mathrm{sl}_3\right) _{\mid A_3\times \mathbf{\mu }_3}=\alpha
^2\beta ^2\left( \beta ^2-\alpha ^2\right) ^2 
\]
\[
c_3\left( W\right) _{\mid A_3\times \mathbf{\mu }_3}=\beta \left( \beta
^2-\alpha ^2\right) 
\]
in $A_{_{A_3\times \mathbf{\mu }_3}}^{*}=A_{A_3}^{*}\otimes A_{\mathbf{\mu }%
_3}^{*}=\mathbf{Z}\left[ \alpha \right] /\left( 3\alpha \right) \otimes 
\mathbf{Z}\left[ \beta \right] /\left( 3\beta \right) $ and then (\ref
{roquadro}) and (\ref{ro'}) prove that $B=1$. \TeXButton{End Proof}
{\endproof}

We define the graded ring 
\[
R^{*}\doteq \frac{\mathbf{Z}\left[ \lambda ,c_3\left( Sym^3E\right) ,\rho
,\chi ,c_6\left( \mathrm{sl}_3\right) ,c_8\left( \mathrm{sl}_3\right)
\right] }{\frak{R}} 
\]
where 
\[
\frak{R}\doteq \left( 3\rho ,3\chi ,3c_8\left( \mathrm{sl}_3\right) ,3\left(
27c_6\left( \mathrm{sl}_3\right) -c_3\left( Sym^3E\right) ^2-4\lambda
^3\right) ,\rho ^2-c_8\left( \mathrm{sl}_3\right) \right) 
\]
and\emph{\ }$\deg \rho =4$,\emph{\ }$\deg \chi =6$.

This is our candidate for $A_{PGL_3}^{*}$. What we do know is that the
canonical morphism 
\[
\pi :R^{*}\longrightarrow A_{PGL_3}^{*} 
\]
is surjective (Th. \ref{generators}).

\begin{remark}
Note that it is immediately clear that $\pi _{\mathbf{Q}}:R^{*}\otimes 
\mathbf{Q}\rightarrow A_{PGL_3}^{*}\otimes \mathbf{Q}$ is an isomorphism. In
fact 
\[
R^{*}\otimes \mathbf{Q}=\frac{\mathbf{Q}\left[ \lambda ,c_3\left(
Sym^3E\right) ,c_6\left( \mathrm{sl}_3\right) \right] }{\left( 27c_6\left( 
\mathrm{sl}_3\right) -c_3\left( Sym^3E\right) ^2-4\lambda ^3\right) }= 
\]
\[
=\mathbf{Q}\left[ \lambda ,c_3\left( Sym^3E\right) \right] . 
\]
Moreover, $\lambda \mapsto 3\alpha _2$ and $c_3\left( Sym^3E\right) \mapsto
27\alpha _3$ via 
\[
A_{PGL_3}^{*}\longrightarrow A_{SL_3}^{*}=\mathbf{Z}\left[ \alpha _2,\alpha
_3\right] 
\]
which is rationally an isomorphism (Prop. \ref{rational}). We will prove in
Proposition \ref{result0} (ii) that more is true: $R^{*}$ and $A_{PGL_3}^{*}$
are isomorphic after inverting $3$.
\end{remark}

We will now establish some properties of the cycle map 
\[
\mathrm{cl}:A_{PGL_3}^{*}\longrightarrow H^{*}\left( \mathrm{B}PGL_3,\mathbf{%
Z}\right) 
\]
and of Totaro's refined cycle map 
\[
\widetilde{\mathrm{cl}}:A_{PGL_3}^{*}\longrightarrow MU^{*}\left( \mathrm{B}%
PGL_3\right) \otimes _{MU^{*}}\mathbf{Z.} 
\]

\begin{remark}
\label{totaro}In \cite{KY} Kono and Yagita proved that in the
Atiyah-Hirzebruch spectral sequence for Brown-Peterson cohomology at the
prime $3$ (\cite{W}) 
\[
E_2^{pq}=H^p\left( \mathrm{B}PGL_3,BP^q\right) \Longrightarrow
BP^{p+q}\left( \mathrm{B}PGL_3\right) 
\]
the $E_\infty $-term is generated as a $BP^{*}$-module by the top row i.e.
by 
\[
im\left( BP^{*}\left( \mathrm{B}PGL_3\right) \longrightarrow H^{*}\left( 
\mathrm{B}PGL_3,\mathbf{Z}_{\left( 3\right) }\right) \right) . 
\]
As a consequence, the natural map 
\[
\underline{\mathrm{cl}}:MU^{*}\left( \mathrm{B}PGL_3\right) \otimes _{MU^{*}}%
\mathbf{Z}\longrightarrow H^{*}\left( \mathrm{B}PGL_3,\mathbf{Z}\right) 
\]
is injective.
\end{remark}

We have the following result\footnote{%
A stronger version of (i) will be proved in Theorem \ref{riassunto}.}:

\begin{proposition}
(i) \label{result0}$\mathrm{cl}$ and $\widetilde{\mathrm{cl}}$ are injective
after inverting $3$;

(ii) $\pi $ is an isomorphism after inverting $3$.
\end{proposition}

\TeXButton{Proof}{\proof} (i) $A_{PGL_3}^{*}$ has only $3$-torsion and $\ker 
\mathrm{cl}$ is torsion (Section 2). Therefore $\mathrm{cl}=\underline{%
\mathrm{cl}}\circ \widetilde{\mathrm{cl}}$ is injective after inverting $3$
and the same is true for $\widetilde{\mathrm{cl}}$.

(ii) It is enough to prove that for any prime $p\neq 3$, the composition%
\footnote{$\left( \cdot \right) _{\left( p\right) }$ denotes localization at
the prime $p$.} 
\[
\left( R^{*}\right) _{\left( p\right) }\stackrel{\pi _{(p)}}{\longrightarrow 
}\left( A_{PGL_3}^{*}\right) _{\left( p\right) }\stackrel{\mathrm{cl}_{(p)}}{%
\longrightarrow }H^{*}(\mathrm{B}PGL_3,\mathbf{Z}_{\left( p\right) }) 
\]
is injective. Leray spectral sequence with $\mathbf{Z}_{\left( p\right) }$%
-coefficients: 
\[
E_2^{pq}=H^p\left( \mathrm{B}PGL_3,H^q\left( \mathrm{B}\mathbf{\mu }_3,%
\mathbf{Z}_{\left( p\right) }\right) \right) \Longrightarrow H^{p+q}\left( 
\mathrm{B}SL_3,\mathbf{Z}_{\left( p\right) }\right) 
\]
collapses at the $E_2$-term since $H^{*}\left( \mathrm{B}\mathbf{\mu }_3,%
\mathbf{Z}_{\left( p\right) }\right) =\mathbf{Z}_{\left( p\right) }$,
concentrated in degree zero, thus yielding an ''edge'' isomorphism
(coinciding with the pullback): 
\[
\varphi _{\left( p\right) }:H^{*}\left( \mathrm{B}PGL_3,\mathbf{Z}_{\left(
p\right) }\right) \simeq H^{*}\left( \mathrm{B}SL_3,\mathbf{Z}_{\left(
p\right) }\right) =\mathbf{Z}_{\left( p\right) }\left[ \alpha _2,\alpha
_3\right] \text{.} 
\]
Now, consider the commutative diagram 
\begin{equation}
\begin{tabular}{ccccc}
$\left( R^{*}\right) _{\left( p\right) }$ & $\stackrel{\pi _{\left( p\right)
}}{\longrightarrow }$ & $\left( A_{PGL_3}^{*}\right) _{\left( p\right) }$ & $%
\stackrel{\mathrm{cl}_{\left( p\right) }}{\longrightarrow }$ & $H^{*}\left( 
\mathrm{B}PGL_3,\mathbf{Z}_{\left( p\right) }\right) $ \\ 
&  & $^{\phi _{\left( p\right) }}\downarrow $ &  & $\downarrow ^{\varphi
_{\left( p\right) }}$ \\ 
&  & $\left( A_{SL_3}^{*}\right) _{\left( p\right) }$ & $\stackrel{}{%
\stackunder{\mathrm{cl}_{SL_3,\left( p\right) }}{\widetilde{\longrightarrow }%
}}$ & $H^{*}\left( \mathrm{B}SL_3,\mathbf{Z}_{\left( p\right) }\right) $%
\end{tabular}
\label{comm}
\end{equation}
and observe that for $p\neq 3$, 
\[
\left( R^{*}\right) _{\left( p\right) }=\frac{\mathbf{Z}_{\left( p\right)
}\left[ \lambda ,c_3\left( Sym^3E\right) ,c_6\left( \mathrm{sl}_3\right)
\right] }{\left( 27c_6\left( \mathrm{sl}_3\right) -c_3\left( Sym^3E\right)
^2-4\lambda ^3\right) }=\mathbf{Z}_{\left( p\right) }\left[ \lambda
,c_3\left( Sym^3E\right) \right] . 
\]
Since, as we already computed, $\phi \circ \pi \left( \lambda \right)
=3\alpha _2$, $\phi \circ \pi \left( c_3\left( Sym^3E\right) \right)
=27\alpha _3$, commutativity of (\ref{comm}) concludes the proof. 
\TeXButton{End Proof}{\endproof}

The stronger result we can prove about $\widetilde{\mathrm{cl}}$ is the
following

\begin{theorem}
\label{riassunto}Totaro's refined cycle class map 
\[
\widetilde{\mathrm{cl}}:A_{PGL_3}^{*}\longrightarrow MU^{*}(\mathrm{B}%
PGL_3)\otimes _{MU^{*}}\mathbf{Z} 
\]
is surjective (and has $3$-torsion kernel).
\end{theorem}

\TeXButton{Proof}{\proof} $\ker \widetilde{\mathrm{cl}}$ is $3$-torsion
since it is torsion and $A_{PGL_3}^{*}$ has only $3$-torsion. So we are left
to prove surjectivity of $\widetilde{\mathrm{cl}}$. To do this, we first
prove that $\widetilde{\mathrm{cl}}$ is surjective (thus an isomorphism by
Prop. \ref{result0} (i)) after inverting $3$ and then that $\widetilde{%
\mathrm{cl}}$ is surjective when localized at the prime $3$.

$\underline{\mathrm{cl}}_{PGL_3}$ is an isomorphism after inverting $3$
since $H^{*}\left( \mathrm{B}PGL_3,\mathbf{Z}\left[ \frac 13\right] \right) $
is torsion free

(\cite{To1})\footnote{%
We briefly sketch the argument. Since the differentials in the
Atiyah-Hirzebruch spectral sequence 
\[
F_2^{pq}=H^p\left( BPGL_3,MU^q\right) \Rightarrow MU^{p+q}\left(
BPGL_3\right) 
\]
are always torsion, they must be $0$ if $3$ is inverted since there is only $%
3$-torsion (recall that $MU^{*}$ is torsion-free). Therefore $F_2^{pq}$
collapses when $3$ is inverted.}. So it is enough to prove that \textrm{cl}$%
_{PGL_3}$ is surjective when $3$ is inverted. Now, in the commutative
diagram 
\[
\begin{tabular}{cc}
$A_{PGL_3}^{*}\left[ \frac 13\right] \stackrel{\mathrm{cl}_{PGL_3}\left[
\frac 13\right] }{\longrightarrow }$ & $H_{PGL_3}^{*}\left[ \frac 13\right] $
\\ 
$^\varphi \downarrow \quad $ & $^{\varphi ^{\prime }}\downarrow \quad $ \\ 
$A_{SL_3}^{*}\left[ \frac 13\right] \stackunder{\mathrm{cl}_{SL_3}\left[
\frac 13\right] }{\longrightarrow }$ & $H_{SL_3}^{*}\left[ \frac 13\right] $%
\end{tabular}
\]
$\varphi ^{\prime }$ is an isomorphism since the corresponding Leray
spectral sequence 
\[
E_2^{pq}=H^p\left( \mathrm{B}PGL_3,H^q\left( \mathrm{B}\mathbf{\mu }_3,%
\mathbf{Z}\right) \right) \Rightarrow H^{p+q}\left( \mathrm{B}SL_3,\mathbf{Z}%
\right) 
\]
collapses after inverting $3$, and \textrm{cl}$_{SL_3}$ is an isomorphism
even without inverting $3$. On the other hand, $\varphi $ is injective
because $\Phi :A_{PGL_3}^{*}\rightarrow A_{SL_3}^{*}$ has $3$-torsion kernel
and is surjective since 
\begin{eqnarray*}
&&\ c_2\left( \mathrm{sl}_3\right) \stackrel{\Phi }{\longmapsto }6\alpha _2
\\
&&\ c_2\left( Sym^3E\right) \stackrel{\Phi }{\longmapsto }15\alpha _2 \\
&&\ c_3\left( Sym^3E\right) \stackrel{\Phi }{\longmapsto }27\alpha _3.
\end{eqnarray*}
Therefore $\mathrm{cl}_{PGL_3}\left[ \frac 13\right] $ is an isomorphism too.

So it remains to prove that the localization at the prime $3$%
\[
\left( \widetilde{\mathrm{cl}}_{PGL_3}\right) _{\left( 3\right) }:\left(
A_{PGL_3}^{*}\right) _{\left( 3\right) }\longrightarrow MU^{*}(\mathrm{B}%
PGL_3)\otimes _{MU^{*}}\mathbf{Z}_{\left( 3\right) } 
\]
is surjective. By \cite{Q}, 
\[
MU^{*}(\mathrm{B}PGL_3)\otimes _{MU^{*}}\mathbf{Z}_{\left( 3\right) }\simeq
BP^{*}(\mathrm{B}PGL_3)\otimes _{BP^{*}}\mathbf{Z}_{\left( 3\right) } 
\]
where $BP^{*}\left( X\right) $ denotes the Brown-Peterson cohomology of $X$
localized at the prime $3$ and 
\[
BP^{*}=BP^{*}\left( pt\right) =\mathbf{Z}_{\left( 3\right) }\left[
v_1,\ldots ,v_n,\ldots \right] \longrightarrow \mathbf{Z}_{\left( 3\right) } 
\]
($\deg v_i=-2\left( 3^i-1\right) $) sends each $v_i$ to zero (see also \cite
{W}). Kono and Yagita computed $BP^{*}(\mathrm{B}PGL_3)$ in \cite{KY}, Th.
4.9, as a $BP^{*}$-module; it is a quotient of the following $BP^{*}$-module 
\[
\left( BP^{*}\mathbf{Z}_{\left( 3\right) }\left[ \left[ \widetilde{y_2}%
\right] \right] \widetilde{y_2}^2\oplus BP^{*}\oplus BP^{*}\mathbf{Z}%
_{\left( 3\right) }\left[ \left[ \widetilde{y_8}\right] \right] \widetilde{%
y_8}\right) \otimes \mathbf{Z}_{\left( 3\right) }\left[ \left[ \widetilde{%
y_{12}}\right] \right] 
\]
and, if 
\[
r:BP^{*}(\mathrm{B}PGL_3)\stackrel{s}{\longrightarrow }H^{*}\left( \mathrm{B}%
PGL_3,\mathbf{Z}_{\left( 3\right) }\right) \stackrel{j_{\bullet }}{%
\longrightarrow }H^{*}\left( \mathrm{B}PGL_3,\mathbf{Z}/3\right) 
\]
(where $s$ is the natural map of generalized cohomology theories and $%
j_{\bullet }$ is induced by $j:\mathbf{Z}_{\left( 3\right) }\rightarrow 
\mathbf{Z}/3$), $r$ has kernel $BP^{<0}\cdot BP^{*}(\mathrm{B}PGL_3)$ and 
\begin{eqnarray*}
r\left( \widetilde{y_2}^2\right) &=&y_2^2\equiv c_2\left( \mathrm{sl}%
_3\right) \\
r\left( \widetilde{y_8}\right) &=&y_8^{} \\
r\left( \widetilde{y_{12}}\right) &=&y_{12}\equiv c_6\left( \mathrm{sl}%
_3\right) ,
\end{eqnarray*}
$y_8\in H^8\left( \mathrm{B}PGL_3,\mathbf{Z}/3\right) $ being the same as in
part (IV) of the proof of Th. \ref{result1}. So we only need to show that $%
\widetilde{y_8}$ is in the image of $\left( \widetilde{\mathrm{cl}}%
_{PGL_3}\right) _{\left( 3\right) }$. By part (IV) of the proof of Th. \ref
{result1}, $y_8$ is in the image of 
\[
j_{\bullet }\circ \left( \mathrm{cl}_{PGL_3}\right) _{\left( 3\right)
}:\left( A_{PGL_3}^4\right) _{\left( 3\right) }\longrightarrow H^8\left( 
\mathrm{B}PGL_3,\mathbf{Z}/3\right) , 
\]
and this concludes the proof since $r$ has kernel $BP^{<0}\cdot BP^{*}(%
\mathrm{B}PGL_3)$. \TeXButton{End Proof}{\endproof}

\begin{remark}
We wish to point out that we do not know whether $\ker \widetilde{\mathrm{cl}%
}$ is zero or not. Moreover, since $\mathrm{cl}\left( \chi \right) =0$ and $%
\underline{\mathrm{cl}}$ is injective (Remark \ref{totaro}), we also have $%
\widetilde{\mathrm{cl}}\left( \chi \right) =0$. Therefore, if Totaro's
conjecture was true (i.e $\widetilde{\mathrm{cl}}$ was an isomorphism) we
should have $\chi =0$; but, again, we are not able to prove whether $\chi =0$
or not.
\end{remark}

%%%%%%%%%%%%%%%%%%%%%%%%%%% End PGL_35.tex %%%%%%%%%%%%%%%%%%%%%%%%%%%%%%%
}

\QSubDoc{Include PGL_36}{%%%%%%%%%%%%%%%%%%%%%%%% Start PGL_36.tex %%%%%%%%%%%%%%%%%%%%%%%%%%%%%

\LaTeXparent{PGL_3}
\ChildStyles{amssymb} 
\ChildDefaults{chapter:5,page:1}

\section{Appendix. A cohomology-independent proof that $A_{PGL_3}^{*}$ is
not generated by Chern classes}

Here we give an alternative proof of Theorem \ref{result1} which is
independent of Kono-Mimura-Shimada's results on the $\mathbf{Z}/3$%
-cohomology of $\mathrm{B}PGL_3$ and deals only with Chow rings with no
reference to cohomology. However, for the same reason, the following proof
does not yield any direct information on the cycle or refined cycle map.

The notations are those of the previous sections.

\begin{proposition}
\label{RRing}The representation ring of $PGL_3$ is generated by $\left\{ 
\mathrm{sl}_3,Sym^3E,Sym^3E^{\vee }\right\} $.
\end{proposition}

\TeXButton{Proof}{\proof} The exact sequence 
\[
1\rightarrow \mathbf{\mu }_3\longrightarrow SL_3\longrightarrow
PGL_3\rightarrow 1
\]
induces an exact sequence of character groups 
\[
0\rightarrow \widehat{T_{PGL_3}}\equiv \widehat{T}\longrightarrow \widehat{%
T_{SL_3}}\stackrel{\pi }{\longrightarrow }\mathbf{Z}/3\rightarrow 0
\]
where $\widehat{T_{SL_3}}=\mathbf{Z}^3/\mathbf{Z}$ ($\mathbf{Z}%
\hookrightarrow \mathbf{Z}^3$ diagonally) and $\pi :\left[
n_1,n_2,n_3\right] \mapsto \left[ n_1+n_2+n_3\right] $. Then 
\[
\mathbf{Z}\left[ \widehat{T}\right] \hookrightarrow \mathbf{Z}\left[ 
\widehat{T_{SL_3}}\right] =\mathbf{Z}\left[ x_1,x_2,x_3\right] \diagup
\left( x_1x_2x_3-1\right) 
\]
is the subring generated by monomials $x_1^{n_1}x_2^{n_2}x_3^{n_3}$ with $%
n_1+n_2+n_3\equiv 0$ $\limfunc{mod}3$. Therefore 
\[
R\left( PGL_3\right) =\left( R\left( T\right) \right) ^{S_3}=\left( \mathbf{Z%
}\left[ \widehat{T}\right] \right) ^{S_3}\hookrightarrow R\left( SL_3\right)
=\left( \mathbf{Z}\left[ \widehat{T_{SL_3}}\right] \right) ^{S_3}=\mathbf{Z}%
\left[ s_1,s_2\right] 
\]
(where $s_i$ is the $i$-th elementary symmetric function on the $x_i$'s) is
the subring generated by $\left\{ s_1^3,s_1s_2,s_2^3\right\} $. Then to
prove the Proposition it is enough to compute $\mathrm{sl}_3,Sym^3E$ and $%
Sym^3E^{\vee }$ in terms of $s_1$ and $s_2$ \emph{in }$R\left( SL_3\right) $.

If $E$ is the standard representation and $\mathbf{1}$ the trivial one
dimensional representation of $SL_3$, we have 
\[
E=x_1+x_2+x_3=s_1 
\]
\[
E^{\vee }=x_1^{-1}+x_2^{-1}+x_3^{-1}=x_1x_2+x_1x_3+x_2x_3=s_2 
\]
\[
\mathrm{sl}_3=E\otimes E^{\vee }-\mathbf{1}=s_1s_2-1; 
\]
so 
\[
Sym^3E=s_1^3-2s_1s_2+1,\text{ }Sym^3E^{\vee }=s_2^3-2s_1s_2+1 
\]
and we conclude. \TeXButton{End Proof}{\endproof}

\begin{corollary}
\label{corRRing}The Chern subring $A_{Ch,PGL_3}^{*}$ of $A_{PGL_3}^{*}$,
generated by Chern classes of representations, is generated by $\left\{
c_i\left( \mathrm{sl}_3\right) ,c_j\left( Sym^3E\right) \right\} _{i,j\geq 0}
$.
\end{corollary}

\begin{theorem}
$\rho $ is not in the Chern subring of $A_{PGL_3}^{*}$.
\end{theorem}

\TeXButton{Proof}{\proof} By Prop. \ref{RRing}, 
\[
R\left( PGL_3\right) =\mathbf{Z}\left[ \mathrm{sl}_3,Sym^3E,Sym^3E^{\vee
}\right] 
\]
and since $\mathrm{sl}_3$ (respectively, $Sym^3E$) is isomorphic to the
regular $A_3\times \mathbf{\mu }_3$-representation minus the trivial one
(respectively, plus the trivial one), we have 
\begin{equation}
c_i\left( \mathrm{sl}_3\right) _{\mid A_3\times \mathbf{\mu }_3}=c_j\left(
Sym^3E\right) _{\mid A_3\times \mathbf{\mu }_3}=0,\text{ }i,j=1,2,3,4.
\label{app1}
\end{equation}
Now recall (Section 3) that $\rho $ is a lift to $A_{PGL_3}^{*}$ of 
\[
\psi \left( \alpha c_3\left( W\right) \right) \in A_{A_3\ltimes T}^{*}\left(
Diag_{\mathrm{sl}_3}^{*}\right) 
\]
where 
\[
\psi :A_{A_3\ltimes T}^{*}\longrightarrow A_{A_3\ltimes T}^{*}\left( Diag_{%
\mathrm{sl}_3}^{*}\right) 
\]
is the (surjective) pullback. So, the image of $\rho $ under the restriction 
\[
A_{PGL_3}^{*}\longrightarrow A_{A_3\ltimes T}^{*}
\]
is of the form $\alpha c_3\left( W\right) +\xi $, for some $\xi \in \ker
\left( \psi \right) $.

Now, let us suppose $\rho $ is in the Chern subring $A_{Ch,PGL_3}^{*}$. By (%
\ref{app1}), we have 
\[
\alpha c_3\left( W\right) +\xi \in \ker \left( \varphi :A_{A_3\ltimes
T}^{*}\longrightarrow A_{A_3\times \mathbf{\mu }_3}^{*}\right) .
\]
From the commutative diagram 
\[
\begin{tabular}{ccc}
$A_{A_3\ltimes T}^{*}$ & $\stackrel{\varphi }{\longrightarrow }$ & $%
A_{A_3\times \mathbf{\mu }_3}^{*}$ \\ 
$^\psi \downarrow $ &  & $\downarrow ^{\psi ^{\prime }}$ \\ 
$A_{A_3\ltimes T}^{*}\left( Diag_{\mathrm{sl}_3}^{*}\right) $ & $\stackunder{%
\phi }{\longrightarrow }$ & $A_{A_3\times \mathbf{\mu }_3}^{*}\left( Diag_{%
\mathrm{sl}_3}^{*}\right) $%
\end{tabular}
\]
we get 
\[
\psi ^{\prime }\left( \alpha c_3\left( W\right) \right) =0.
\]
Therefore, if we show that $\alpha c_3\left( W\right) $ is not in the kernel
of $\psi ^{\prime }$, we will have proved that $\rho $ cannot be in the
Chern subring of $A_{PGL_3}^{*}$. To do this, let us consider the two
localization sequences\footnote{%
Here $Diag_{\mathrm{sl}_3}$ are the diagonal matrices in $\mathrm{sl}_3$ and
we identify $A_{A_3\times \mathbf{\mu }_3}^{*}\simeq A_{A_3}^{*}\otimes A_{%
\mathbf{\mu }_3}^{*}$ with 
\[
\frac{\mathbf{Z}\left[ \alpha \right] }{\left( 3\alpha \right) }\otimes 
\frac{\mathbf{Z}\left[ \beta \right] }{\left( 3\beta \right) }. 
\]
}: 
\begin{equation}
A_{A_3\times \mathbf{\mu }_3}^{*}\stackrel{\cdot (-\alpha ^2)\otimes 1}{%
\longrightarrow }A_{A_3\times \mathbf{\mu }_3}^{*}\left( Diag_{\mathrm{sl}%
_3}\right) \stackrel{p}{\longrightarrow }A_{A_3\times \mathbf{\mu }%
_3}^{*}\left( Diag_{\mathrm{sl}_3}\smallsetminus \left\{ 0\right\} \right)
\rightarrow 0  \label{app2}
\end{equation}
\begin{equation}
A_{\mathbf{\mu }_3}^{*}\simeq A_{A_3\times \mathbf{\mu }_3}^{*}\stackrel{}{%
\left( Z\right) \stackrel{j_{*}}{\longrightarrow }}A_{A_3\times \mathbf{\mu }%
_3}^{*}\left( Diag_{\mathrm{sl}_3}\smallsetminus \left\{ 0\right\} \right) 
\stackrel{q}{\longrightarrow }A_{A_3\times \mathbf{\mu }_3}^{*}\left( Diag_{%
\mathrm{sl}_3}^{*}\right) \rightarrow 0  \label{app3}
\end{equation}
where we used that 
\[
Z\simeq A_3\times \mathbf{C}^{*},
\]
$A_3\ltimes T$-equivariantly. Since (Section 3), 
\[
W\simeq \mathbf{C}_{\chi ,\mathbf{\mu }_3}\boxtimes \mathbf{C}_{perm,A_3}^3
\]
as $A_3\times \mathbf{\mu }_3$-representations (where $\mathbf{C}_{\chi ,%
\mathbf{\mu }_3}$ is the $\mathbf{\mu }_3$-representation of character $\chi
=\exp (i2\pi /3)$ and $\mathbf{C}_{perm,A_3}^3$ is the $A_3$-permutation
representation), its Chern roots are 
\[
\left\{ \beta +\alpha ,\beta -\alpha ,\beta \right\} 
\]
and then 
\begin{equation}
\alpha c_3(W)_{\mid A_3\times \mathbf{\mu }_3}=(\beta ^2-\alpha ^2)\alpha
\beta .  \label{app4}
\end{equation}
By (\ref{app4}) and (\ref{app2}), it is enough to prove that $j_{*}=0$.

Let us consider the pullback $E$ of $\mathbf{C}_{\chi ,\mathbf{\mu }_3}$ to $%
X=Diag_{sl_3}\smallsetminus \left\{ 0\right\} $ as an $A_3\times \mathbf{\mu 
}_3$-equivariant vector bundle, with $A_3$ acting trivially on $\mathbf{C}%
_{\chi ,\mathbf{\mu }_3}$ and $A_3\times \mathbf{\mu }_3$ acting as usual on 
$X$ (i.e. $\mathbf{\mu }_3$ acting trivially and $A_3$ by permutations). We
have 
\[
j^{*}\left( c_1\left( E\right) \right) \equiv c_1(E)_{\mid \mu _3}=\beta . 
\]
But we also have $j_{*}\left( 1\right) =0,$ since 
\[
Z=D^{-1}\left( \left\{ 0\right\} \right) 
\]
where 
\begin{eqnarray}
D &:&X\rightarrow \mathbf{A}^1  \nonumber \\
\left( \lambda _1,\lambda _2,\lambda _3\right) &\longmapsto &\left( \lambda
_1-\lambda _2\right) \left( \lambda _1-\lambda _3\right) \left( \lambda
_2-\lambda _3\right)
\end{eqnarray}
is the square root of the discriminant (which is $A_3\times \mathbf{\mu }_3$%
-equivariant). So $j_{*}=0$ and we conclude. \TeXButton{End Proof}{\endproof}%
\ 

%%%%%%%%%%%%%%%%%%%%%%%%% End PGL_36.tex %%%%%%%%%%%%%%%%%%%%%%%%%%%%%%
}

\QSubDoc{Include PGL_3bib}{%%%%%%%%%%%%%%%%%%%%%%%%% Start PGL_3bib.tex %%%%%%%%%%%%%%%%%%%%%%%%%%%%%

\LaTeXparent{PGL_3}

%%%%%%%%%%%%%%%%%%%%%%%%%% End PGL_3bib.tex %%%%%%%%%%%%%%%%%%%%%%%%%%%%%%
}

\end{document}